\documentclass{article}

\usepackage[hidelinks]{hyperref}
\usepackage{xcolor}
\usepackage{xspace}
\usepackage{amsmath}
\usepackage{amsthm}
\usepackage{amssymb}
\usepackage{booktabs}

\usepackage[
    style=trad-abbrv,   
    backend=biber,      
    giveninits=true,    
    url=true,           
    doi=true            
    ]{biblatex}
\usepackage{subcaption}
\usepackage{dashbox}
\usepackage{geometry}

\usepackage{dashbox}

\usepackage{algorithm}
\usepackage{algpseudocode}

\newcommand{\algrule}[1][.2pt]{\par\vskip.5\baselineskip\hrule height #1\par\vskip.5\baselineskip}

\newcommand{\ignore}[1]{}
\newcommand{\NaN}{\texttt{NaN}\xspace}
\newcommand{\NaNs}{\texttt{NaN}s\xspace}
\newcommand{\Inf}{\texttt{Inf}\xspace}
\newcommand{\Infs}{\texttt{Inf}s\xspace}
\newcommand{\EV}{\texttt{ExVal}\xspace}

\newcommand{\Comments}{1}
\ifnum\Comments=1
    \usepackage[disable]{todonotes}
\else
    \usepackage{todonotes}
\fi

\usepackage{algorithm}
\usepackage{algpseudocode}

\ifnum\Comments=1
    \usepackage[disable]{todonotes}
\else
    \usepackage{todonotes}
\fi

\definecolor{darkgreen}{rgb}{0,0.5,0}
\definecolor{purple}{rgb}{1,0,1}

\newcommand{\kibitz}[2]{\ifnum\Comments=1\textcolor{#1}{#2}\fi}
\newcommand{\jackson}[1]{\kibitz{blue}      {[Jackson: #1]}}
\newcommand{\jwd}[1]{\kibitz{orange}      {[JIM D: #1]}}

\definecolor{darkperiwinkle}{RGB}{102, 95, 209}
\newcommand{\ejr}[1]{\kibitz{darkperiwinkle}{ejr: #1}}

\newcommand{\makeblasname}[2]{\expandafter\newcommand\csname #1\endcsname{\ensuremath{\operatorname{\mathsf{#2}}}\xspace}}
\makeblasname{gemm}{GEMM}
\makeblasname{trsm}{TRSM}
\makeblasname{trmm}{TRMM}
\makeblasname{gemmtr}{GEMMTR}
\makeblasname{syrk}{SYRK}
\makeblasname{syrtwok}{SYR2K}
\makeblasname{axpy}{AXPY}
\makeblasname{blasdot}{DOT}

\makeblasname{uplo}{UPLO}
\makeblasname{trans}{TRANS}
\makeblasname{larfb}{LARFB}
\makeblasname{larzb}{LARZB}

\newcommand{\opeval}[1]{[ #1 ]}
\newcommand{\lopeval}[1]{\left[ #1 \right]}
\newcommand{\Xhat}{\hat{X}}
\newcommand{\xhat}{\hat{x}}
\newcommand{\Chat}{\hat{C}}

\newcommand{\dogemm}[2]{\gemm\opeval{\ensuremath{#1 \cdot #2}}}
\newcommand{\ldogemm}[2]{\gemm\lopeval{\ensuremath{#1 \cdot #2}}}
\newcommand{\dotrsm}[2]{\trsm\opeval{\ensuremath{#1^{-1} #2}}}
\newcommand{\dotrsmT}[2]{\trsm\opeval{\ensuremath{#1^{-T} #2}}}
\newcommand{\ldotrsm}[2]{\trsm\lopeval{\ensuremath{#1^{-1} #2}}}
\newcommand{\dosyrk}[2]{\syrk\opeval{#1 -  #2 \cdot #2^T}}

\newcommand{\rowmax}{\xspace\ensuremath{\operatorname{\mathsf{rowmax}}}\xspace}
\newcommand{\colmax}{\ensuremath{\operatorname{\mathsf{colmax}}}\xspace}
\newcommand{\diag}{\operatorname{\mathsf{diag}}}

\NewDocumentCommand{\matentry}{m O{} m m}{\ensuremath{{#2{#1}}(#3,#4)}}
\NewDocumentCommand{\matblock}{m O{} m m}{\ensuremath{{#2{#1}}_{#3,#4}}}
\NewDocumentCommand{\vecentry}{m O{} m}{\ensuremath{{#2{#1}}(#3)}}
\NewDocumentCommand{\vecaccess}{m O{} m}{\ensuremath{{#2{#1}}_{#3}}}
\newcommand{\midx}[3]{\matentry{#1}{#2}{#3}}

\newcommand{\dimall}{:}
\newcommand{\dimrange}[2]{{#1}{:}{#2}}

\newcommand{\idxset}[1]{\mathbb{#1}}
\newcommand{\iI}{\idxset{I}}

\newcommand{\failcase}[1]{\textcolor{red}{\textbf{#1}}}
\newcommand{\mfailcase}[1]{\textcolor{red}{\ensuremath{\mathbf{#1}}}}

\title{How to grade the accuracy of the BLAS}

\author{James Demmel\footnote{University of California, Berkeley; demmel@berkeley.edu}
\and Greg Henry\footnote{NVIDIA Corp.; grhenry@nvidia.com}
\and Igor Kozachenko\footnote{University of California, Berkeley; igork@berkeley.edu}
\and Julien Langou\footnote{University of Colorado Denver; julien.langou@ucdenver.edu}
\and Xiaoye Sherry Li\footnote{Lawrence Berkeley National Lab; xsli@lbl.gov}
\and Jason Riedy\footnote{Microsoft Corp.; jason@acm.org}
\and Jackson Vanover\footnote{Unaffiliated; jacksonvanover@gmail.com}}

\date{\today}

\begin{document}

\maketitle

\begin{abstract}
Motivated by accelerating machine learning (ML),
many computer vendors and chip manufacturers are
building accelerators for matrix multiplication,
which save time and energy by operating in the
lower precisions needed for ML. This has in turn
motivated many efforts to use these accelerators
to provide faster matrix multiplication
implementations with the higher
precision required by many other linear algebra
applications. Motivated by the large design space
of algorithms for approximating higher precision,
with significant performance/accuracy tradeoffs,
we provide a benchmark to ``grade the accuracy''
of a matrix multiplication implementation (or the
BLAS more generally), ranging from an ``A'' for
attaining the classic floating point error bound,
to a ``C'' for attaining
a weaker but still useful bound, that is satisfied
by Strassen-like algorithms. We also propose
``ungameable'' tests that vendors or users can run to verify
their promised accuracy, and describe how these
different grades impact the accuracy of applications
like LU, QR, and Cholesky decomposition.
Our test code is publicly available
at \url{github.com/Reference-LAPACK/grading-the-BLAS} for developers
and users, and we also plan to publicly release all our test results.
\end{abstract}

\section{Introduction}
Many computer vendors and chip manufacturers are building
low-precision matrix multiplication accelerators in order to
accelerate machine learning applications
\cite{TPU17,10.1145/3773656.3773670,Hopper,Rubin}.
Given the large and growing gap in speed between these
low-precision implementations compared to traditional
higher-precision SGEMM and DGEMM implementations
\cite{blasnetlib}, there has been much recent effort (and a long
history) on how to use these
low-precision accelerators to provide faster implementations
of SGEMM and DGEMM as well
\cite{Ozaki_Ogita_Oishi_Rump_2012,Henry_Tang_Heinecke_2019,Ootomo_Ozaki_Yokota_2024,Uchino_Ozaki_Imamura_2024,Ootomo_Yokota_2023,10.1145/3773656.3773670,NvidiaPASC25,PASC25,Zhang_2025}.

Here is a partial taxonomy of different matrix multiply algorithms that we consider.
(We consider only real
and not complex matrices, 
and leave the more complicated design
space of complex arithmetic to future work.) 
First, we distinguish classical $O(n^3)$ algorithms from Strassen-like algorithms.
By $O(n^3)$, we mean any algorithm that computes each entry $\matentry{C}{i}{j}$ of $C = A \cdot B$
as a dot product of row $i$ of $A$ with column $j$ of $B$, independent of the
dimensions of $A$ and $B$. Strassen-like algorithms include not just Strassen's
original $O(n^{\log_2 7})$ algorithm \cite{Strassen69}, but the many variants
that have been invented since then~\cite{Coppersmith-Winograd1990,Gall2014,almanduanwilliamsxuxuzhou25,FalconGEMM}. 
Based on \cite{Miller1972},
we know that all Strassen-like algorithms can only satisfy a weaker error bound
than $O(n^3)$ algorithms \cite{Bini_Lotti1980,DemmelDumitriuHoltzKleinberg2007,BallardBensonDruinskyLipshitzSchwartz2016}, details below.

Second, we independently distinguish between different kinds of arithmetic that
may be used. The possibilities include
\begin{enumerate}
    \item Classic floating point, with $A$, $B$, $C$ and all intermediate results using the same floating point format.
    \item Classic floating point, but with $A$, $B$, $C$ and intermediate results using
    possibly different precisions, \textit{e.g.,} higher precision for intermediate results, and
    possibly also $C$, than is used for $A$ and $B$.
    \item Representing entries of $A$ and $B$ as sums of lower-precision floating
    point numbers, each one storing subsets of the  bits of $A$ and $B$.
    These sums may or may not represent each entry of $A$ and $B$ exactly
    \cite{Henry_Tang_Heinecke_2019}.
    \item Representing entries of $A$ and $B$ as sums of scaled short integers, each
    one again storing subsets of the mantissa bits of $A$ and $B$. We cannot store
    a separate exponent for each short integer, rather groups of integers
    (\textit{e.g.,} rows of $A$ or columns of $B$) must share the same exponent, so
    smaller numbers in a group may be stored with lower precision \cite{Uchino_Ozaki_Imamura_2024,Ootomo_Ozaki_Yokota_2024}. In our later numerical
    tests, we refer to the implementation of
    this that we test as the ``Ozaki I''
    algorithm.
    \item 
    Similar to the last approach, floating point numbers can be represented as a sum of scaled integers, but an integer $j$ can
    be represented as a vector of values
    $j\mod m_i$ for a set of pairwise
    relatively prime integers $m_i$, with arithmetic
    done based on the Chinese Remainder Theorem 
    \cite{10.1145/3773656.3773670}.
    In our later numerical tests, we refer to
    the implementation of this that we test
    as the ``Ozaki II'' algorithm.
\end{enumerate}

For brevity, we later refer to 3., 4. and 5. above as
``emulation'' of conventional floating point.

It is also possible, and sometimes likely, to
hybridize the above algorithms. For example,
any performant Strassen-like algorithm will recursively
multiply smaller matrices, until the matrices are small
enough to switch to a faster, $O(n^3)$ algorithm.
And each recursive call could potentially use a
different Strassen-like algorithm \cite{BallardBensonDruinskyLipshitzSchwartz2016, BensonBallardSPPPP15}.
More generally, auto-tuning techniques may customize algorithms for different
ranges of matrix dimensions \cite{FalconGEMM}.
Also, NVIDIA's cuBLAS library~\cite{10.1145/3773656.3773670} 
performs some tests on the input matrices $A$ and $B$
to quickly decide if their emulation approach will
satisfy the strongest error bound we describe below,
and falls back to using conventional floating point
if they cannot guarantee that emulation will be
accurate enough.
See Section~\ref{sec:related-work} for more details.

It is also possible to prescale $A$ and $B$
by diagonal matrices to get $A' = D_1 \cdot A \cdot D_3$ and $B' = D_3^{-1} \cdot B \cdot D_2$,
and postscale their product $C' = A' \cdot B'$
to get the desired $C = D_1^{-1} \cdot C' \cdot D_2^{-1}$.
One can choose the diagonal matrices $D_i$ to significantly improve the error, and error bound, for Strassen-like algorithms \cite{BallardBensonDruinskyLipshitzSchwartz2016}
when $A$ and $B$ have entries with greatly varying magnitudes, for an additional cost
of $O(n^2)$. The same scaling could also improve the
accuracy of approaches using emulation, for
example by narrowing the range of exponents needed
to represent the matrix entries.

Given the many ways to approximate the higher-precision matrix
multiplication (and BLAS more generally) desired by users, and
the strong performance/accuracy tradeoffs among the methods
described above, we think it is important
for vendors to carefully document the accuracy they promise to
provide, and for both vendors and users to be able to verify
accuracy claims using publicly available test code.
Our goal is to provide this test code. In particular,
this test code should not be ``gameable'', \textit{i.e.,} it should not be possible to pass
the tests with a high-speed implementation, while actually not
providing the claimed accuracy in all cases 
(we discuss what gaming means in
more detail in Section~\ref{sec:gaming}).
This is important
since vendor implementations are likely to be proprietary, and
not publicly available (unlike our proposed test code), so that an actual public error analysis
of the source code may not be possible.

\begin{table}[ht]
  \centering
  \caption{Summary of notation and evaluation grades for \gemm and \trsm. See Section~\ref{sec:grading} for an explanation of notation.
}
  \label{tab:notation-grades}
  \begin{tabular}{|c|c|c|} \hline
    Expression  & $C = A \cdot B$ & $ T \cdot X = B$ \\
    FP Evaluation  & $\Chat = \dogemm{A}{B}$ & $\Xhat = \dotrsm{T}{B}$ \\
    Error Expr.            & $E = \Chat - A \cdot B$ & $R = B - T \cdot \Xhat$ \\ \hline 
    ``A'' & $|E| \leq f_A(n) \epsilon |A| \cdot |B|$ & $|R| \leq g_A(n) \epsilon |T| \cdot |\Xhat|$ \\
    ``B'' & $|E| \leq f_B(n) \epsilon \rowmax(|A|) \colmax(|B|)$ & $|R| \leq g_B(n) \epsilon \rowmax(|T|) \colmax(|\Xhat|)$ \\
    ``C'' & $\|E\| \leq f_C(n) \epsilon \|A\| \|B\|$ & $\|R\| \leq g_C(n) \epsilon \|T\| \|\Xhat\|$ \\
    \hline
    
  \end{tabular}
\end{table}

We initially intended to design our tests to ``reverse engineer'' the
algorithm being used, allowing us to know which error analysis applies.
But the large and growing algorithmic design space described above, the
discovery of ever more sophisticated ways to game our proposed
tests\footnote{Including one suggested by
an unknown audience member during a presentation.},
and the inability to predict what new algorithm or game an LLM might imagine
(or hallucinate) tomorrow~\cite{AlphaTensor,chatBLAS2024}, 
led us to change our goal
(see Section~\ref{sec:gaming} for more
discussion). Our goal now is
to reliably assign one of three ``grades'' to an implementation, depending
on what ``worst case'' error bound it satisfies, ranging from an ``A'' for the
error bound satisfied by classical $O(n^3)$ matrix multiplication
using floating point, to a ``C'' for the bound satisfied by
Strassen-like algorithms (without the diagonal scaling described above).
It turns out that our initial ``reverse engineering'' tests, and their
updated versions designed to defend against gaming, will serve our
purpose of assigning grades reliably, with a caveat: given the
possibility of using different algorithms for different matrix
dimensions and different hardware configurations (\textit{e.g.,} number of
processors, available accelerators, \textit{etc}.),
we can only expect (but not prove;
see Section~\ref{sec:gaming}) our grades apply to the matrix sizes
and hardware configurations on which our tests are run.

These ``A'' and ``C'' grades (error bounds) have been widely used in the
numerical analysis literature to analyze linear algebra algorithms,
but not the ``B'' grade. We summarize what is known about the use
of these error bounds, and what problems remain.

\ignore{
We can
do this in many cases, but at least one gaming trick\footnote{Suggested by
an unknown audience member during a presentation.} potentially makes two
quite different implementations indistinguishable, as we describe below.
Nonetheless, as we will see, we can reliably assign the same grade
(a ``B'') to both implementations.
}

\ignore{
It may also be the case that a vendor uses a different algorithm depending
on the matrix size. For example, it is common for a Strassen-like algorithm
to revert to using an $O(n^3)$ algorithm once matrices are small enough.
More generally, auto-tuning techniques may customize algorithms for different
ranges of matrix dimensions. This means that grades may depend on matrix
dimensions. While our tests described below can ``reverse engineer'' some of
these algorithmic changes, we believe users will be most interested in the
``lowest grade'', which describes what they can rely on most generally.
}

The rest of this paper is organized as follows.
Section~\ref{sec:gaming} describes
in more detail what we mean by
``gaming'', and the assumptions
underlying our grading tests.
Section~\ref{sec:grading} describes our \gemm grading scheme, \textit{i.e.,} the
error bounds corresponding to ``A'', ``B'', and ``C'' grades.
(An ``A+'' or even ``A++'' is also possible, but we will not test for these.)
Section~\ref{sec:testing} describes our \gemm grading tests that initially attempted to reverse engineer
the algorithm by a careful choice of ``ungameable'' test matrices, and
how they are used to assign grades.
Section~\ref{sec:testresults} describes the results of our \gemm tests on a variety
of matrix multiplication implementations.  
Section~\ref{sec:trsm-residual-gemm} describes how to
test \trsm.
Section~\ref{sec:otherBLAS} describes how these tests can be extended to other BLAS routines.
Section~\ref{sec:LUQRChol} describes how using BLAS with different grades
impacts the error analyses of the LU, QR, and Cholesky factorizations.
Section~\ref{sec:Exceptions} discusses grading for
consistent exception handling.
Section~\ref{sec:grading_fairly} dispels some common myths about emulation in terms of reliability and why our approach should handle reliability concerns.
Section~\ref{sec:related-work} discusses
related work and on hybrid approaches to get better accuracy. Section~\ref{sec:Conclusions} summarizes our conclusions, and lists open
problems and future work.

\section{What does ``gaming'' mean?}\label{sec:gaming}


``Gaming'' means getting a higher grade
than deserved, by intentionally or
unintentionally using an algorithm that
gets the higher grade on our
difficult test cases, but not all
difficult cases.
We describe
in more detail the reasons we need to
be concerned about this, why we cannot
provably prevent it, why we believe
our tests will work (with high probability),
and why improving our tests over time is
a likely scenario.

To provide some motivation, a software company recently hosted an open
competition, with a leader board, for
anyone to compete in providing the
fastest implementation of batched
QR decomposition on a GPU. Batched QR
decomposition means that one is 
given many matrices simultaneously,
so that one can parallelize across
matrices as well as for a single matrix. There were 15K submissions,
with many beating the best available
implementation from Nvidia. To get
onto the leader board, one had to
pass tests provided by the host,
which were designed to rule out
certain less accurate implementations
(eg Cholesky QR). The implementors
of the fastest submissions were 
invited to an event at the company.
One declined, saying he had neither
read the specification of the problem
nor the code produced by the LLM he
used before submitting it. The next
day several collaborators reported
that they had tried some of the
leading implementations on
problem sizes not in the provided
test code, and they had simply
crashed, because they had
been designed only to work on the
problem sizes in the test code.

The point of this example is to 
emphasize that we should no longer
assume that numerical (or other!) libraries are written and tested by trustworthy teams of professional programmers. 

Ideally, we would like to provide a
proof that a code is correct. But
without access to source or object code, which we are assuming is not
available because it may be proprietary, we cannot provide a
proof. Exhaustive testing is 
obviously impossible because of 
the number of possible floating
point matrix inputs. Even
testing on one input does not
mean it will be correct if run again,
because \linebreak 
(1) if the underlying available
hardware changes dynamically, a different
algorithm could be used, (2) an algorithm could potentially 
``flip a coin'' to decide which algorithm to use on each input, or (3) it could
even randomly decide (with very low probability)
to deliberately get the wrong answer.

Although we do not expect such ``malware'' to
be written, the possibility motivates us
to be clear about what we do assume about
the software provider (AI or human) in order
to design our tests:
We will assume that an ``honest'' provider's motivation
is to maximize performance subject to 
attaining a certain minimum accuracy, i.e.
``grade'' (including ``no accuracy'', 
if speed on
certain cases is all that matters).
Sometimes this will 
mean slowing down on ``hard'' cases,
i.e. different algorithms may be used
on inputs of different difficulties
(see Section~\ref{sec:related-work}, 
\cite{10.1145/3773656.3773670}).
In contrast, a ``dishonest'' (or
``honest but mistaken'') developer
may intentionally (or unintentionally)
just detect
the hard examples in our test code, or a
not-too-large superset, and just run a
slower and more accurate algorithm on these, thereby getting the
desired (incorrect) higher grade while continuing to
be fast on most inputs. 

This possibility (intentional or
unintentional) of detecting the hard
examples in our test code (and running
a slower, more accurate algorithm), 
but not detecting all hard examples, 
is an example of ``gaming''.
Our goal is to make detecting just our
test cases (or a not-too-large superset)
expensive enough to thwart gaming. 
As stated above, 
we do not have a proof of this,
just the intuition that any reasonably 
cheap (say $O(n^2)$) detection algorithm
would be triggered on most, if not all, 
hard cases, and many others (\cite{10.1145/3773656.3773670}
uses an $O(n^3)$ detection algorithm).
So a ``dishonest'' developer would not be
able to claim both their desired higher
grade and (average) good performance,
and an ``honest but mistaken'' developer
would be warned about a lower grade.

As described later, one of our tests
(Test 1c) was motivated by a suggested 
way to game its predecessor (Test 1b) 
proposed by an audience member during a
presentation. Just as the computer security
community does research on both new attacks
and new defenses, we invite the numerical community to consider new ways to ``game'' 
our tests and share their results, to help
us provide even less ``gameable'' tests
in the future.

\section{How to grade matrix multiplication}\label{sec:grading}

We now describe the different ``grades'' (error bounds)
that matrix multiplication can achieve.
For simplicity, we initially make the following assumptions:
(1) We are multiplying $n \times n$ matrices $C = A \cdot B$.
(2) Machine precision $\epsilon$, the overflow threshold $OV$, and the underflow
threshold $UN$ refer to the precision in which the output matrix $C$ is stored
($A$ and $B$ could be stored in the same
or lower precision).
(3) Our error bounds will assume overflow does not occur, and that if underflow
occurs, it is small compared to the error bounds below. (Dealing with exception
handling is future work.)

We note that the use of mixed precision, \textit{e.g.,} using higher
precision before rounding the output $C$ to a lower precision
defined by $\epsilon$, may satisfy different (\textit{e.g.,} stronger)
error bounds than those below. We also leave the formalization
and testing of this larger set of possible error bounds
to future work.

We use \gemm and \trsm for generic operations, and typed names such as DGEMM and DDOT for concrete BLAS interfaces.
We use the notation $\dogemm{A}{B}$ to refer to the result computed by the
algorithm used for $A \cdot B$,
$\| A \|$ is the maximum-entry norm or matrix $A$,
and $f(n)$ is an (at most polynomially) increasing function of the matrix dimension $n$ (details later).
Let $\rowmax(|M|)$ be the \emph{column} vector where entry $i$ is the maximum magnitude entry of row $i$ in $M$,
$\max_j |\matentry{M}{i}{j}|$.
Let $\colmax(|M|)$ be the \emph{row} vector where entry $j$ is the maximum magnitude entry of column $j$ in $M$,
$\max_i |\matentry{M}{i}{j}|$.

The different grades that a matrix multiply algorithm can achieve are as follows. 

\begin{description}
    \item[Grade = ``A'':] $| \dogemm{A}{B} - A \cdot B | \leq f_A(n) \, \epsilon \, |A| \cdot |B|$.
This is the strongest bound, and the test currently used by the BLAS test code\footnote{Albeit with the weak
requirement of $f_A(n) = \epsilon^{-1/2}$.}.
According to~\cite{Miller1972}, any algorithm satisfying this bound must do $O(n^3)$ operations,
and so it rules out Strassen-like algorithms. It may also rule out the integer arithmetic
approaches in~\cite{Ootomo_Ozaki_Yokota_2024,uchino2025high},
depending on range requirements.
\item[Grade = ``C'':] $\| \dogemm{A}{B} - (A \cdot B) \| \leq f_C(n) \, \epsilon \, \| A \| \cdot \| B \|$.
This is the weakest bound, and satisfiable by Strassen-like algorithms as described in~\cite{Bini_Lotti1980,
DemmelDumitriuHoltzKleinberg2007}.
This bound suffices for norm-wise backward error analysis of LU, QR and
other LAPACK routines~\cite{DemmelHigham1992}, but will not preserve some properties that some
users might expect (\textit{e.g.,} column scaling invariance of LU and QR). 
\item[Grade = ``B'':]
$| \dogemm{A}{B} - A \cdot B | \leq f_B(n) \, \epsilon \, \rowmax(|A|) \cdot \colmax(|B|)$.
This grade is “in between” an ``A'' and a ``C''. Like getting an ``A'', but unlike a ``C'',
it is invariant under row scaling of $A$ (changing $A$ to $D_1 \cdot A$ for $D_1$
diagonal) and column scaling of $B$ (changing $B$ to $B \cdot D_2$).
Unlike getting an ``A'', it is not invariant
under simultaneous scaling of columns of $A$ and rows of $B$
(changing $A$ to $A \cdot D_3$ and $B$ to $D_3^{-1} \cdot B$). We use this fact in the design of Test~2 below.
\end{description}

We see that getting an ``A'' is strictly stronger than getting a ``B'', which is
strictly stronger than getting a ``C'', which will guide our test design in the
next section.

It is possible to get a higher grade than ``A'', call it ``A+'', by using the
$O(n^3)$ algorithm with higher precision
to compute all the dot products;
we will use such an implementation
to perform our tests later in
Section~\ref{sec:testresults}.
It is also possible to get an ``A++''
by computing
all entries of $C$ to high relative accuracy, even the tiniest ones, and
perhaps even correctly rounded. There is a literature of algorithms
that do this \cite{Rump2009,DemmelHida2003,Ozaki_Ogita_Oishi_Rump_2012,xblas:toms,Ogita_Rump_Oishi_2005,Rump_Ogita_Oishi_1_2008,Rump_Ogita_Oishi_2_2008}, but because of their
higher cost, we will not consider these
further.
We note that we can also correctly
assign grades
using a ``reference''
implementation that uses the $O(n^3)$ algorithm with
classic floating point, since this leads to an uncertainty
in $f_A(n)$ of at most $\pm n$, or $\pm \log_2 n$ if sums are
computed using a binary tree.


We note that by choosing $D_1$ so that the rows of $D_1\cdot A$ have unit norms,
and $D_2$ so that the columns of $B\cdot D_2$ have unit norms, performing Strassen to multiply
$\Chat = (D_1\cdot A)\cdot (B\cdot D_2)$, and then forming $C = D_1^{-1}\cdot \Chat \cdot D_2^{-1}$, we can get Strassen to earn a ``B'' instead of a ``C'' for an extra $O(n^2)$ cost~\cite{dumitrescu1998accuracy}.
Zero rows of $A$ and columns of $B$ must be
handled separately,  by zeroing out the corresponding rows and columns of C; this
leads to a gaming possibility described in the next section.

$D_3$ can be chosen so that for each $i$,
column $i$ of $A$ and row $i$ of $B$ have the same norm
(let $\matentry{D_3}{i}{i}=\sqrt{\|\matentry{B}{i}{\dimall}\|/\|\matentry{A}{\dimall}{i}\|}$, or
set both $\matentry{A}{\dimall}{i}$ and $\matentry{B}{i}{\dimall}$ to zero if either
one is already zero). This provably results in
$\|A \cdot D_3\| \cdot \|D_3^{-1} \cdot B\| \leq \|A\| \cdot \|B\|$, lowering the attainable error bound for grade ``C'' \cite{BallardBensonDruinskyLipshitzSchwartz2016}.
One can also apply ``outer scaling'' (by $D_1$ and $D_2$)
and ``inner scaling'' (by $D_3$) multiple times in an
alternating fashion to get the best of both worlds,
see \cite[Sec. 6]{BallardBensonDruinskyLipshitzSchwartz2016}
for an analysis. Finally, we note that the diagonal
entries of $D_i$ may be limited to powers of 2, to
avoid roundoff.

Note that the error bound for grade ``A'' does not
change when doing inner or outer scaling (ignoring
overflow and underflow), unlike grades ``B'' and ``C''.
This motivates us to have some test cases that let us
detect whether the algorithm uses inner or outer scaling,
since this would confirm that the algorithm is
{\em not} using an $O(n^3)$ algorithm with classic
floating point.

Given a grade ``A'', ``B'', or ``C'', we
also want to estimate $f_A(n)$, $f_B(n)$
and $f_C(n)$, which of course can depend
on more dimensions in the rectangular case.
Upper bounds on $f_X(n)$ have been
published for many algorithms
\cite{Higham2002,DemmelDumitriuHoltzKleinberg2007,Uchino_Ozaki_Imamura_2024}
but to our knowledge
only some are attainable in practice, 
with specially constructed 
$A$ and $B$,
and with assumptions like order of summation.
In this paper we will report on 
measured values of $f_X(n)$ using test
matrices described in the next section,
and leave using examples with even
larger $f_X(n)$ to future work. 


\section{Proposed test cases}\label{sec:testing}

We describe tests that assume accurate (grade ``A'')
versions of $A \cdot B$  and $|A| \cdot |B|$ are available, either from a reference version of \gemm as described above,
or an ``A+'' (or ``A++'') implementation. We describe a sequence of
tests, from simple to more complicated, that were motivated
by possible ways to game them. All these tests are still
useful, because they only provide ``true positives'', \textit{i.e.,}
examples where an algorithm fails to satisfy a particular
error bound, \textit{i.e.,} get a particular grade, and so must get a lower grade. Gaming can only cause false negatives.

The first sequence of tests (Test\_1a through Test\_1c) were originally intended to identify Strassen-like
algorithms. As stated
above, they can provide ``true positives'' for
Strassen. But as more sophisticated gaming techniques were considered, and countermeasures proposed, we
realized that they might only be able to distinguish
$O(n^3)$ algorithms using classic floating point from
alternatives, either Strassen-like or using emulation.
So these test are still valuable in distinguishing
which algorithms should {\em not} get an ``A''.

\begin{description}
\item[Test\_1a:] The simplest way to identify Strassen-like routines is to choose $A = randn(n)$
and $B=randn(n)$, and set a few randomly chosen rows of $A$ and columns of $B$ to zero.
Any $O(n^3)$ algorithm (with any arithmetic) will compute $A\cdot B$ with the same
rows and columns of $C$ equal to zero, but Strassen will not with high
probability (w.h.p), because of rounding errors during cancellation.

However, this very simple test can be gamed using the idea at the end of
Section~\ref{sec:grading}. If an implementation of Strassen tries to scale
the rows of $A$ and columns of $B$ to have unit norms to try to earn a ``B'',
it will detect zero rows and columns, which can be zeroed out in the product $C$,
at $O(n^2)$ cost, and so Strassen could pass Test\_1a. We try to address this weakness with Test\_1b.

\item[Test\_1b:]
Start with $A = randn(n)$ and $B = randn(n)$. Pick a random subset of $1<k<n$ rows of $A$,
labeled $(i_1, \ldots, i_k)$, and $k$ random columns of $B$, labeled $(j_1, \ldots, j_k)$. For each $m\in [1:k]$, pick a
random proper subset $S_m$ of $[1:n]$, and zero out entries $\matentry{A}{i_m}{r}$ for all $r\in S_m$,
and zero out entries $\matentry{B}{s}{j_m}$ for all $s$ not in $S_m$. This means that the dot
product $\matentry{A}{i_m}{\dimall}\cdot \matentry{B}{\dimall}{j_m} = \matentry{(A\cdot B)}{i_m}{j_m}$ is exactly zero. If each $\matentry{(A\cdot B)}{i_m}{j_m}$
is computed as a conventional dot product, it will be a sum of zeros and so exactly zero,
in any precision. But with any Strassen-like algorithm, getting zero will require
exact cancellation, which is unlikely. And getting
exactly $k$ zeros is even more unlikely.
Computing any nonzero value violates the requirement for an ``A''.

Obviously other sparsity patterns for $A$ and $B$ might work too, but using “minimal”
sparsity as in Test\_1b would seem to maximize the chance of recognizing a Strassen-like
algorithm.
Since performant implementations of Strassen-like algorithms will not recurse down to
$n=1$, but call an $O(n^3)$ algorithm instead on small enough matrices, one could try to 
use Test\_1b 
to “reverse engineer” the size $n_0$ of matrices at or below which an $O(n^3)$
algorithm is used, by finding the largest matrix that still passes Test\_1b.
But of course this might depend on how the
Strassen implementation deals with
non-square, non-powers-of-two dimensions.

A failure to compute the expected $k$ zero entries is a clear signal that a
Strassen-like algorithm is being used. In other words, there are no false positive
signals for Strassen.
However, there is another way for a Strassen-like algorithm to ``game'' Test\_1b,
by setting computed entries of $C$ less than an error bound to zero\footnotemark[1].
Since our test code will be public, a vendor could read our test code and
precisely decide what this threshold would be. By first scaling the rows of $A$ and
columns of $B$ to have unit norms (as described in 
Section~\ref{sec:grading}), and setting small values of
$\Chat = (D_1 \cdot A) \cdot (B \cdot D_2)$ to zero this would make the behavior
hard to distinguish from an $O(n^3)$ algorithm using
emulation, \textit{e.g.,}
\cite{Uchino_Ozaki_Imamura_2024}, which earns a ``B''.

\ignore{
More tests to identify $O(n^3)$ and usual FP:
\begin{enumerate*}
\item A = rand(n,n), B = rand(n,n), multiply some rows of A, cols of B
by a very large constant x (x**2 does not overflow), so each entry of C =A*B is accurate using $O(n^3)$ algorithm, even with int8.
Strassen will get very large errors in some entries.
\item $\sqrt(OV)^2$ in some locations, $\sqrt(UN)^2$ in some, to test dynamic range,
eliminate use of short exponent approximations.
\item (0*y + y*0 + (x-1)*(x+1) - x*x) scaled so shorter precision, limited exponent
range would be wrong.
\item make sure \Inf propagates just to a whole row or whole column, to detect Strassen.
\item cases where intermediate terms in Strassen would overflow, but not $O(n^3)$
\end{enumerate*}
} 

\item[Test\_1c:] We strengthen Test\_1b so that only an $O(n^3)$ algorithm using
conventional floating point can pass it. We take the test matrices $A$ and $B$
proposed for Test\_1b, and in some subset of the $k$ pairs of
rows $\matentry{A}{i_m}{\dimall}$ and columns $\matentry{B}{\dimall}{j_m}$, we choose a
random $k_m$ and
set $\matentry{A}{i_m}{k_m}=\matentry{B}{k_m}{j_m}=x$, so
that the correct value of $\matentry{C}{i_m}{j_m}=x^2$.
If we choose the nonzero entries in row $i_m$ of $A$
and column $j_m$ of $B$ that are multiplied by zeros
to be close to $\sqrt{OV}$, $x$ to be close to
$\sqrt{UN}$, and the other entries of $A$ and $B$
to be $O(1)$, then scaling $D_1 \cdot A$
(resp. $B \cdot D_2$) to have unit row (resp. column)
norms will make the $x^2$ values underflow to zero, causing an incorrect value $\matentry{C}{i_m}{j_m}=0$.
Scaling by $D_3$ will not change the matrices (much)
since all the column norms of $A$ and row norms of $B$
are close to $\sqrt{OV}$.
A Strassen-like algorithm will not compute $x^2$
accurately w.h.p., whether or not it sets tiny entries
of $C$ to zero.
And the use of emulation will set the entries near
$\sqrt{UN}$ to zero, so that an $O(n^3)$ algorithm using
emulation will also compute 0 instead of $x^2$.
In other words, Test\_1c should reliably distinguish
an $O(n^3)$ algorithm using classic floating point,
which deserves an ``A'' grade, from the alternatives.
We return to distinguishing ``B'' from ``C'' grades
below.
\end{description}

Having even more ``ungameable'' tests will make our test
code more reliable, in case someone invents a game
we haven't thought of yet. So here is a second
sequence of tests, again designed to be progressively
less gameable.

\begin{description}
\item[Test\_2a:]
Let $x = rand(1,n) + 1$ be a random row vector with entries in the range [1,2].
Let $D$ be a diagonal matrix with entries equal to powers of two, ranging
from close to overflow to close to underflow, with roughly evenly spaced
exponents. Let $y = x\cdot D$ and $z^T = D^{-1}\cdot x^T$, so that $y\cdot z^T = x\cdot x^T = O(n)$ would be
computed to high relative accuracy using standard floating point, but an $O(n^3)$ algorithm with emulation
would result in a high relative error unless the entire exponent
range is represented. Now let each row of $A$ equal $y$ and each column of
$B$ equal $z^T$, and test that each entry is accurate, by comparing to the
floating point result.

Test\_2a is easily gamable, for $O(n^2)$ extra cost, by replacing
$A$ by $A\cdot D_3$ and $B$ by $D_3^{-1}\cdot B$, where the diagonal matrix $D_3$ is chosen
so that the i-th column of A and i-th row of B each has (about) the same norm,
effectively changing each row of $A\cdot D_3$ back to $x$ and each column of $D_3^{-1}\cdot B$
back to $x^T$. Note that this diagonal scaling  will
also greatly improve the accuracy
of Strassen, by reducing $\|A\| \cdot \|B\| = O(OV^2)$
to $\|A\cdot D_3\|\cdot \|D_3^{-1}\cdot B\| = O(1)$.

\item[Test\_2b:] 
We can try to prevent the gaming
of Test\_2a as follows: for each $i \in [1:n]$,
circularly shift each row $i$ of $A$ right by $i$, and circularly shift each column $i$ of $B$
down by $i$. Then each row and column of $A$ and of $B$ has the same norm, so any
diagonal scaling is pointless. To make this pattern
difficult to identify, we further randomly permute
the rows and columns of $A$, replacing $A$ by
$A' = P_1 \cdot A \cdot P_3$,
where $P_1$ and $P_3$ are random permutations, and
similarly randomly permute the rows and columns of $B$,
replacing it by $B' = P_3^T \cdot B \cdot P_2$. This maintains
the constant row and column norms of $A$ and $B$.
Since all the entries of $A'$ and $B'$ are positive,
an $O(n^3)$ algorithm using conventional floating point
will compute all entries of $A' \cdot B'$ to high
relative accuracy. In contrast, neither emulated arithmetic nor Strassen-like algorithms will attain
high relative accuracy for all entries,
in particular for the tiniest entries, which are on
the diagonal of $A \cdot B$, or the entries of
$A' \cdot B'$ corresponding to nonzero entries of
$P_1 \cdot P_2$.
This clearly determines whether an ``A'' is deserved or not.

\ignore{
\item[OLD Test\_2b:] We can try to prevent the gaming of Test\_2a as follows: for each $i \in [1:n]$,
circularly shift each row of $A$ right by $i$, and circularly shift each column of $B$
down by $i$. Then each row and column of $A$ and of $B$ has the same norm, so any
diagonal scaling is pointless, and the diagonal entries of $A\cdot B$ should all equal
$x\cdot x^T$ to high relative accuracy in floating point. We should also reduce the range of
numbers in the diagonal matrix $D$ used to compute $y=x\cdot D$ to go from roughly $\sqrt{UN}$ to
$\sqrt{OV}$, to avoid underflow and overflow. \ejr{Is this needed in the new test? Haven't thought about it, just noticed in my code.}
} 

\item[Test\_2c:] 
A possible, but expensive
and unlikely way to game Test\_2b would be to
first check if all the entries of $A'$ have the same
sign, and similarly for $B'$. In this rare case, the
algorithm could decide to switch to an $O(n^3)$
algorithm with conventional floating point, to guarantee high relative accuracy in each entry.
To prevent
this, we can pick a small random subset $S_A$ of random rows of $A'$ and negate some of their entries in each row.
We then pick a random subset $S_B$ of columns of $B'$,
and negate some of their entries in each column.
Finally, we let $S$ be a diagonal matrix with
random $\pm 1$ diagonal entries, and replace $A'$
by $A'' = A' \cdot S$ and $B'$ by $B'' = S \cdot B'$.
Then all the entries of $A'' \cdot B''$ not in rows
$S_A$ or columns $S_B$ will be sums of positive numbers,
and only all computed with high relative accuracy if
an ``A'' grade is deserved.
Detecting this pattern
and only using an $O(n^3)$ algorithm with conventional
floating point when it is detected would be even more
expensive than gaming Test\_2b. 

\ignore{
\item[OLD Test\_2c:] Unfortunately, Test\_2b is still gamable,
by using conventional floating point to compute just
the diagonal entries of $A \cdot B$ to high accuracy,
for an extra $O(n^2)$ cost. Instead, we can randomize the
locations of these $n$ accurate entries, so that it is
too expensive to detect their locations and compute them
with extra precision. Simply choose random permutation
matrices $P_1$ and $P_2$, and replace $A$ by $P_1 \cdot A$,
and $B$ by $B \cdot P_2$. The entries of the product
to test for accuracy are the nonzero entries of $P_1 \cdot P_2$.
} 
\end{description}

So far, we claim that Test\_1c and Test\_2c
can reliably distinguish algorithms deserving
an ``A'' from those that do not. (The other tests
are worth keeping too, because more testing means
more reliability.)
It remains to distinguish the ``B''s from the ``C''s.
This is relatively easy:

\begin{description}
\item[Test\_3:] Choose the entries of $A$ and $B$ to
be random numbers in the range $[1,2]$, so that $A \cdot B$ is easy to compute accurately. Let $D_1$
and $D_2$ be diagonal matrices with  random
diagonal entries
whose entries range from roughly $\sqrt{UN}$
to $\sqrt{OV}$,
and multiply $D_1 \cdot A$ times
$B \cdot D_2$. Each entry of the product should
equal $D_1 \cdot (A \cdot B) \cdot D_2$ to high
relative accuracy if and only if a ``B'' grade is
deserved.
\end{description}

\ignore{
{\em Should I add the test I proposed, to compare
$A \cdot B$ to
$P_1^T \cdot [ (P_1 \cdot A) \cdot (B \cdot P_2)] \cdot P_2^T$,
where $P_1$ and $P_2$ are random permutations?
The idea was that Strassen would get more varied
answers than any $O(n^3)$ algorithm, modulo the usual
questions about reproducible floating point summation. Here is a draft:}
} 

We mention one more pair of tests, initially
designed just to identify Strassen-like 
algorithms, but which also demonstrate the 
possibility of distinguishing emulation-based
algorithms from conventional $O(n^3)$ floating
point algorithms.
\begin{description}
\item[Test\_4a:] This test is intended 
to detect a Strassen-like algorithm, but again may be gameable,
or sensitive to the underlying implementation. The idea is to
compare the results of computing
$A \cdot B$ and 
$P_1^T \cdot ((P_1 \cdot A) \cdot (B \cdot P_2)) \cdot P_2^T$,
where $P_1$ and $P_2$ are random permutations. 
An $O(n^3)$ algorithm, using any arithmetic, should get
identical answers, assuming (1) the sums are computed in the
same order (which is not guaranteed in a parallel
computing environment \cite{ahrens2020reproblas}),
and (2) emulation only uses the
entries of row $i$ of $A$ (resp. column $j$ of $B$)
to choose the representation of row $i$ of $A$
(resp. column $j$ of $B$). On the other hand, a Strassen-based algorithm
will necessarily make quite different rounding errors
during the recursion. The most straightforward approach
is to use the examples in Test\_2c, and repeat for multiple
random permutations $P_1$ and $P_2$.
\ignore{
We can try to increase the number of entries of the product
to test for reproducibility as follows.
Start with the entries of $A$ and $B$ being random positive
numbers between 1 and 2, so any $O(n^3)$ algorithm will get
each entry of the product to high relative accuracy.
We can increase the difference in
the rounding errors by changing $A$ to $D_1 \cdot A$
and $B$ to $B \cdot D_2$ where $D_1$ and $D_2$ have
diagonal entries with a wide range of magnitudes, but
this is gameable by diagonal scaling. We can try to
avoid this as follows: Choose a random $i$ and $j$,
a large scalar $z$ close to the largest entries of
$A$ and $B$, and let $e$ be a
$n/2 \times 1$ vector of ones.
Then let columns $i$ and $j$ of $A$ be $[0 \cdot e; z \cdot e]$
and $[z \cdot e; 0 \cdot z]$, resp., and
rows $j$ and $i$ of $B$ be the transposes of these vectors, resp.
These modified $A$ and $B$ will have nearly identical row
and column norms, making diagonal scaling pointless,
and their product will add $z^2$ just to the upper right
and lower left $n/2 \times n/2$ submatrices of the product,
leaving significant changes in roundoff error visible in
the other two corners.
} 
\item[Test\_4b:] We modify Test\_4a to compare the 
results of $A \cdot B$ and 
$(A \cdot P_3) \cdot (P_3^T \cdot B)$,
where
$P_3$ is a random permutation, that should also not affect the
answer in exact arithmetic. It will however change the order of
summation in dot products computed using conventional floating point.
\end{description}

As we will see later in 
Table~\ref{tab:test3_and_test4},
Test\_4b is the only test that can by itself
(so far)
distinguish $O(n^3)$ algorithms using
conventional floating point, emulation based algorithms, and Strassen-based algorithms.

All these tests can be further modified
to refine their ability to detect
algorithmic details as follows
(we leave this as future work):
\begin{enumerate}
\item In case an emulation approach
groups rows of $A$ or columns of $B$ into
blocks of consecutive entries, and
chooses a single exponent to represent
all these entries, having large and
small entries stored close together will
increase the chance of losing precision
in the smaller entries. This suggests
choosing the diagonal matrices in
Tests 2a, 2b and 2c to have large and
small entries (nearly) adjacent, to see if
this impacts the precision of the entries
of $C$ that should be accurate.
\item In case an emulation approach uses
a floating point format with
a narrower exponent range, the use of
diagonal matrices with entries ranging
from $\sqrt{UN}$ to $\sqrt{OV}$ may result
in overflows, or if diagonal scaling is used,
underflows that cause inaccuracies in
entries that should be accurate. In this
case, one could try replacing
$UN$ and $OV$ with values
that are closer together in magnitude to
see if the inaccuracy goes away, and
so report this by possibly giving
different grades
depending on the range of values in
$A$ and $B$.
\end{enumerate}

\ignore{ 
We mention one more possibility
to create tests that are even more resistant to gaming. We note that very aggressive gaming is more likely
to be used by an AI attempting to pass our tests, rather than a human developer, as has been observed in other situations.
Our idea is to randomly ``mix'' all the above
tests into one matrix multiplication:
1) Pick random matrix dimensions $A^{m \times k}$ and $B^{k \times n}$.
2) Choose a random subset of tests 1c, 2c, 3 and 4b. 
3) For each test in the subset, select a random, nonoverlapping
subset of rows of $A$ and columns of $B$.
4) For each test in the subset, insert
the test matrices 
into the selected subsets of nonoverlapping rows 
of $A$ and columns of $B$, and test the
accuracy of the corresponding nonoverlapping rectangular subblock of $C$. 
Since the goal of the implementor
is to attain a certain grade while running as fast as possible, 
the cost of discovering these
randomly inserted test cases, and using 
a custom ``gamed'' solution that is fast and accurate enough only on these test
cases, is likely to be quite high,
even more strongly discouraging gaming.
}
Here is another way
to create tests that are even more resistant to gaming. We note that very aggressive gaming is more likely
to be used by an AI attempting to pass our tests, rather than a human developer, as has been observed in other situations.
Our idea is to randomly ``mix'' all the above
tests into one matrix multiplication:
1) Pick random matrix dimensions $A^{m \times k}$ and $B^{k \times n}$.
2) Choose a random subset of tests 1c, 2c, 3 and 4b. 
3) For each test in the subset, select a random, nonoverlapping
subset of rows of $A$ and columns of $B$.
4) For each test in the subset, insert
the test matrices 
into the selected subsets of nonoverlapping rows 
of $A$ and columns of $B$, and test the
accuracy of the corresponding nonoverlapping rectangular subblock of $C$.
Since the goal of the implementor
is to attain a certain grade while running as fast as possible, 
the cost of discovering these
randomly inserted test cases, and using 
a custom ``gamed'' solution that is fast and accurate enough only on these test
cases, is likely to be quite high,
even more strongly discouraging gaming.

To support both reproducibility and gaming avoidance, our publicly
available test code provides default seeds for random number
generation, so default test matrices are bitwise reproducible,
but allows users the option of providing their own seeds. We will use
this feature to use our own (secret) random number seeds to
conduct tests with even less opportunity for gaming, and advise
others to do so as well.

\ignore{
{\em One further topic to consider is whether we should
modify these test cases to make sure that large
and tiny entries of $A$ are nearby when they are in the
same row, and large and tiny entries of $B$ are nearby
when they are in the same column, in case this
impacts how emulation is done, \textit{i.e.,} make it more
likely that tiny entries are not preserved.}
} 

\subsection{Test cases for evaluating $f_X(n)$}

We first consider an implementation of
\gemm that deserves an ``A''. In this case we can use the 
classical error analysis of dot products \cite{Higham2002}
to get the bound $f_A(n) = n$ in the case of square
$A$ and $B$. Given rectangular
matrices $A^{m \times k}$ and $B^{k \times n}$, then
$f_A(m,k,n) = k$. This bound holds independent of the
order of summation, but constructing an example that
attains it requires some assumptions both about the
order of summation and rounding mode. Here we assume
the summation is done sequentially as follows: 
$s = 0$; for $l = 1:k$, $s = s + \matentry{A}{i}{l} \cdot \matentry{B}{l}{j}$,
using the default round-to-nearest-even rounding mode.
We will construct a single dot product $x \cdot y^T$ that
nearly attains this bound, and then let every row of $A$ 
equal $x$, and every column of $B$ equal $y^T$,
so every entry of $A \cdot B$ is identical.

We assume we  are working in FP64 arithmetic, with 
$\epsilon = 2^{-53}$.   
Then we let $\vecentry{x}{1}=\vecentry{y}{1}=1$,
and for each $i>1$, we choose $\vecentry{x}{i}$ and $\vecentry{y}{i}$ to be
powers of 2 whose product is $\vecentry{x}{i} \cdot \vecentry{y}{i} = \epsilon$.
Since we round to nearest even, 
$fl(1+ \vecentry{x}{i}\cdot \vecentry{y}{i}) = fl(1+\epsilon) = 1$, so
the true dot product is $1+(n-1)\epsilon$ while the
computed dot product is $1$. The error $(n-1)\epsilon$
nearly attains the upper bound of $n\epsilon (1+O(n\epsilon))$. This yields $f_A(n) \approx n$.
It also yields $f_B(n) \approx n$, but as we will see
below, $f_B(n)$ can be much larger for a GEMM that
earns an ``A''.

If the summation is done using a binary tree instead,
then $f_A(n)$ is at most $\log_2(n)$. When $n = 2^k$ is
a power of 2, we can construct a worst case example 
as follows. We let $\vecentry{x}{i}$ and $\vecentry{y}{i}$ denote the inputs
at the leaves of the tree, for $i=1$ to $n=2^k$.
We set $\vecentry{x}{1}=\vecentry{y}{1}=1$ as before, then
$\vecentry{x}{2} \cdot \vecentry{y}{2} = \epsilon$,
$\vecentry{x}{3} \cdot \vecentry{y}{3} = \vecentry{x}{4} \cdot \vecentry{y}{4} = \epsilon/2$,
$\vecentry{x}{5} \cdot \vecentry{y}{5} = ... = \vecentry{x}{8} \cdot \vecentry{y}{8} = \epsilon/4$,
... ,
$\vecentry{x}{2^{k-1}+1} \cdot \vecentry{y}{2^{k-1}+1} = ... = \vecentry{x}{2^k} \cdot \vecentry{y}{2^k} = \epsilon/2^{k-1}$.
We see that the true dot product is $1 + k \epsilon$,
but the computed dot product is 1.

Next, we consider a \gemm that gets a ``B'' grade.
As above, we construct vectors $x$ and $y$, and set
each row of $A$ to $x$, and each column of $B$ to $y^T$,
so that all entries of $A \cdot B$ are equal to 
$x \cdot y^T$. We let $x = [1, 2^{-k}, 1, 2^{-k}, ...]$
and $y = [2^{-k}, 1, 2^{-k}, 1, ...]$ so that
$x \cdot y^T = n2^{-k}$, assuming $n$ is even.
If $k$ is large enough, 
the ``slices'' of
$A$ and $B$ \cite{Uchino_Ozaki_Imamura_2024} 
containing the $2^{-k}$ entries will not be used,
and so $\dogemm{A}{B} = 0$, an error of $n2^{-k}$
in each entry, which can be close to 
$\max_k |\matentry{A}{i}{k}| \cdot \max_k |\matentry{B}{k}{j}|\cdot n \cdot \epsilon$
for all $i$ and $j$. This is close to the error bound in
section 5 of \cite{Uchino_Ozaki_Imamura_2024}. However,
if slices that are exactly 0 are not used, this could
mean that the slices containing $2^{-k}$ are in fact
used, requiring the following modification.
Assuming that a nonzero in any row or column of a slice
causes it to be retained, we can modify the first
row of $A$ and first column of $B$ to contain 
$[1, 2^{-1}, 2^{-2}, 2^{-3}, ...]$ to make sure all
slices are nonzero (if the matrix dimensions are small,
we could use $[1, 2^{-s}, 2^{-2s}, 2^{-3s}, ...]$
for an appropriate small integer $s$).

Next, we return to the question of how large $f_B(n)$ can
be for an algorithm that earns an ``A''.
We give a concrete example where $f_B(n) \approx .25 n^2$.
We note that $f_B(n)$ can be a factor of $n$ larger than
$f_A(n)$ because the components of
$\rowmax (|A|)\colmax(|B|)$ can be a factor of $n$ smaller
than $|A| \cdot |B|$.
We suppose $A$ and $B$ are both $n \times n$, where
each row of $A$ is identical to the row vector $x^T$,
and each column of $B$ is identical to the column
vector $y$, so we only need to consider the
error in the dot product $x^T y$. The idea is to
choose $x_i = 1 - e_i$ and $y_i = 1 + e_i$, where
$1 \gg e_i > 0$ is chosen as large as possible so 
that when $x_i \cdot y_i = 1 - e_i^2$ is added 
to the sum $\sum_{j=1}^{i-1} x_j \cdot y_j$, $e_i^2$ rounds away. Initially
$e_1 = 2^{-27}$ so that $1 - 2^{-2 \cdot 27}$ rounds
to 1 in FP 64, and $e_i$ gradually increases.
Concretely, we may choose $x_2 = \cdots =x_4 = 1-e_1$
and $y_2 =  \cdots = y_4 = 1 + e_1$,
$x_5 = \cdots = x_9 = 1 - 2 e_1$,
$y_5 = \cdots = y_9 = 1 + 2 e_1$,
$x_{10} = \cdots = x_{16} = 1 - 3 e_1$,
$y_{10} = \cdots = y_{16} = 1 + 3 e_1$, ... ,
$x_{k^2 +1} = \cdots = x_{(k+1)^2} = 1 - k e_1$,
$y_{k^2 +1} = \cdots = y_{(k+1)^2} = 1 + k e_1$.
One can confirm that each $-(i \cdot e_1)^2$ rounds away,
so $fl(\sum_{i=1}^{k^2} x_i \cdot y_i) = k^2$, with an
absolute error of $\approx .25n^2 \epsilon$ for large $n$,
so with $f_A(n) \approx .25n$, and $f_B(n) \approx .25 n^2$
(note that we assume that we sum in increasing order
of subscripts).

Constructing such (nearly) worst case examples for algorithms
with ``C'' grades is future work.

\section{Experimental Results for DGEMM}\label{sec:testresults}

This section is organized as follows. We first discuss the testing platform and setup used for our experiments. We then demonstrate the effectiveness of the described tests and their robustness to the discussed gaming strategies. Finally, we apply these tests to production DGEMM implementations to demonstrate how they might be used to generate insights for BLAS users.

\subsection{Setup}

{{C\nolinebreak[4]\hspace{-.05em}\raisebox{.4ex}{\tiny\bf ++}}} source code implementing these experiments is available at \url{github.com/Reference-LAPACK/grading-the-BLAS}. Except where noted, all experiments were executed on a single dedicated compute node on the Perlmutter machine at the National Energy Research Scientific Computing Center at Lawrence Berkeley National Laboratory\cite{perlmutter} (SUSE Linux, single AMD EPYC 7763 CPU, single NVIDIA A100 GPU) with binaries compiled using GCC 13.2.1 and CUDA Toolkit 12.9.41.

\subsection{Demonstrating the Effectiveness and Robustness of the Tests}
\label{sec:testresultspart1}

Given DGEMM implementations with known grades, we demonstrate the effectiveness of the discussed tests in discriminating between them, even when confronted with adversarial gaming strategies.

The grade ``A'' DGEMM is a naive triple loop implementation of DGEMM that calculates each entry of the output matrix as a single dot product. The ``A+'' implementation is the ``A'' implementation with 80-bit extended precision (\texttt{long double} in {{C\nolinebreak[4]\hspace{-.05em}\raisebox{.4ex}{\tiny\bf ++}}}, GCC 13.2.1, x86-64 SUSE Linux, \texttt{-O0} optimization).

The grade ``B'' DGEMMs are implementations of the Ozaki I \cite{Ootomo_Ozaki_Yokota_2024} and Ozaki II \cite{uchino2025high} schemes which emulate DGEMM using INT8 Tensor cores. The implementations are taken from the repos containing the source code accompanying the aforementioned research papers describing these contributions, found at \url{https://github.com/enp1s0/ozIMMU} and \url{https://github.com/RIKEN-RCCS/GEMMul8} respectively. These implementations provide parameters that adjust their performance/precision envelope, thereby affecting the value of $f(n)$ seen in the bounds that determine each grade. In order to maximize accuracy, our tested implementations are parameterized to use the maximum number of slices / moduli available. Furthermore, we use the ``accurate'' mode of Ozaki II described in \cite{uchino2025high}.

The grade ``C'' DGEMM is an implementation of the original Strassen's algorithm \cite{Strassen69} with a basecase of N=2.

For clarity, all tables and figures describe experiment results for square \gemm{}s of dimension 32. We empirically verified that the observed behavior also holds for square \gemm{}s of dimensions 64, 128, 256, 512, and 1024.

\subsubsection{Test 1 Series}

Recall that Test 1 is designed to distinguish grade ``A'' \gemm implementations from grade ``B'' and ``C'' implementations and that subsequent tests in the series are designed to be robust in the face of gaming strategies involving the scaling of rows/columns of A/B respectively to have unit norm (``outer scaling'') and the zeroing of elements of C below a specified threshold.

We parameterize our implementations of the tests thusly: For Test\_1a, we randomly select 50\% of the rows/columns of A/B to be zero which, for a square $n \times n$ \gemm, should yield $n/2$ rows/columns of zeros, \textit{i.e.,} $3n^2/4$ total zeros. For Test\_1b, we randomly select 50\% of the rows/columns of A/B and we make those selected rows/columns complimentarily 50\% sparse; for a square $n \times n$ \gemm, this should yield $n/2$ zeros. For Test\_1c, we start with the same parameterization as Test\_1b. Nonzero entries in the complimentarily sparse rows/columns are then scaled to be on the order of $2^{500}$. We then randomly select 50\% of the targeted rows/columns in which to add a small nonzero on the order of $2^{-500}$. For a square $n \times n$ \gemm, this should yield $n/4$ zeros and $n/4$ small nonzeros on the order of $2^{-1000}$.

\begin{figure}[p]
    \centering
    \begin{subfigure}{0.7\textwidth}
        \centering
        \includegraphics[width=0.55\textwidth]{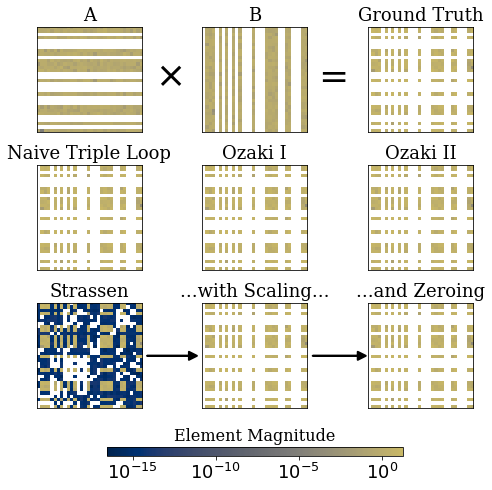}
        \hfill
        \raisebox{3em}{\fbox{\includegraphics[width=0.42\textwidth]{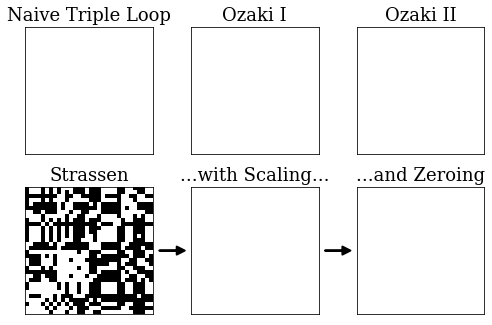}}}
        \caption{Test\_1a}
        \label{fig:Test_1a}
    \end{subfigure}
    \begin{subfigure}{0.7\textwidth}
        \centering
        \includegraphics[width=0.55\textwidth]{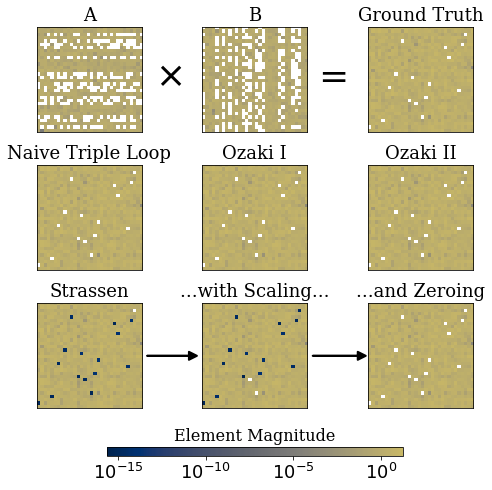}
        \hfill
        \raisebox{3em}{\fbox{\includegraphics[width=0.42\textwidth]{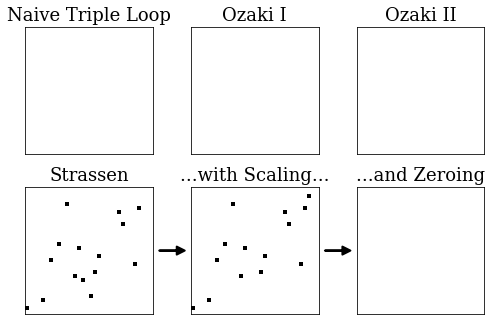}}}
        \caption{Test\_1b}
        \label{fig:Test_1b}
    \end{subfigure}
    \begin{subfigure}{0.7\textwidth}
        \centering
        \includegraphics[width=0.55\textwidth]{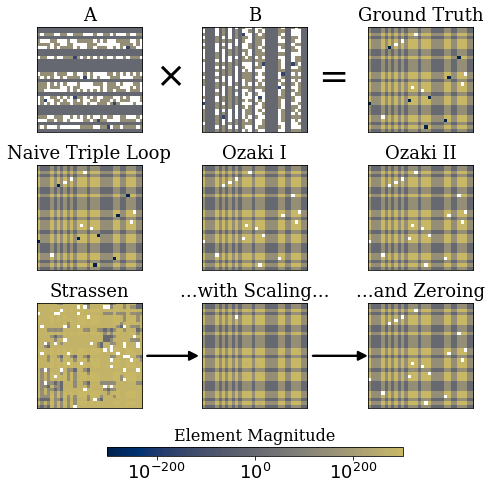}
        \hfill
        \raisebox{3em}{\fbox{\includegraphics[width=0.42\textwidth]{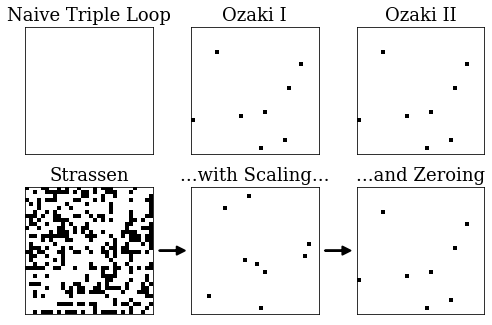}}}
        \caption{Test\_1c}
        \label{fig:Test_1c}
    \end{subfigure}
    \caption{Test\_1 series. The left plots depict element magnitudes with exact zeros in white. The right plots depict in black elements whose sign $(+,-,0)$ differs from the ``A+'' implementation. Square matrix size n=32. }
    \label{fig:test1_series}
\end{figure}

The signal for the Test 1 series is the loss of expected zeros and small nonzeros described above and henceforth referred to as ``target values''.  DGEMMs deserving of an ``A'' must accurately compute all target values with expected zeros remaining exactly zero and small nonzero values remaining nonzero and small with the correct sign. 
The results are summarized in Table \ref{tab:test1_series} with greater detail offered in Figure \ref{fig:test1_series}.

\begin{table}
    \centering
    \caption{Fraction of incorrect target values across the Test 1 series, for $n=32$. For Test\_1a and Test\_1b, target values are the zero elements induced by the sparsity structure of A and B. For Test\_1c, target values are these zero elements plus the very small nonzero $x^2$ values introduced to thwart the zeroing gaming strategy. Nonzero table entries are highlighted in red. DGEMM implementations deserving an ``A'' must accurately compute all target values with expected zeros remaining exactly zero and small nonzero values remaining nonzero and small.}
    \label{tab:test1_series}
    \begin{tabular}{lrrr}
        \toprule
        & \multicolumn{3}{c}{Fraction of wrong target values} \\
        \cmidrule(lr){2-4}
        DGEMM Implementation & Test\_1a & Test\_1b & Test\_1c \\
        \midrule
        Naive Triple Loop              & 0 / 768                            & 0 / 16                           &  0 / 16                            \\
        Ozaki I                        & 0 / 768                            & 0 / 16                           & \textcolor{red}{\textbf{8 / 16}}   \\
        Ozaki II                       & 0 / 768                            & 0 / 16                           & \textcolor{red}{\textbf{8 / 16}}   \\
        Strassen                       & \textcolor{red}{\textbf{524 / 768}} & \textcolor{red}{\textbf{15 / 16}} & \textcolor{red}{\textbf{14 / 16}} \\
        Strassen + Outer Scaling             & 0 / 768                            & \textcolor{red}{\textbf{14 / 16}} & \textcolor{red}{\textbf{15 / 16}}   \\
        Strassen + Outer Scaling + Zeroing   & 0 / 768                            & 0 / 16                           & \textcolor{red}{\textbf{8 / 16}} \\
        \bottomrule
    \end{tabular}
\end{table}

First, we note how Test\_1a correctly distinguishes the Strassen implementation whose distinct rounding errors prevent perfect cancellation and populate the expected zero rows and columns with small nonzero entries, a pattern visible in Figure \ref{fig:Test_1a}. Furthermore, we see how the addition of outer scaling allows the Strassen implementation to game the test.

Moving to Test\_1b, Figure \ref{fig:Test_1b} depicts how the complimentary sparsity pattern of select rows / columns of the A/B matrices thwarts that outer scaling gaming strategy: Strassen and Strassen+Scaling fail to compute nearly all of the 16 expected zeros, yielding small nonzeros at 15 and 14 of them, respectively. The addition of zeroing elements below a specified threshold ($1000 \, \epsilon \, \|C\|_{\max}$ for these experiments) successfully games the test.

Finally, the additional anti-gaming measures adopted for Test\_1c successfully thwart the aforementioned outer scaling and zeroing strategies for Strassen while also distinguishing the Ozaki I and Ozaki II schemes from the grade ``A'' implementation. Figure \ref{fig:Test_1c} offers more detail:
Strassen alone loses the output structure entirely, both dropping target values and even flipping the sign of many of the non-target values. Strassen+Scaling recovers the overall structure but corrupts 14 of the 16 zeros and very small nonzero $x^2$ values into very large elements, a failure mode which the zeroing threshold cannot correct. Lastly, both emulation schemes flush the eight small $x^2$ values to zero.

\clearpage
\newgeometry{margin=0.5in}
\clearpage
\restoregeometry

\subsubsection{Test 2 Series}

Recall that Test 2 is designed to distinguish grade ``A'' \gemm implementations from grade ``B'' and ``C'' implementations and that subsequent tests in the series are designed to be robust in the face of the gaming strategy involving the scaling of columns/rows of A/B respectively to have roughly the same norms (``inner scaling'').

All input matrices in the Test 2 series have the property that, for some $k$, their values span from $2^{-k}$ to $2^{k}$ with exponents spaced linearly over $[-k,k]$. For Test\_2a, $k = 1000$. For Test\_2b and Test\_2c, $k=500$. For Test\_2c, we randomly select 25\% of the rows/columns of A/B and randomly select 25\% of their values to negate.

\begin{figure}[p]
    \centering
    \begin{subfigure}{\textwidth}
        \centering
        \includegraphics[width=0.48\textwidth]{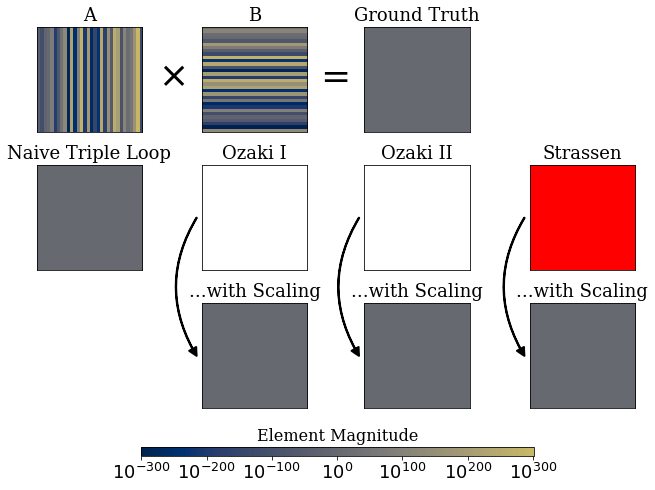}
        \hfill
        \raisebox{0.5em}{\fbox{\includegraphics[width=0.48\textwidth]{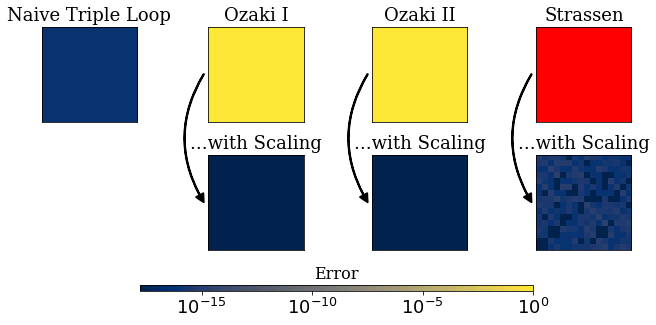}}}
        \caption{Test\_2a}
        \label{fig:test_2a}
    \end{subfigure}
    \begin{subfigure}{\textwidth}
        \centering
        \includegraphics[width=0.48\textwidth]{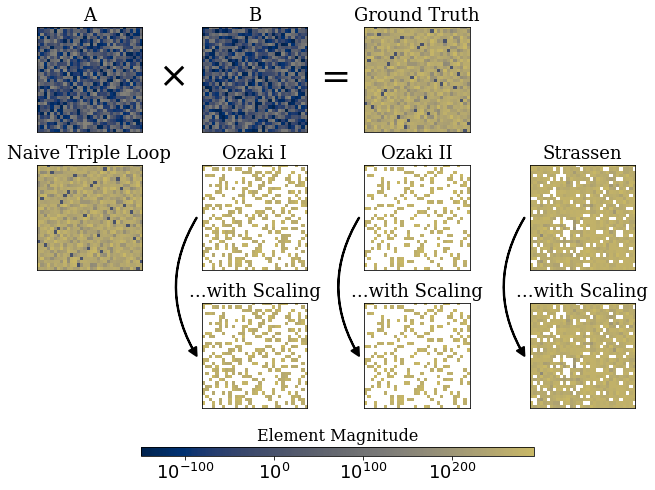}
        \hfill
        \raisebox{0.5em}{\fbox{\includegraphics[width=0.48\textwidth]{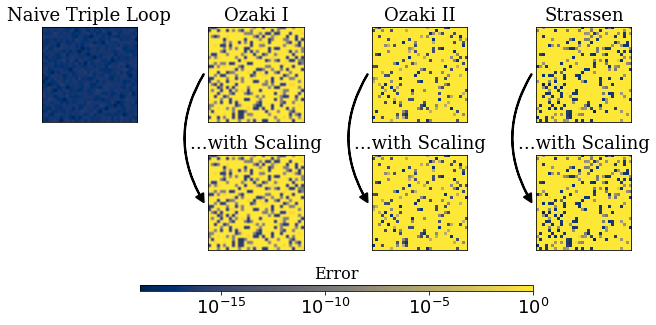}}}
        \caption{Test\_2b}
        \label{fig:test_2b}
    \end{subfigure}
    \caption{Test\_2 series. The left plots depict element magnitudes and the right plots depict the error with respect to the ``A+'' implementation.  For clarity of visualization, error is calculated as $|C - C_{ref}| / (|C| + |C_{ref}|)$ so that, for finite operands, the error lies in the interval $[0,1]$. We map corner cases of $0/0$ to $0$ and $\EV/\EV$ to \NaN where $\EV \in \{\NaN, \pm\Inf\}$. Zero is depicted in white. \NaN is depicted in red. Square matrix size n=32.}
    \label{fig:test2_series}
\end{figure}

We take the signal for the Test 2 series to be the componentwise max of $f_A(n) = |\dogemm{A}{B} - A \cdot B \,| / (\epsilon \, |A| \cdot |B|)$ 
where the ``A+'' implementation serves as the ground truth. For Test\_2c which includes random negation of input elements, the componentwise max is taken over all entries that are sums of positive numbers.
DGEMM implementations deserving an ``A'' must satisfy the worst-case bound $f_A(n) < n$. The results are summarized in Table \ref{tab:test2_series} with greater detail offered in Figure \ref{fig:test2_series}.

To start, Test\_2a successfully distinguishes Ozaki I, Ozaki II, and Strassen from the grade ``A'' implementation. In Figure \ref{fig:test_2a}, we can see how the difficulty caused by the wide dynamic range of $A$ and $B$ manifests in the output of the non-``A'' DGEMM implementations.

The all-zero output observed for both Ozaki I and Ozaki II reflects a property shared by both schemes. Each row of $A$ and column of $B$ is quantized using a single per-row or per-column scaling factor that fails to represent the wide dynamic range and flushes the smaller values to zero. This yields INT8 representations of $A$ and $B$ with surviving nonzero supports at disjoint indices, leading to dot products that are exactly zero.

The all-\NaN output for Strassen stems from two types of floating-point exceptions. First, the intermediate products, \textit{e.g.,} $M_1 = (\matblock{A}{1}{1} + \matblock{A}{2}{2}) \cdot (\matblock{B}{1}{1} + \matblock{B}{2}{2})$, cause overflows to \Inf. Second, combining those infinite intermediate products to form the result, \textit{e.g.,} $\matblock{C}{1}{1} = M_1 + M_4 - M_5 + M_7$, causes invalid exceptions resulting in \NaN values. Once generated, these \NaNs propagate throughout the recursion.

These failure modes are alleviated by the addition of inner scaling, thus allowing Ozaki I and Ozaki II to game Test\_2a. We suspect that the lower $f_A(n)$ for the emulation schemes with inner scaling is due to their error-free summation of INT8 operands. We note that Strassen with scaling still fails to get an ``A'', but comes close, with $f_A(32) = 35.6 > n= 32$, motivating our stronger Test\_2b.

Test\_2b adds circular shifting and random permuting of the input matrices to defeat the inner scaling strategy. Finally, Test\_2c adopts a negation strategy to thwart an expensive but unlikely means of gaming Test\_2b which we have not implemented; results mirror those of Test\_2b. In Figure \ref{fig:test_2b}, we can see how the difficulty caused by circular shifting manifests in the output of the non-``A'' DGEMM implementations: many zeros where there should be none and the loss of all small nonzeros.

\clearpage
\newgeometry{margin=0.5in}
\begin{table}
    \centering
    \caption{ $f_A(n)$ calculated with respect to the ``A+'' DGEMM across the Test 2 series. Failing cases are highlighted in red and bolded. DGEMM implementations deserving of an ``A'' must satisfy the worst-case bound $f_A(n) < n$; here $n=32$.}
    \label{tab:test2_series}
    \begin{tabular}{llll}
        \toprule
        & \multicolumn{3}{c}{ $f_A(n)$} \\
        \cmidrule(lr){2-4}
        DGEMM Implementation & Test\_2a & Test\_2b & Test\_2c \\
        \midrule
        Naive Triple Loop              & $1.59 \times 10^{0}$                                       & $3.27 \times 10^{0}$                                       & $1.65 \times 10^{0}$                                       \\
        Ozaki I                        & \mfailcase{9.01 \times 10^{15}}             & \mfailcase{9.01 \times 10^{15}}             & \mfailcase{9.01 \times 10^{15}}             \\
        Ozaki I + Inner Scaling              & $2.85 \times 10^{-2}$                                      & \mfailcase{9.01 \times 10^{15}}             & \mfailcase{9.01 \times 10^{15}}             \\
        Ozaki II                       & \mfailcase{9.01 \times 10^{15}}             & \mfailcase{2.02 \times 10^{19}}             & \mfailcase{2.02 \times 10^{19}}             \\
        Ozaki II + Inner Scaling             & $2.85 \times 10^{-2}$                                      & \mfailcase{2.02 \times 10^{19}}             & \mfailcase{2.02 \times 10^{19}}             \\
        Strassen                       & \failcase{\NaN}                               & \mfailcase{3.00 \times 10^{290}}            & \mfailcase{6.31 \times 10^{289}}            \\
        Strassen + Inner Scaling             & \mfailcase{3.56 \times 10^{1}}              & \mfailcase{3.00 \times 10^{290}}            & \mfailcase{6.31 \times 10^{289}}            \\
        \bottomrule
    \end{tabular}
\end{table}

\clearpage
\restoregeometry

\subsubsection{Test 3 and the Test 4 series}

Recall that Test 3 is designed to distinguish between grade ``B''  and grade ``C'' \gemm implementations and the Test 4 series is designed to distinguish Strassen-like algorithms.

In our implementation of Test\_3, we use scaling matrices with diagonal entries of the form $2^k$ with $k$ linearly-spaced between $-500$ and $500$, arranged in alternating large and small magnitudes to maximize the probability of pathological cancellation using a Strassen-like algorithm. The signal for Test\_3 is the componentwise max of $f_B(n) = |\dogemm{A}{B} - A \cdot B \,| / (\epsilon \, \rowmax(|A|) \cdot \colmax(|B|))$ where the ``A+'' implementation serves as the ground truth. DGEMM implementations deserving a ``B'' must satisfy the worst-case bound $f_B(n) < n^2$.

In our implementation of the Test 4 series, we use the input matrices from Test\_2c and apply 10 random permutations.
We take the signal for the Test 4 series to be the maximum componentwise relative error of the output values over those 10 random permutations with respect to the output of the same DGEMM implementation under test without permutation. Relative error is calculated as $|C_{perm} - C_{ref}| / |C_{ref}|$.

\begin{figure}
    \centering
    \includegraphics[width=0.57\textwidth]{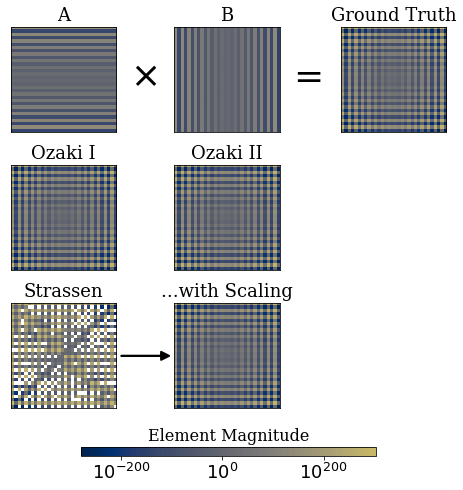}    
    \hfill    
    \fbox{\includegraphics[width=0.4\textwidth]{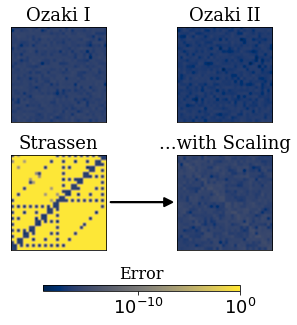}}
    \caption{Results of Test\_3. The left plot depicts element magnitudes and the right plot depicts the error with respect to the ``A+'' implementation.  For clarity of visualization, error is calculated as $|C - C_{ref}| / (|C| + |C_{ref}|)$ so that, for finite operands, the error lies in the interval $[0,1]$. Zero is depicted in white. Square matrix size n=32.}
    \label{fig:Test_3}
\end{figure}

Experimental results are presented in Table \ref{tab:test3_and_test4} with additional detail on Test 3 provided in Figure \ref{fig:Test_3}.

\begin{table}
    \centering
    \caption{Results for Test\_3 and the Test\_4 series.}
    \label{tab:test3_and_test4}
    \begin{tabular}{lcll}
        \toprule
        & Max componentwise $f_B(n)$ & \multicolumn{2}{c}{Max componentwise relative error} \\
        \cmidrule(lr){2-2} \cmidrule(lr){3-4}
        DGEMM Implementation & Test\_3 & Test\_4a & Test\_4b \\
        \midrule
        Naive Triple Loop              & $8.19 \times 10^{1}$                          & $0.0$                              & $6.22 \times 10^{-16}$ \\ 
        Ozaki I                        & $1.63 \times 10^{2}$                          & $0.0$                              & $0.0$                              \\
        Ozaki II                       & $3.48 \times 10^{1}$                          & $0.0$                              & $0.0$                              \\
        Strassen                       & \mfailcase{4.52 \times 10^{283}} & \failcase{\Inf} & \failcase{\Inf} \\
        Strassen + Outer Scaling       & \mfailcase{1.20 \times 10^{3}}   & \failcase{\Inf} & \failcase{\Inf} \\
        \bottomrule
    \end{tabular}
\end{table}
    
We see that Test\_3 successfully distinguishes the ``B''-grade emulation-based DGEMMs from the ``C''-grade Strassen implementation, whose  $f_B(n)$ blows up to $10^{283}$ and whose output contains many erroneous zeros. As expected, the addition of outer scaling dramatically improves Strassen's accuracy on Test\_3, nearly meriting a ``B'' with its max componentwise $f_B(32) = 1200 > 32^2 = 1024$. 

We see that the ``outer'' permutations of Test\_4a distinguish between the Strassen implementations (either
``B'' or ``C'') and the non-Strassen implementations.
The ``inner'' permutations of Test\_4b go further: the naive triple loop incurs the expected rounding errors; the emulation methods incur no rounding errors since their summations are exact; and both Strassen implementations incur the same expected catastrophic errors as Test\_4a.
In other words, Test\_4b by itself distinguishes between 
a conventional $O(n^3)$ floating point algorithm
(an ``A'' grade), an $O(n^3)$ emulation based algorithm
(Ozaki I and Ozaki II, a ``B'' grade), and 
a Strassen-based algorithm (whether it gets a ``B'' or ``C'').
This is useful for verifying the results of other tests.

\subsection{Grading DGEMM Implementations Across Various Scaling Factors and Problem Sizes}
\label{sec:testresultspart2}

Having demonstrated their effectiveness and robustness, we now apply these tests to demonstrate how they might be used to generate insights for BLAS users.

BLAS libraries can contain multiple DGEMM backends with distinct correctness profiles; in such cases, a call to that library's DGEMM API will invoke a decision procedure which will dispatch the calculation to a suitable backend based on some feature of the inputs (\textit{e.g.,} \cite{10.1145/3773656.3773670}). Consequently, tests should sweep the space of relevant features to determine what grades are achieved for the parameterizations of interest. For this demonstration, we employ a sweep of both the dimension of the \gemm, restricted to the square case for simplicity,
and the scaling factor used to control the dynamic range of A and B. Dimensions are powers of two between 32 and 1024, inclusive. Scaling factors $s$ are multiples of 100 between 0 and 500, inclusive, where the values of the input matrices range from $2^{-s}$ to $2^{s}$. Note that $2^{50  0}$ is close to the  
square root of the FP64 overflow threshold of
just under $2^{512}$, and
$2^{-500}$ is close to the square root of the
underflow threshold.
      
Note that, given the possible dependence of the grade on matrix dimensions and scale factors of $A$ and $B$,
our publicly available test code will let
the user sweep whatever subset of matrix dimensions (both
square and rectangular)
and dynamic ranges of $A$ and $B$ they find important,
and report both the results (including the grade) for
all the cases tested, and the overall grade (the lowest
grade over all the test cases). For example,
a performant Strassen implementation will shift to an
$O(n^3)$ algorithm for matrix dimensions below a threshold,
and so likely get a higher grade for smaller matrices.

\begin{figure}[p]
    \centering
    \begin{subfigure}{\textwidth}
        \centering
        \includegraphics[width=0.5\textwidth]{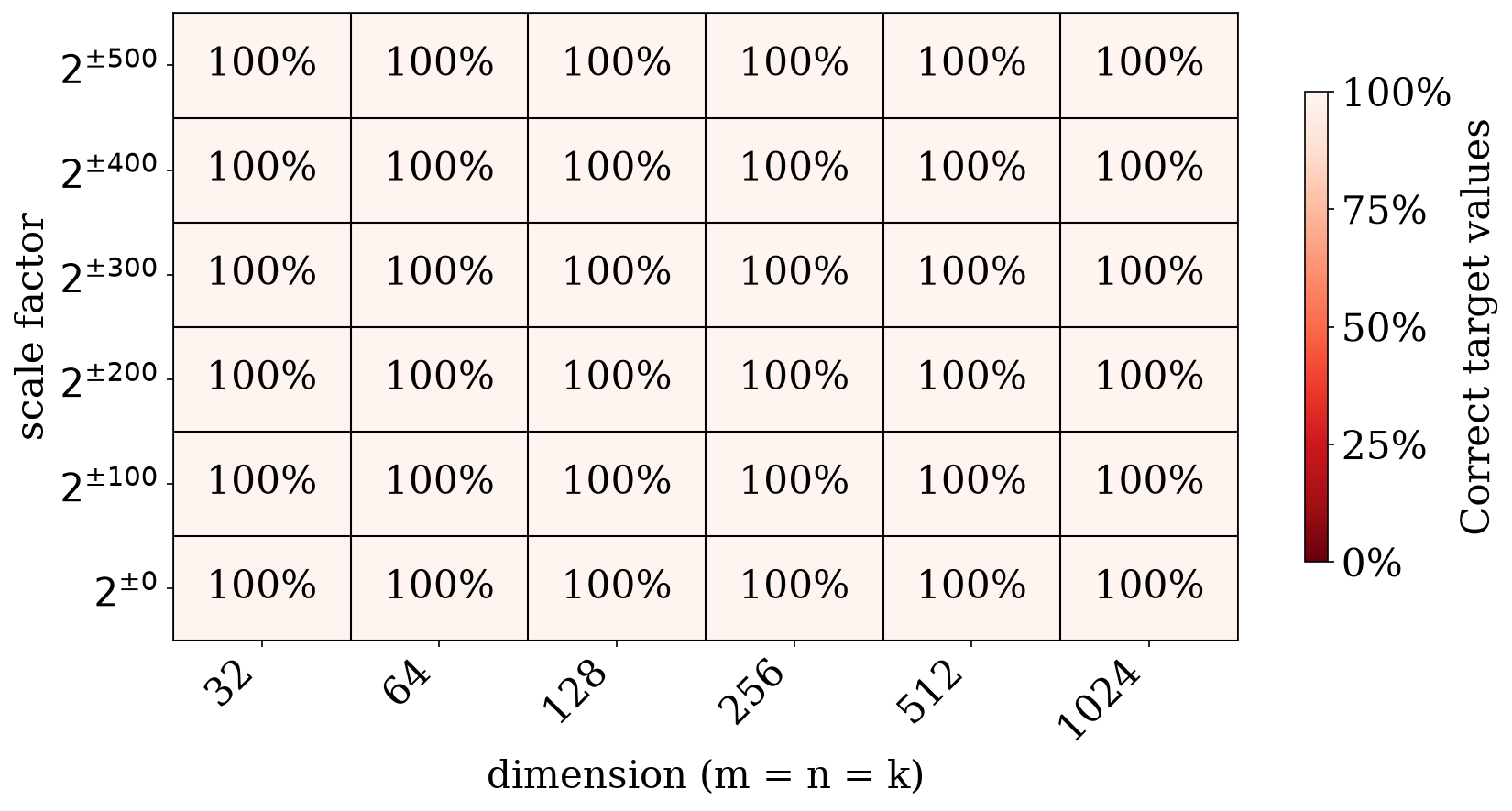}
        \caption{Test\_1c Sweep}\vspace{6mm}
    \end{subfigure}
    \begin{subfigure}{\textwidth}
        \centering
        \textcolor{lightgray}{\dashbox{\includegraphics[width=\textwidth]{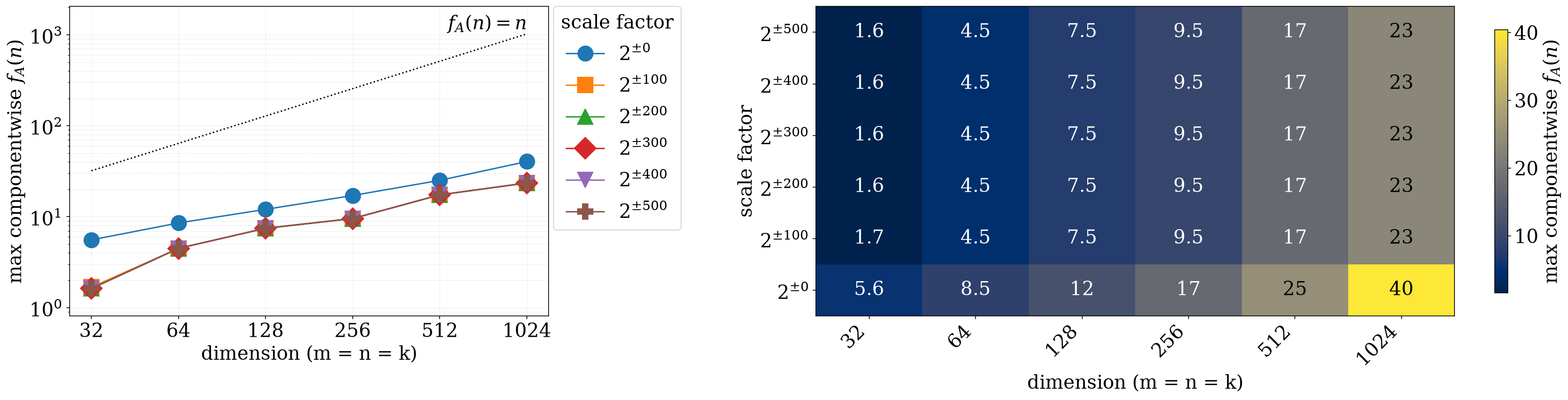}}}
        \vspace{-10mm}\caption{Test\_2c Sweep}\vspace{6mm}
    \end{subfigure}
    \begin{subfigure}{\textwidth}
        \centering
        \textcolor{lightgray}{\dashbox{\includegraphics[width=\textwidth]{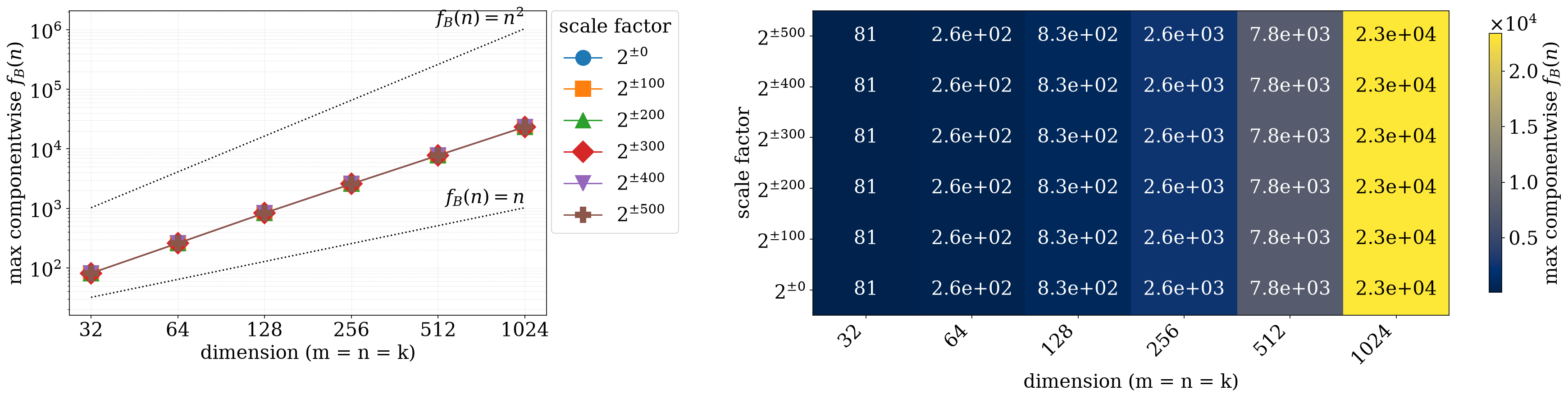}}}
        \vspace{-10mm}\caption{Test\_3 Sweep}\vspace{6mm}
    \end{subfigure}
    \begin{subfigure}{\textwidth}
        \centering
        \textcolor{lightgray}{\dashbox{\includegraphics[width=\textwidth]{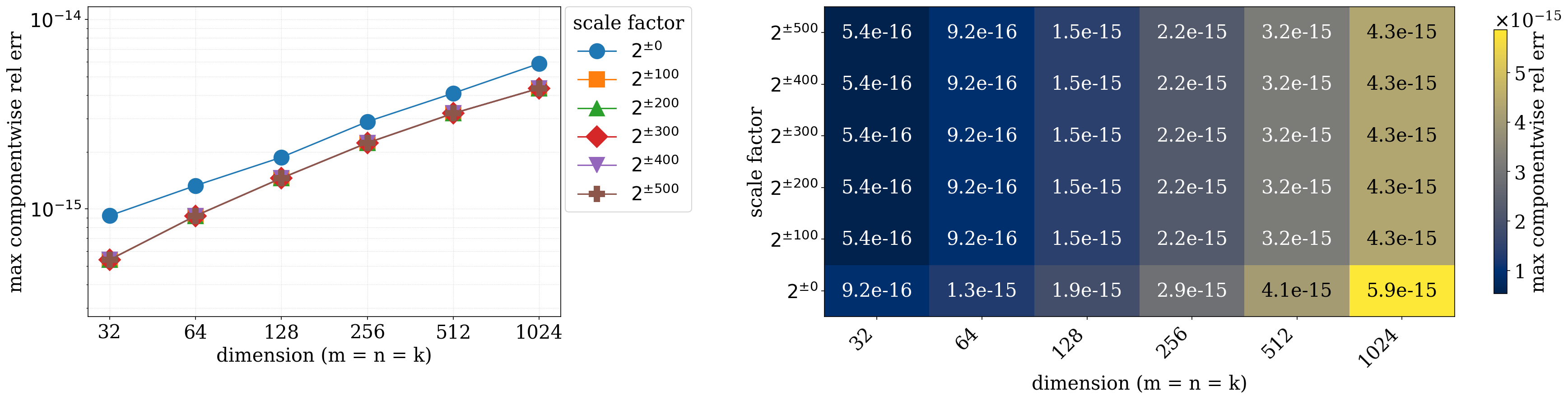}}}
        \vspace{-10mm}\caption{Test\_4b Sweep}\vspace{6mm}
    \end{subfigure}
    \caption{Naive Triple Loop}
    \label{fig:sweeps_naive_triple_loop}
\end{figure}

\begin{figure}[p]
    \centering
    \begin{subfigure}{\textwidth}
        \centering
        \includegraphics[width=0.5\textwidth]{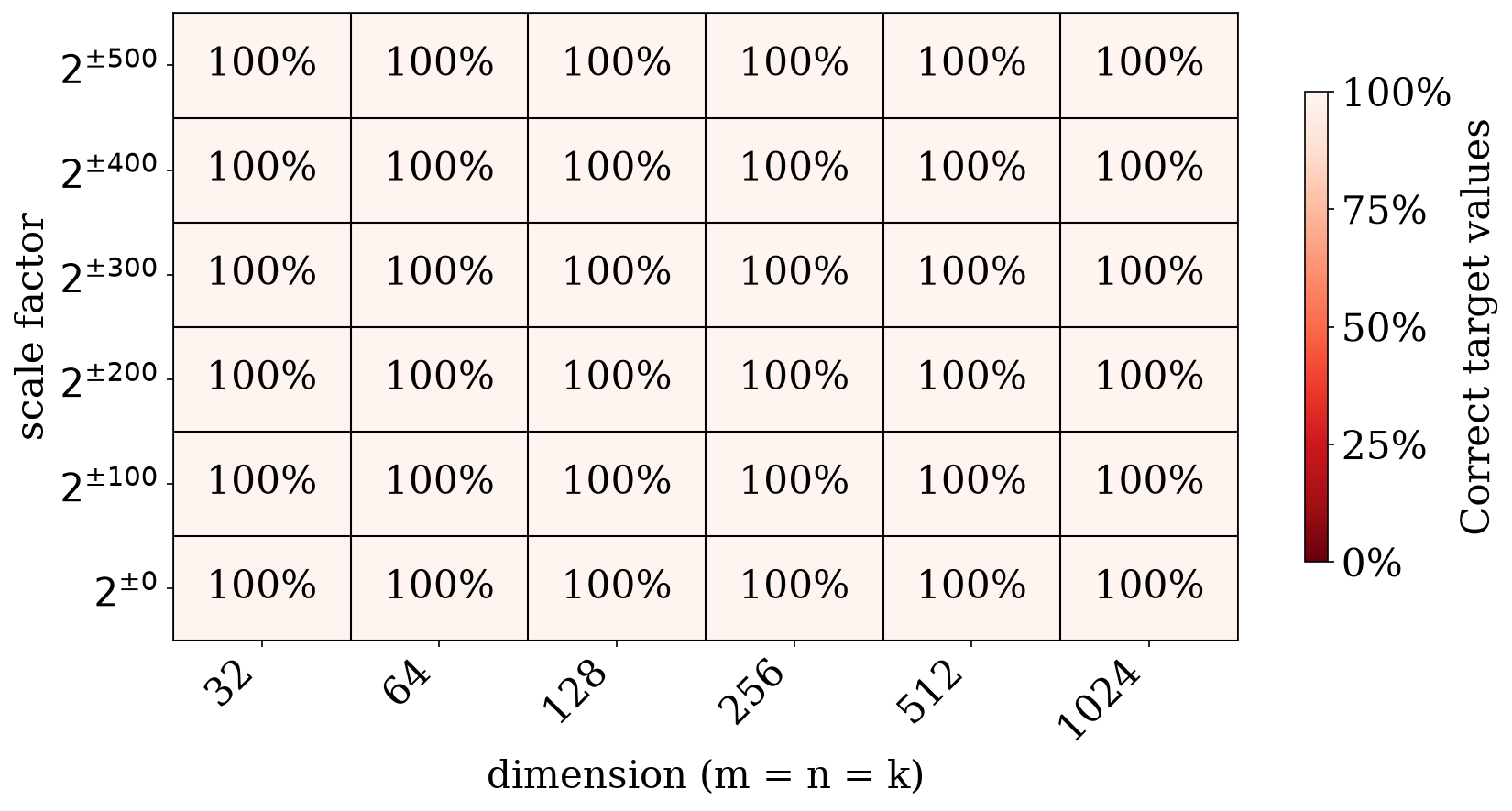}
        \caption{Test\_1c Sweep}\vspace{6mm}
    \end{subfigure}
    \begin{subfigure}{\textwidth}
        \centering
        \textcolor{lightgray}{\dashbox{\includegraphics[width=\textwidth]{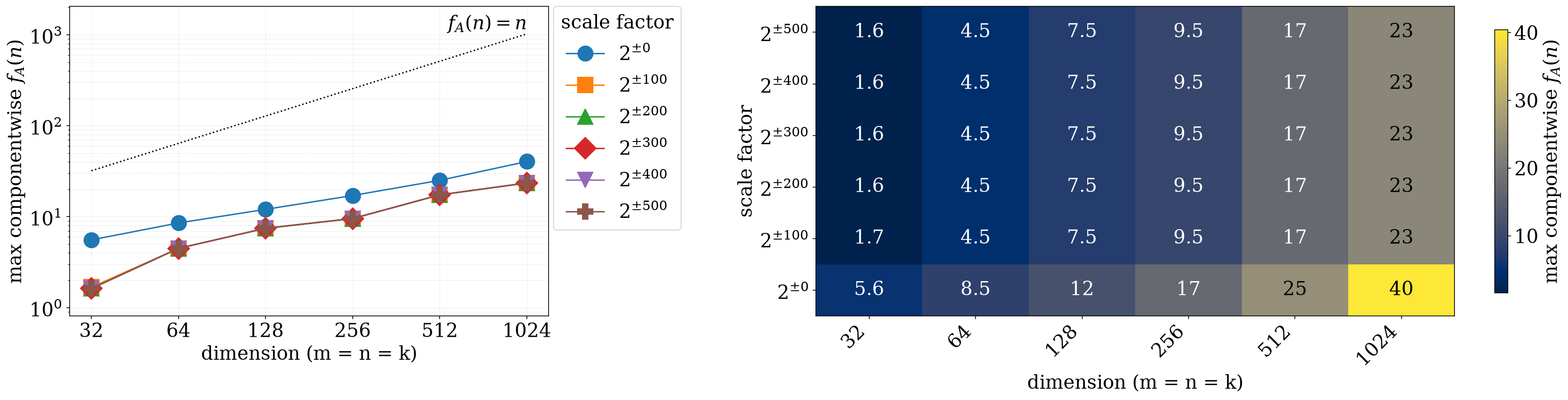}}}
        \vspace{-10mm}\caption{Test\_2c Sweep}\vspace{6mm}
    \end{subfigure}
    \begin{subfigure}{\textwidth}
        \centering
        \textcolor{lightgray}{\dashbox{\includegraphics[width=\textwidth]{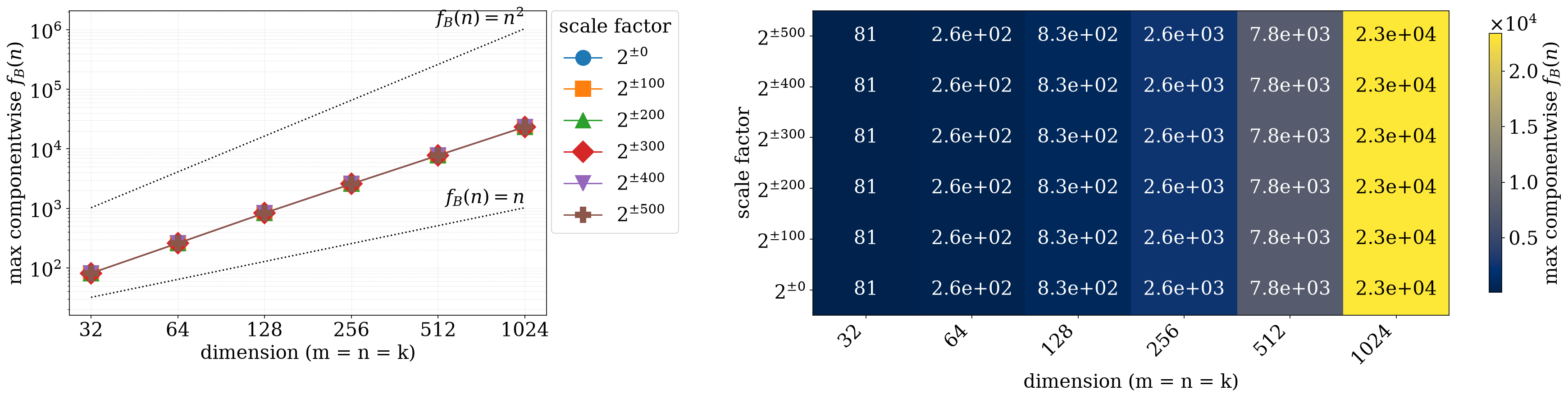}}}
        \vspace{-10mm}\caption{Test\_3 Sweep}\vspace{6mm}
    \end{subfigure}
    \begin{subfigure}{\textwidth}
        \centering
        \textcolor{lightgray}{\dashbox{\includegraphics[width=\textwidth]{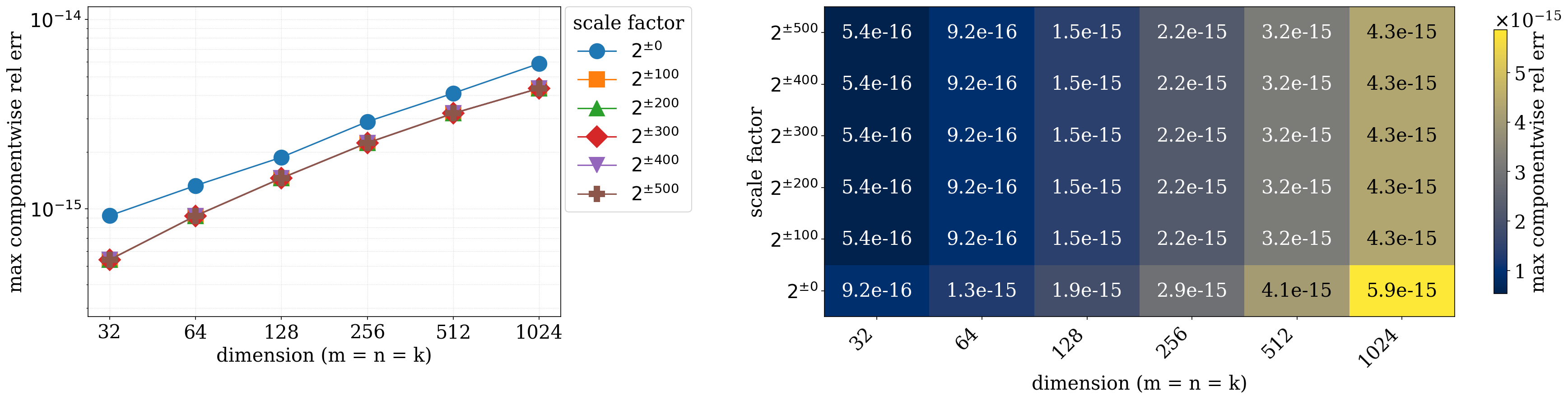}}}
        \vspace{-10mm}\caption{Test\_4b Sweep}\vspace{6mm}
    \end{subfigure}
    \caption{Netlib BLAS}
    \label{fig:sweeps_netlib_blas}
\end{figure}

\begin{figure}[p]
    \centering
    \begin{subfigure}{\textwidth}
        \centering
        \includegraphics[width=0.5\textwidth]{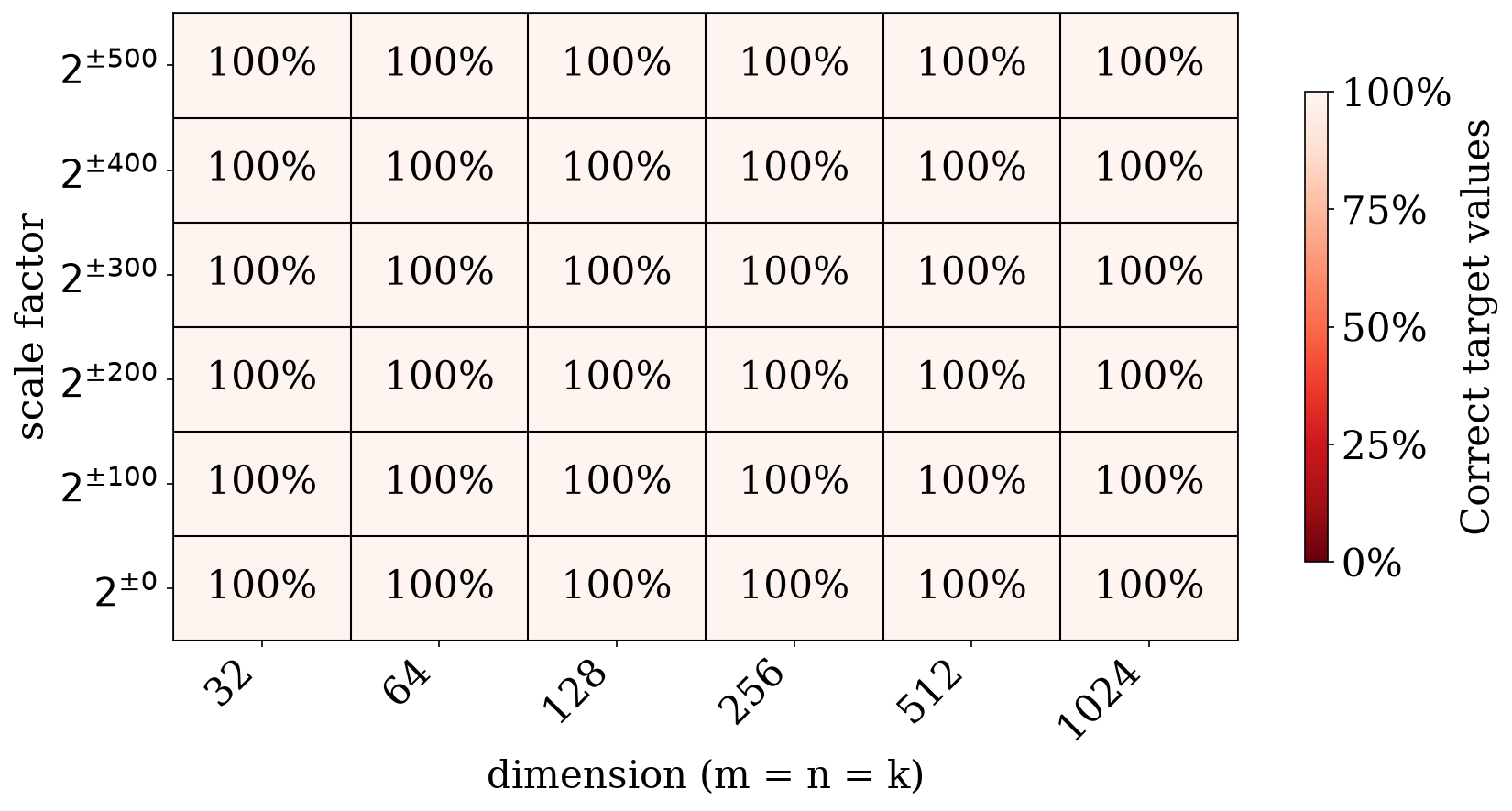}
        \caption{Test\_1c Sweep}\vspace{6mm}
    \end{subfigure}
    \begin{subfigure}{\textwidth}
        \centering
        \textcolor{lightgray}{\dashbox{\includegraphics[width=\textwidth]{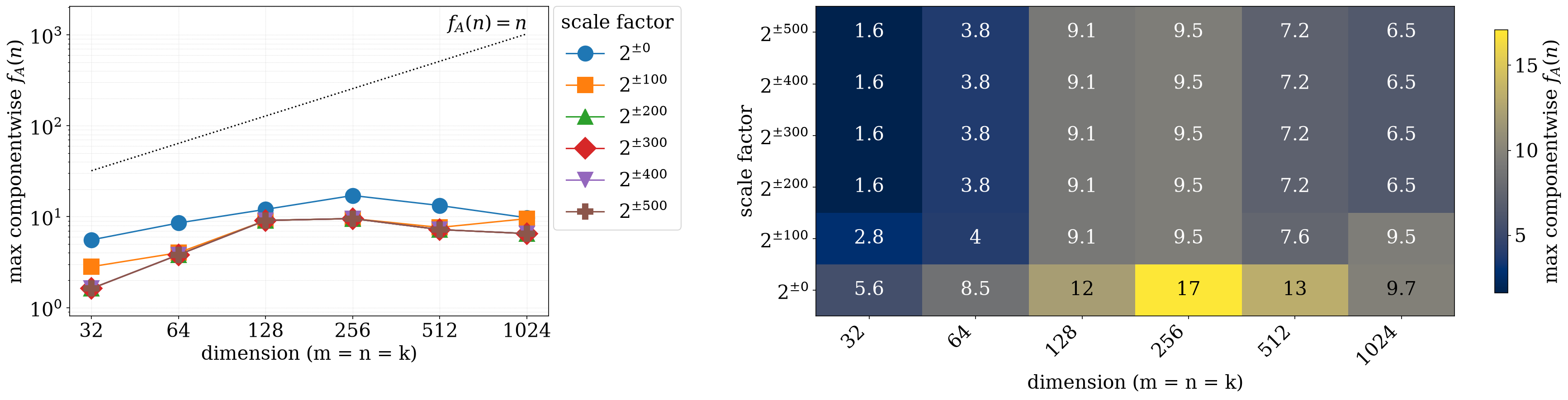}}}
        \vspace{-10mm}\caption{Test\_2c Sweep}\vspace{6mm}
    \end{subfigure}
    \begin{subfigure}{\textwidth}
        \centering
        \textcolor{lightgray}{\dashbox{\includegraphics[width=\textwidth]{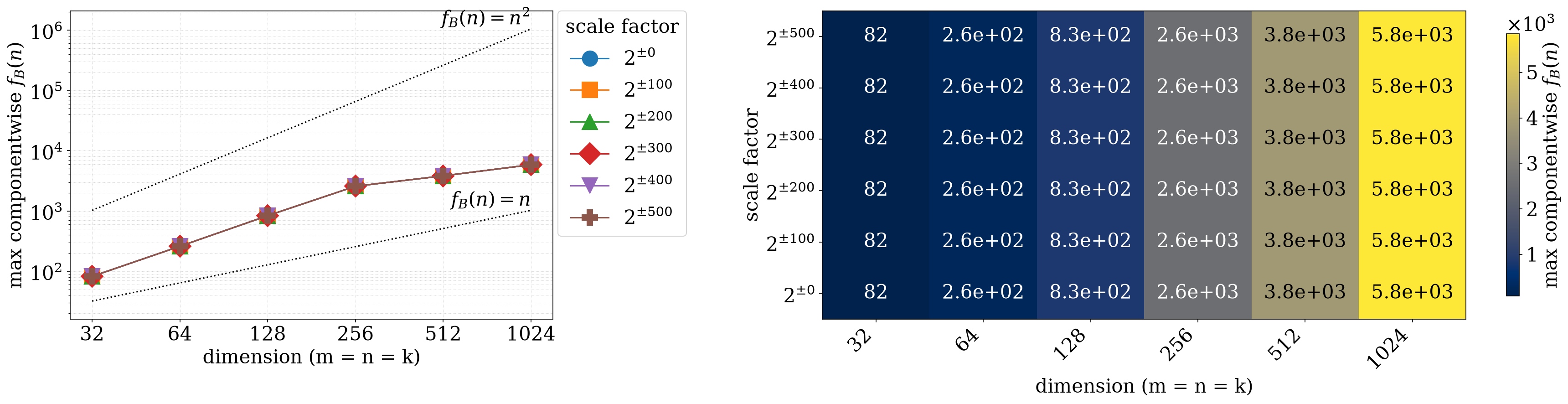}}}
        \vspace{-10mm}\caption{Test\_3 Sweep}\vspace{6mm}
    \end{subfigure}
    \begin{subfigure}{\textwidth}
        \centering
        \textcolor{lightgray}{\dashbox{\includegraphics[width=\textwidth]{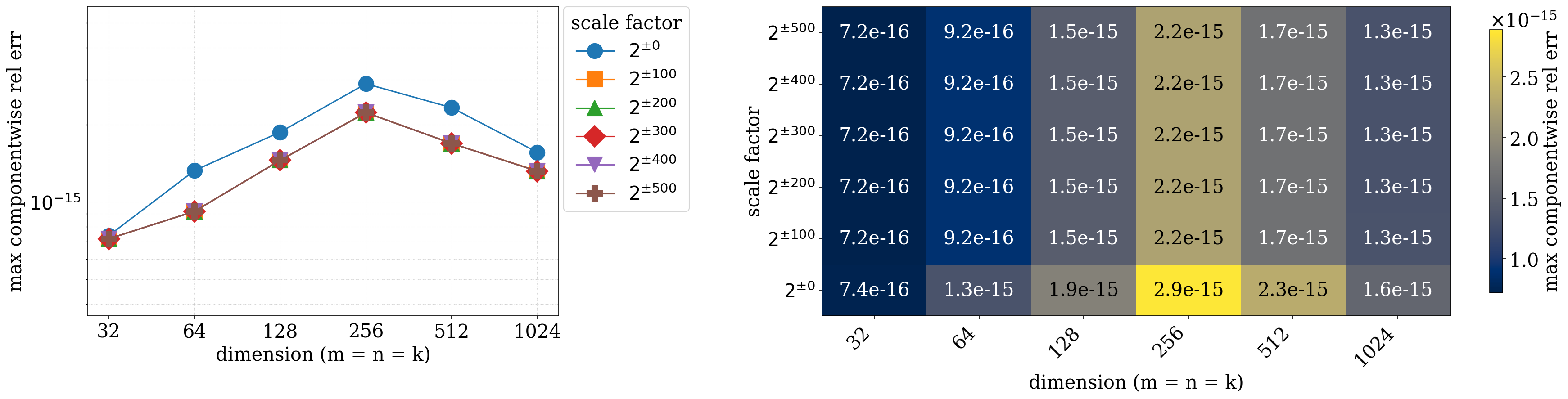}}}
        \vspace{-10mm}\caption{Test\_4b Sweep}\vspace{6mm}
    \end{subfigure}
    \caption{BLIS}
    \label{fig:sweeps_blis}
\end{figure}

\begin{figure}[p]
    \centering
    \begin{subfigure}{\textwidth}
        \centering
        \includegraphics[width=0.5\textwidth]{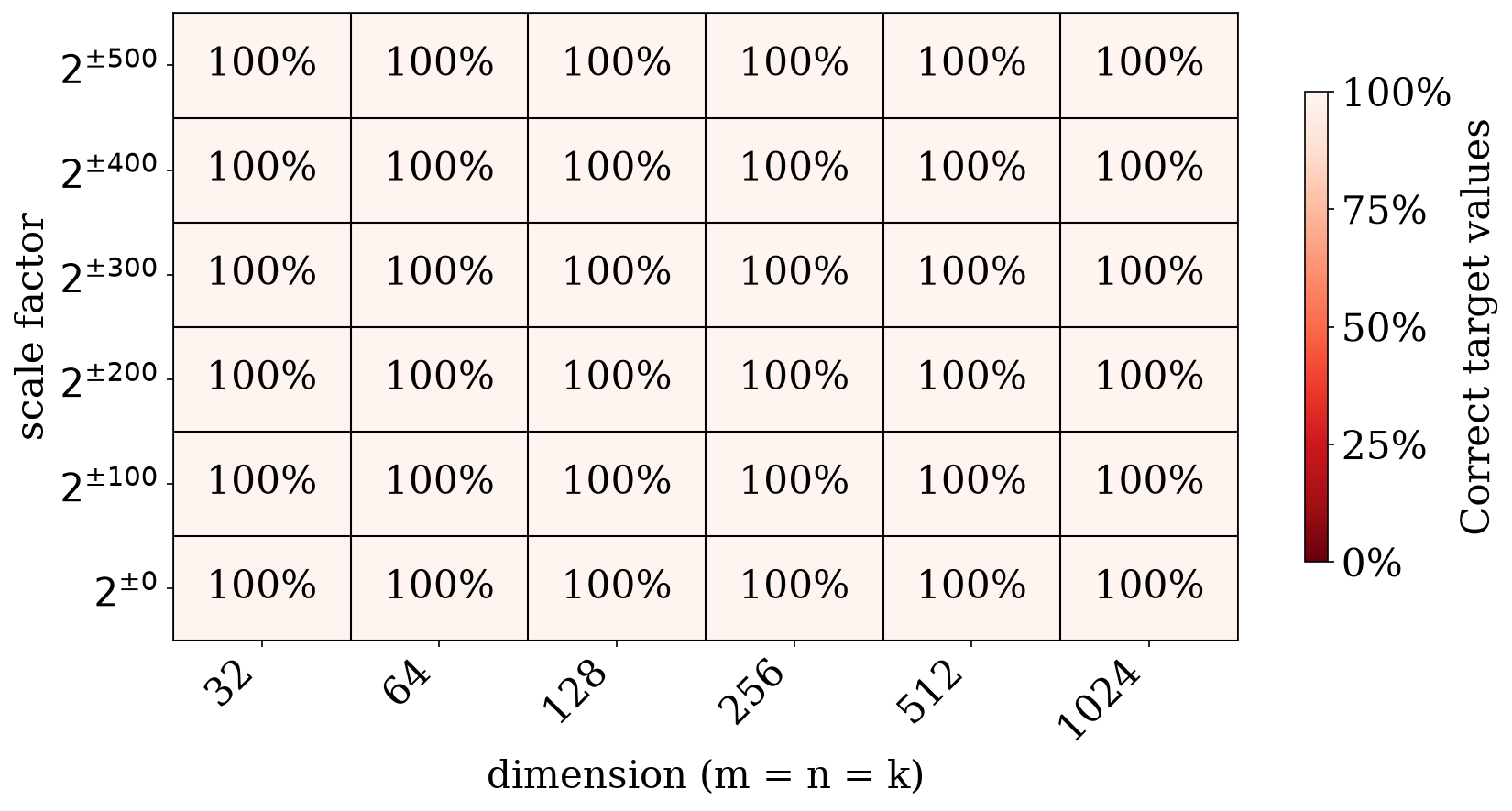}
        \caption{Test\_1c Sweep}\vspace{6mm}
    \end{subfigure}
    \begin{subfigure}{\textwidth}
        \centering
        \textcolor{lightgray}{\dashbox{\includegraphics[width=\textwidth]{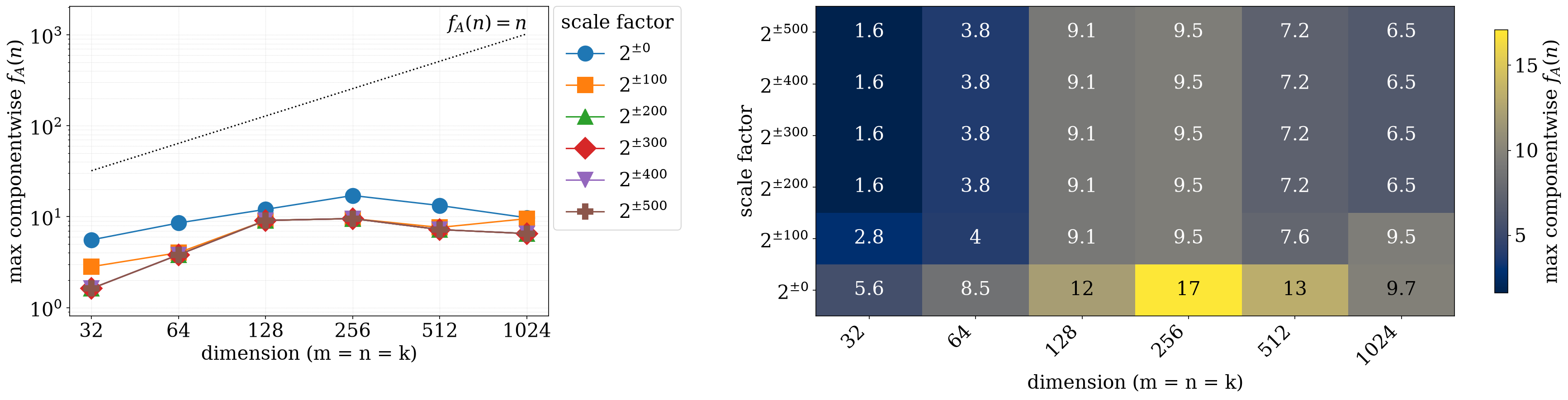}}}
        \vspace{-10mm}\caption{Test\_2c Sweep}\vspace{6mm}
    \end{subfigure}
    \begin{subfigure}{\textwidth}
        \centering
        \textcolor{lightgray}{\dashbox{\includegraphics[width=\textwidth]{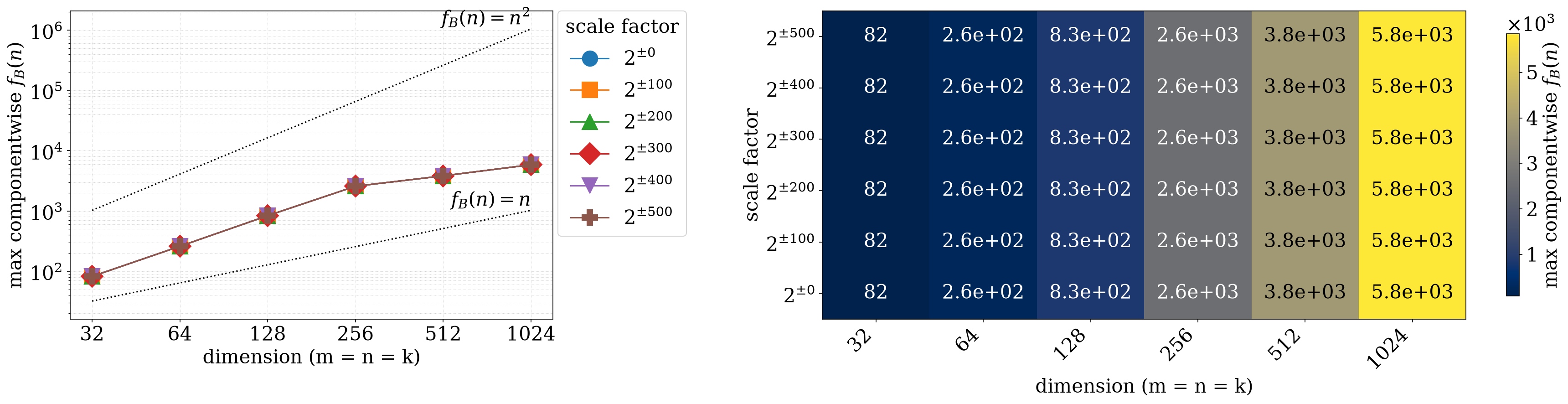}}}
        \vspace{-10mm}\caption{Test\_3 Sweep}\vspace{6mm}
    \end{subfigure}
    \begin{subfigure}{\textwidth}
        \centering
        \textcolor{lightgray}{\dashbox{\includegraphics[width=\textwidth]{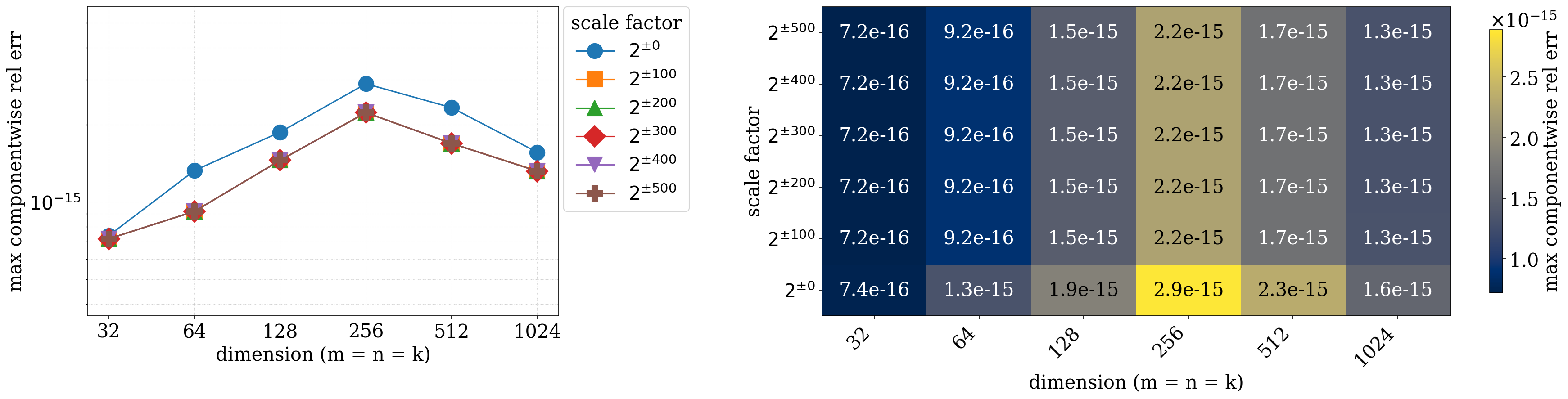}}}
        \vspace{-10mm}\caption{Test\_4b Sweep}\vspace{6mm}
    \end{subfigure}
    \caption{OpenBLAS}
    \label{fig:sweeps_openblas}
\end{figure}

\begin{figure}[p]
    \centering
    \begin{subfigure}{\textwidth}
        \centering
        \includegraphics[width=0.5\textwidth]{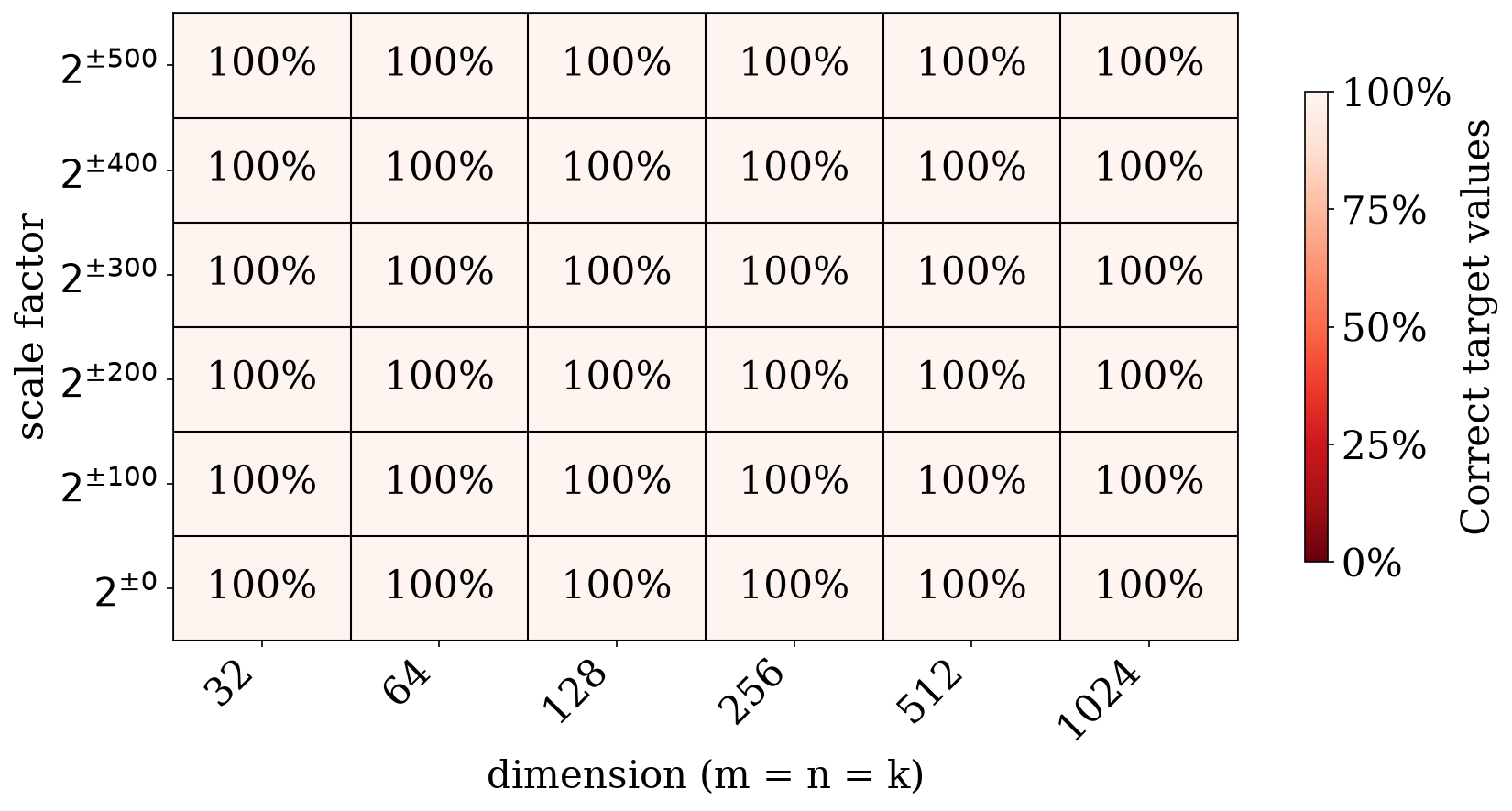}
        \caption{Test\_1c Sweep}\vspace{6mm}
    \end{subfigure}
    \begin{subfigure}{\textwidth}
        \centering
        \textcolor{lightgray}{\dashbox{\includegraphics[width=\textwidth]{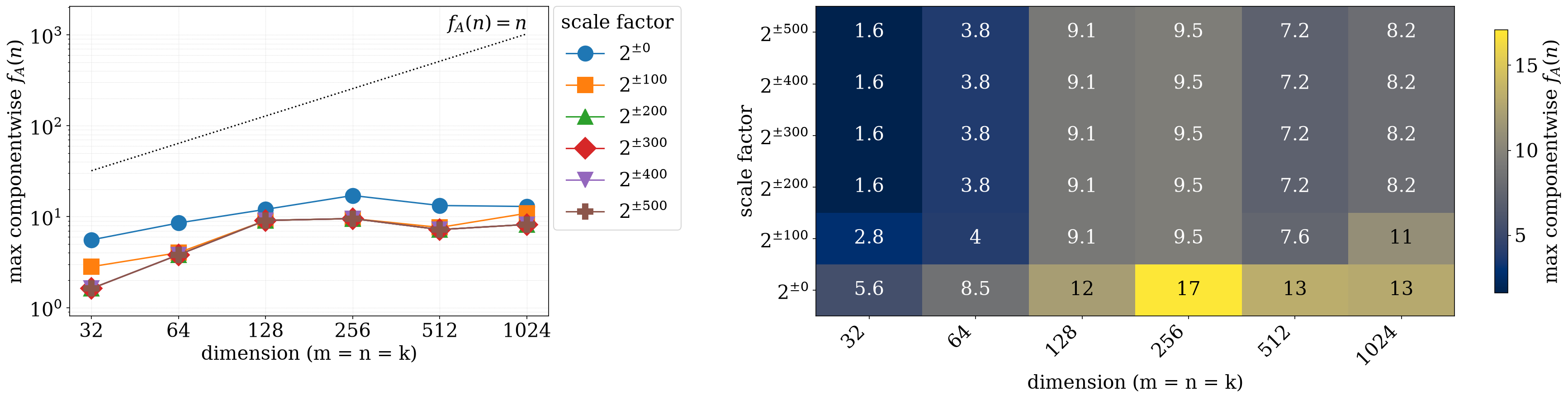}}}
        \vspace{-10mm}\caption{Test\_2c Sweep}\vspace{6mm}
    \end{subfigure}
    \begin{subfigure}{\textwidth}
        \centering
        \textcolor{lightgray}{\dashbox{\includegraphics[width=\textwidth]{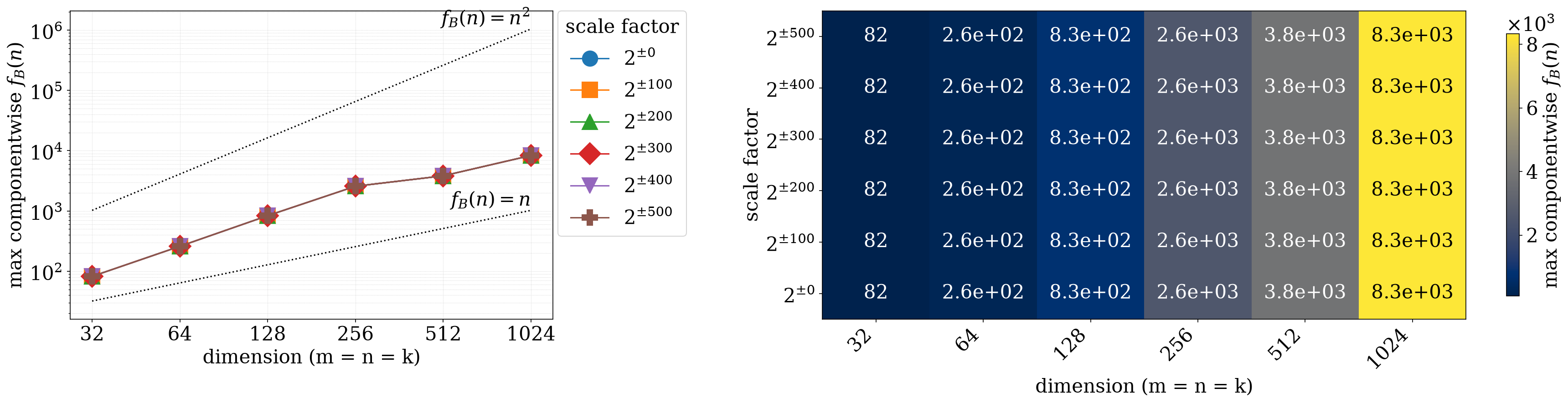}}}
        \vspace{-10mm}\caption{Test\_3 Sweep}\vspace{6mm}
    \end{subfigure}
    \begin{subfigure}{\textwidth}
        \centering
        \textcolor{lightgray}{\dashbox{\includegraphics[width=\textwidth]{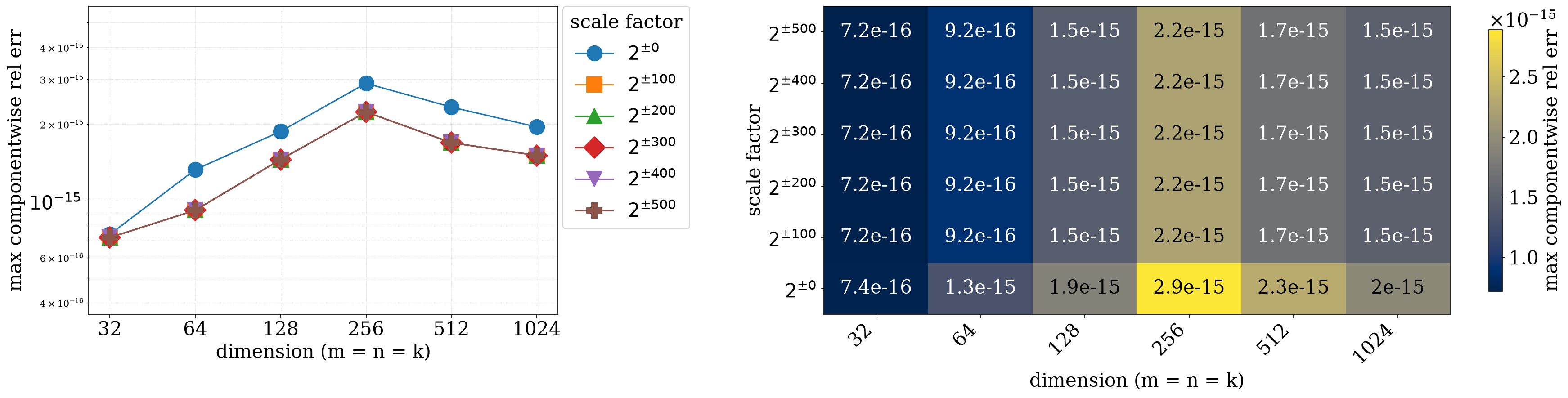}}}
        \vspace{-10mm}\caption{Test\_4b Sweep}\vspace{6mm}
    \end{subfigure}
    \caption{oneMKL 2026.0 sequential and parallel multi-threaded (OpenMP 4.5, GNU pthread 2.39, 32 threads, single oneMKL thread per core). The plots represent both sequential and parallel multi-threaded oneMKL implementations due to identical experimental data.}
    \label{fig:sweeps_onemkl_openmp_gnu_thread}
\end{figure}

\begin{figure}[p]
    \centering
    \begin{subfigure}{\textwidth}
        \centering
        \includegraphics[width=0.5\textwidth]{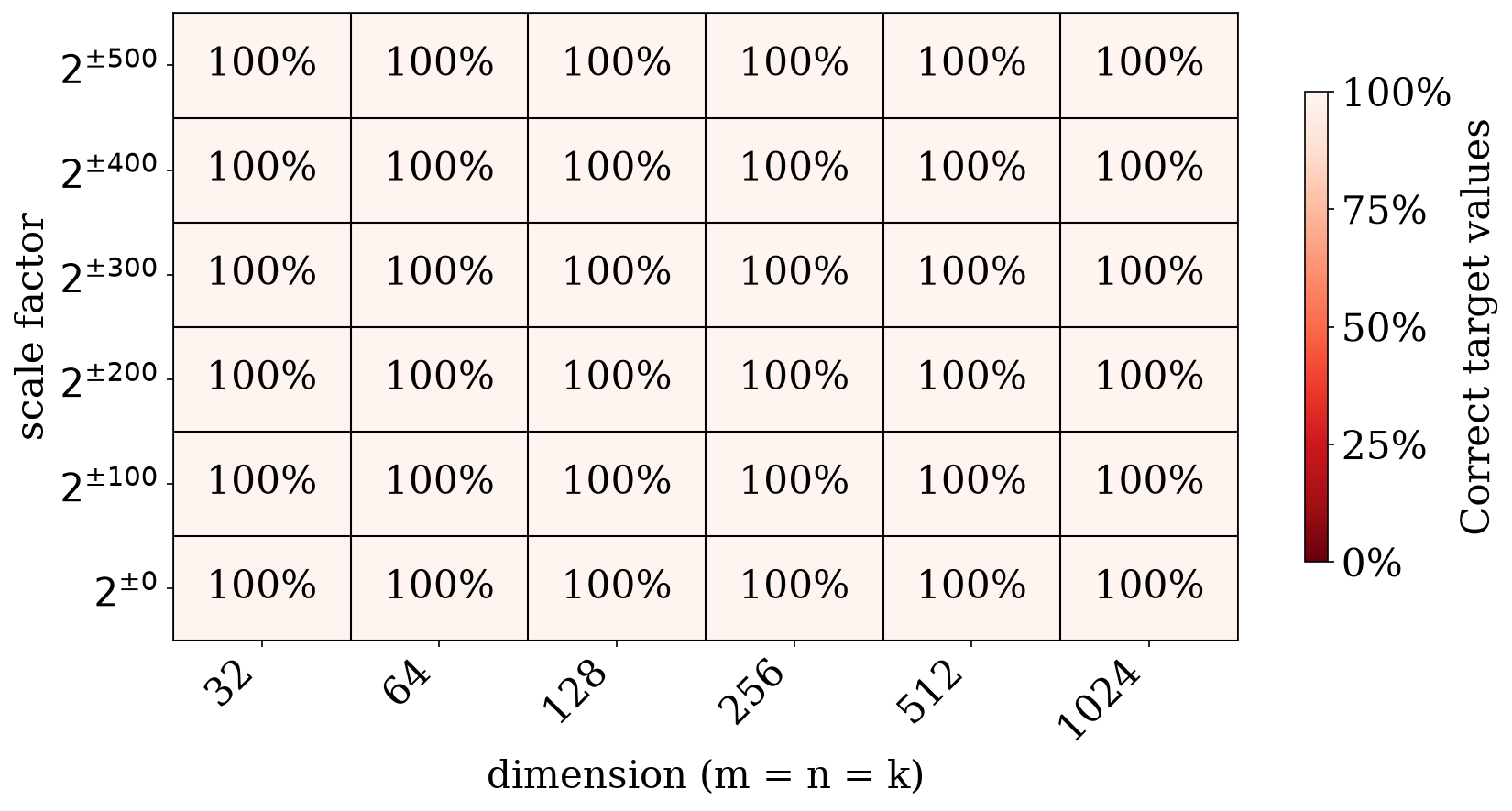}
        \caption{Test\_1c Sweep}\vspace{6mm}
    \end{subfigure}
    \begin{subfigure}{\textwidth}
        \centering
        \textcolor{lightgray}{\dashbox{\includegraphics[width=\textwidth]{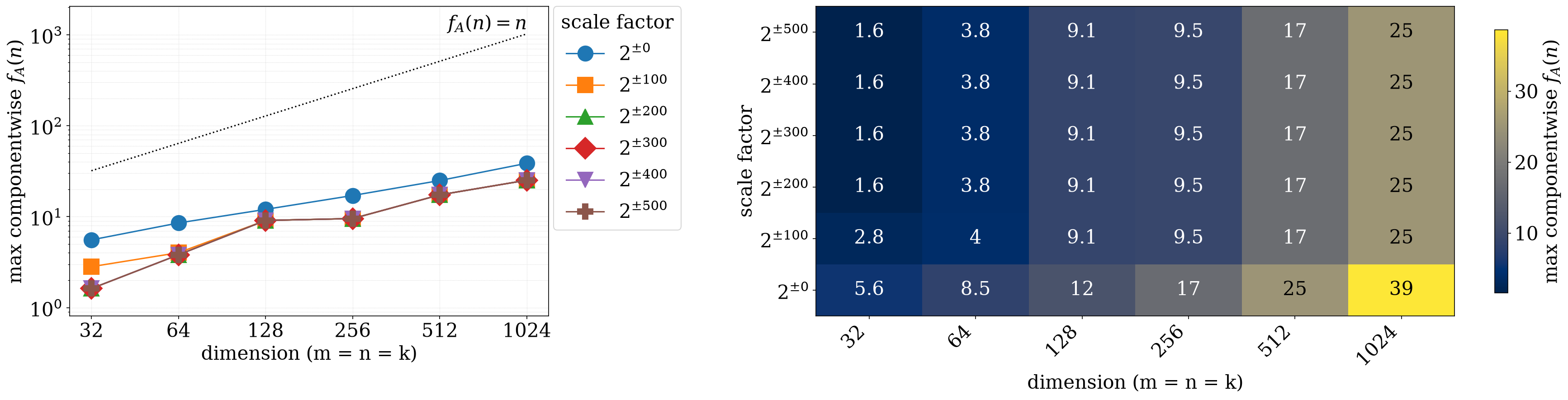}}}
        \vspace{-10mm}\caption{Test\_2c Sweep}\vspace{6mm}
    \end{subfigure}
    \begin{subfigure}{\textwidth}
        \centering
        \textcolor{lightgray}{\dashbox{\includegraphics[width=\textwidth]{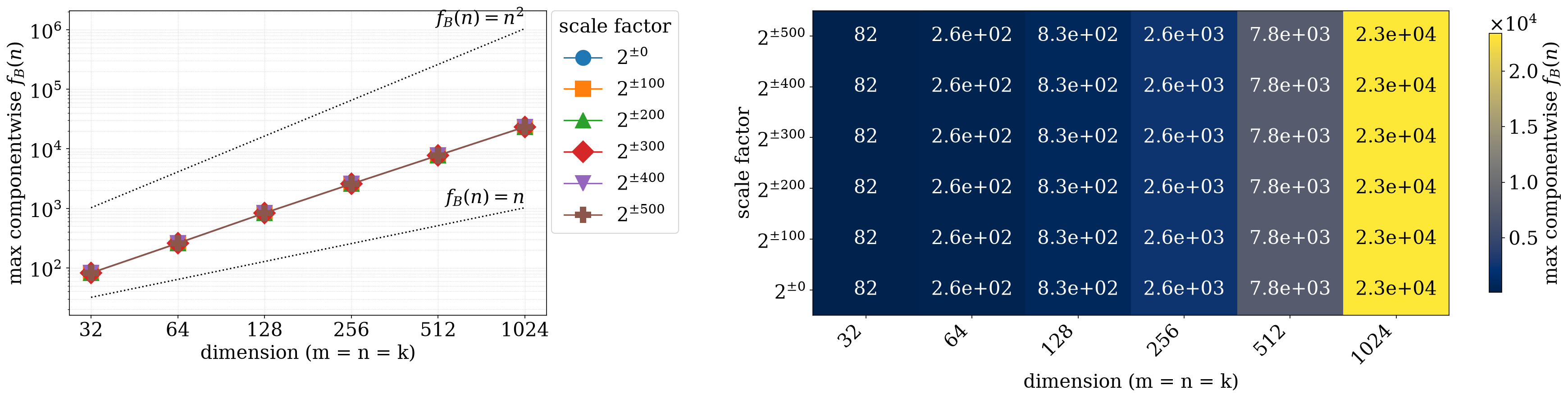}}}
        \vspace{-10mm}\caption{Test\_3 Sweep
        }\vspace{6mm}
    \end{subfigure}
    \begin{subfigure}{\textwidth}
        \centering
        \textcolor{lightgray}{\dashbox{\includegraphics[width=\textwidth]{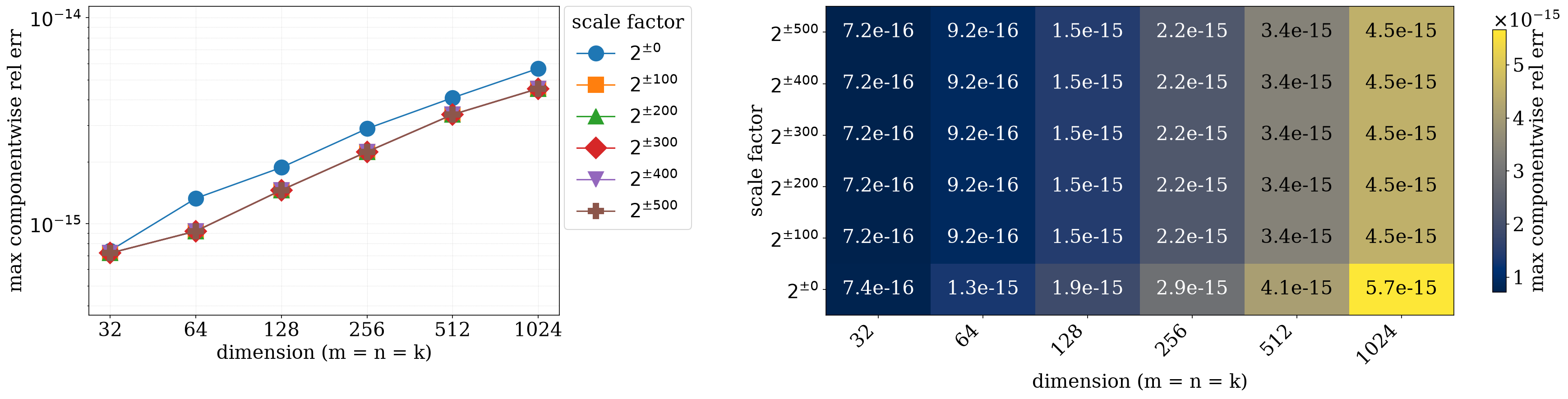}}}
        \vspace{-10mm}\caption{Test\_4b Sweep}\vspace{6mm}
    \end{subfigure}
    \caption{cuBLAS}
    \label{fig:sweeps_cublas}
\end{figure}

\begin{figure}[p]
    \centering
    \begin{subfigure}[t]{0.48\textwidth}
        \centering
        \includegraphics[width=\textwidth]{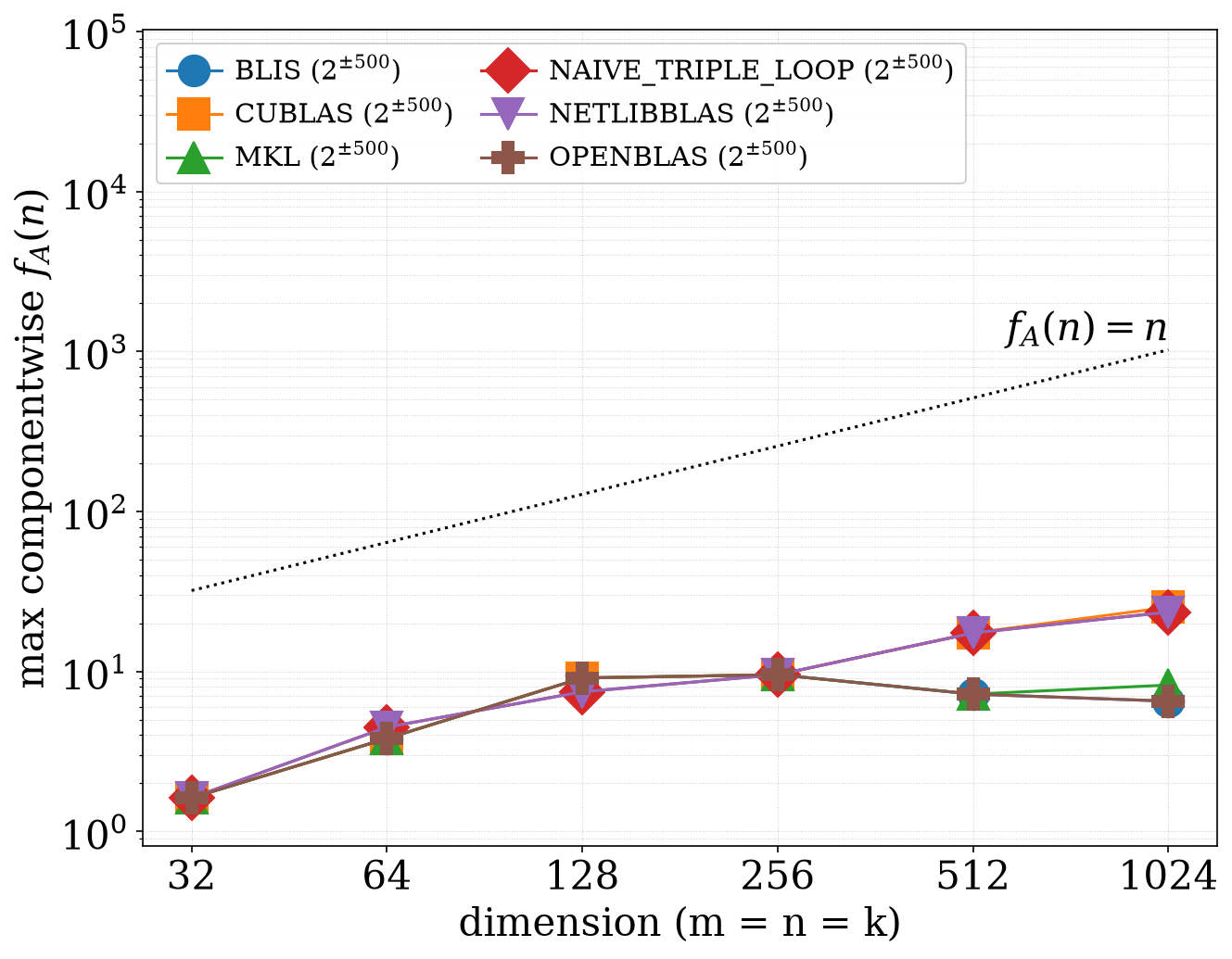}
        \caption{Test\_2c}
        \label{fig:aggregate_f_of_n_test2c}
    \end{subfigure}\hfill
    \begin{subfigure}[t]{0.48\textwidth}
        \centering  
        \includegraphics[width=\textwidth]{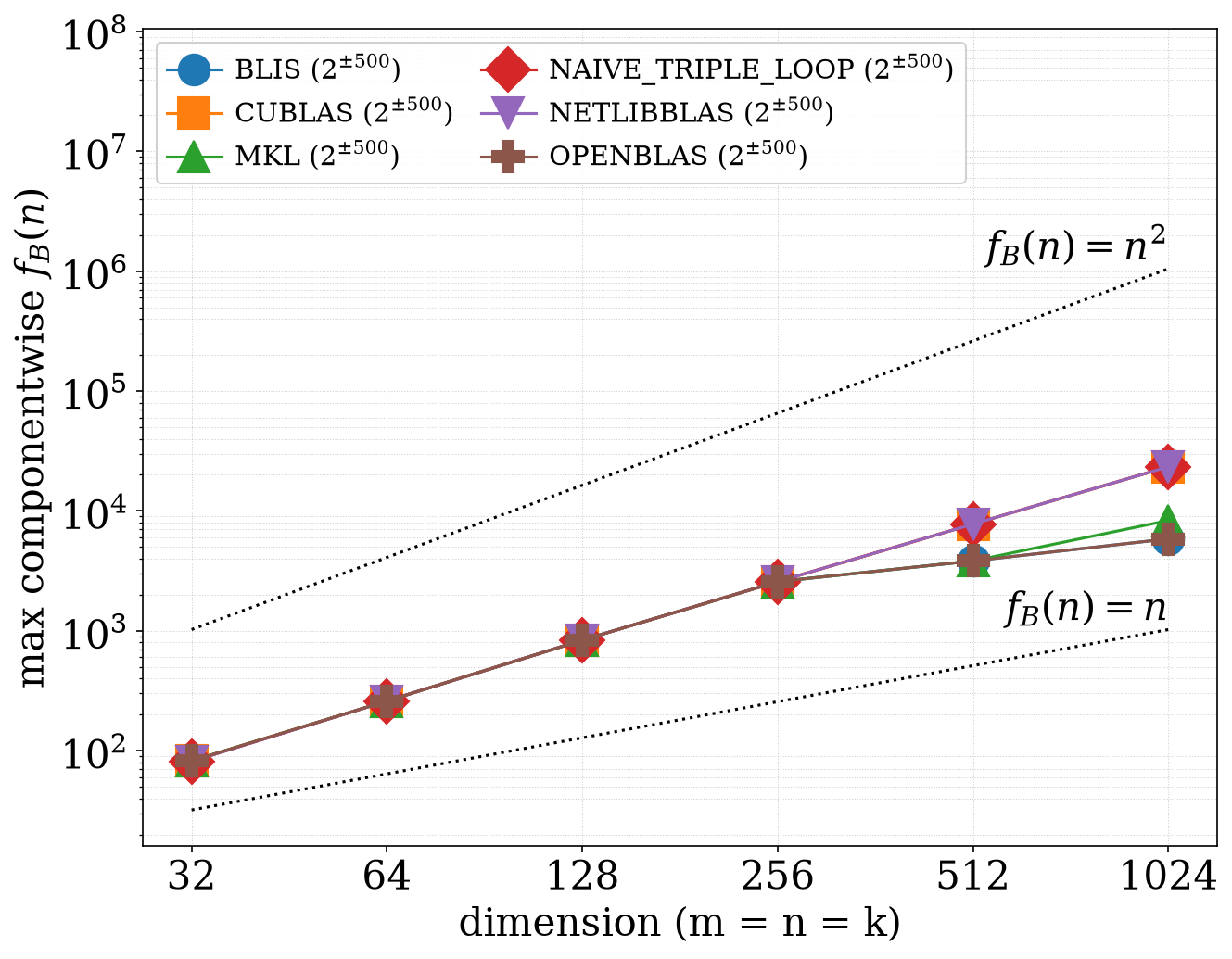}
        \caption{Test\_3}
        \label{fig:aggregate_f_of_n_test3}
    \end{subfigure}
    \caption{Test\_2c and Test\_3 results across production implementations, across a range of $n$, and for the largest scaling factor $2^{\pm 500}$.}
    \label{fig:aggregate_f_of_n}
\end{figure}

\begin{figure}[p]
    \centering
    \begin{subfigure}{\textwidth}
        \centering
        \includegraphics[width=0.5\textwidth]{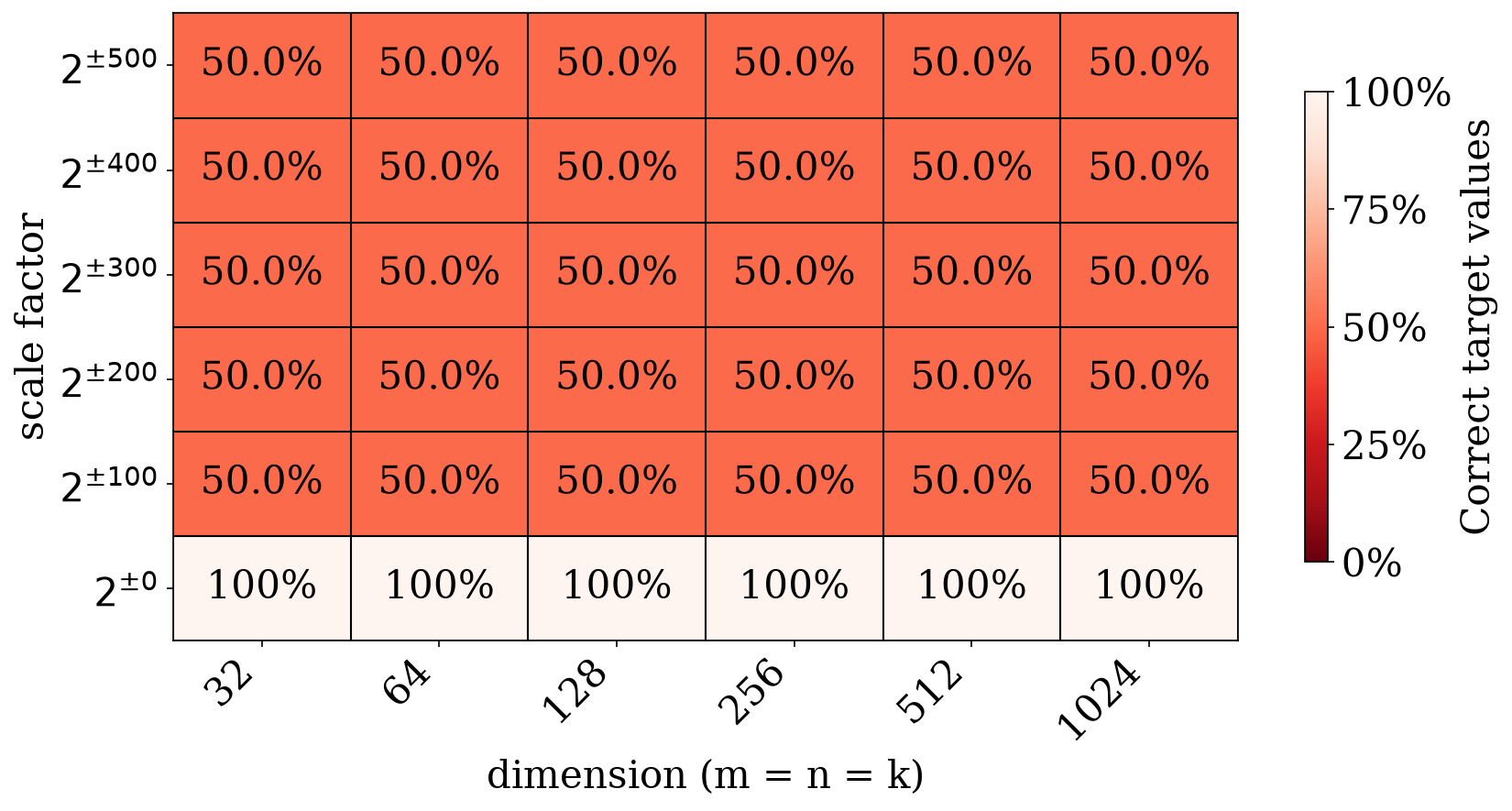}
        \caption{Test\_1c Sweep}\vspace{6mm}
    \end{subfigure}
    \begin{subfigure}{\textwidth}
        \centering
        \textcolor{lightgray}{\dashbox{\includegraphics[width=\textwidth]{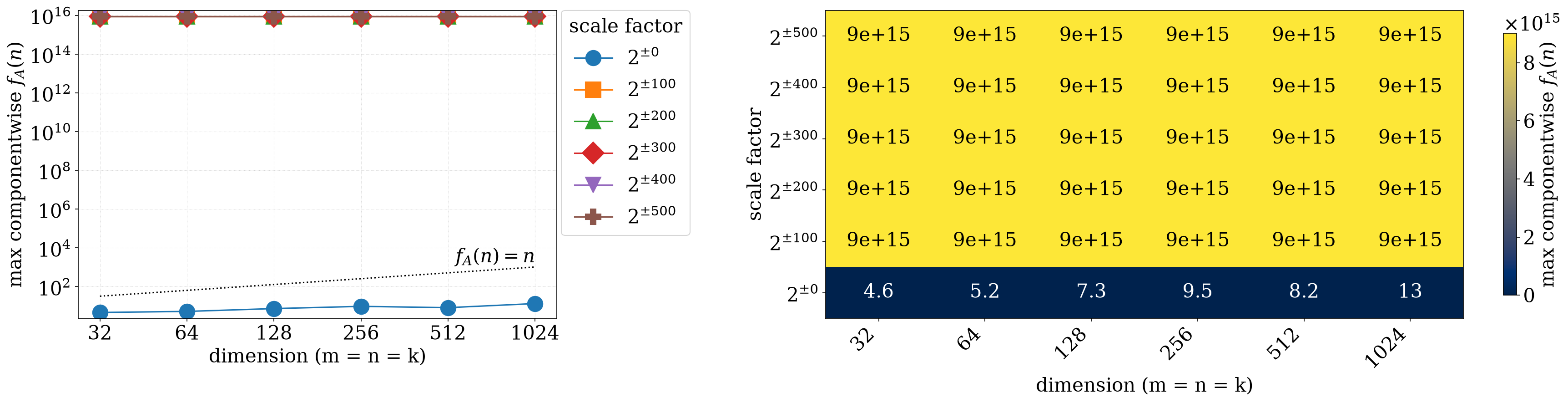}}}
        \vspace{-10mm}\caption{Test\_2c Sweep}\vspace{6mm}
    \end{subfigure}
    \begin{subfigure}{\textwidth}
        \centering
        \textcolor{lightgray}{\dashbox{\includegraphics[width=\textwidth]{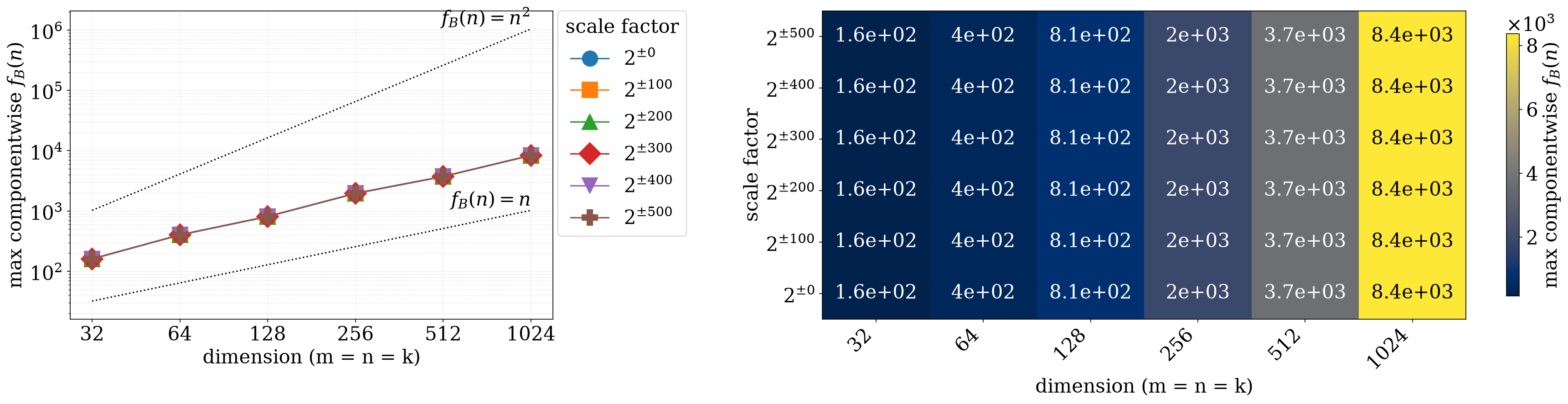}}}
        \vspace{-10mm}\caption{Test\_3 Sweep}\vspace{6mm}
    \end{subfigure}
    \begin{subfigure}{\textwidth}
        \centering
        \textcolor{lightgray}{\dashbox{\includegraphics[width=\textwidth]{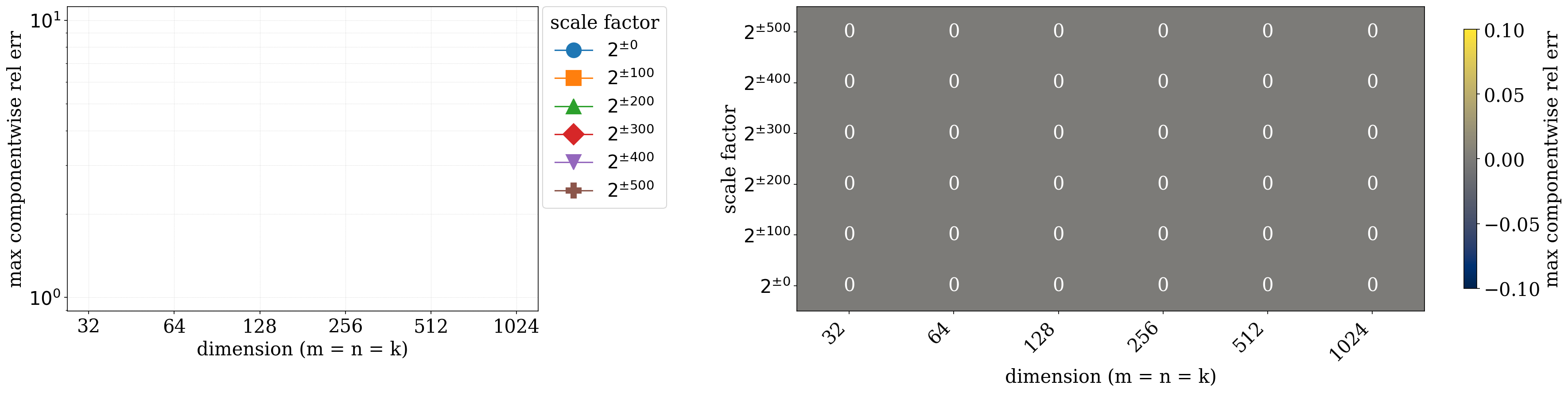}}}
        \vspace{-10mm}\caption{Test\_4b Sweep}\vspace{6mm}
    \end{subfigure}
    \caption{Ozaki I}
    \label{fig:sweeps_ozaki1}
\end{figure}

\begin{figure}[p]
    \centering
    \begin{subfigure}{\textwidth}
        \centering
        \includegraphics[width=0.5\textwidth]{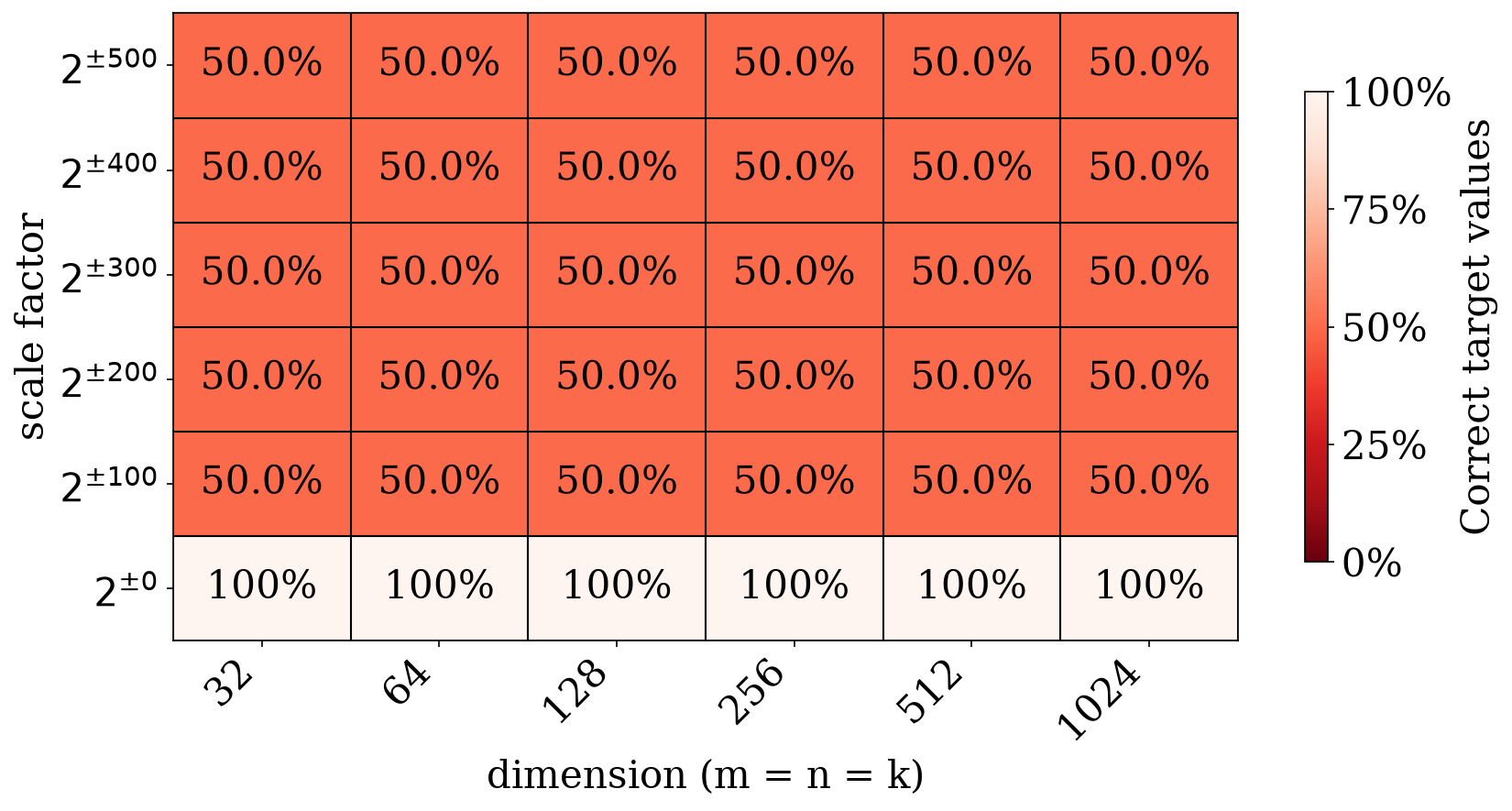}
        \caption{Test\_1c Sweep}\vspace{6mm}
    \end{subfigure}
    \begin{subfigure}{\textwidth}
        \centering
        \textcolor{lightgray}{\dashbox{\includegraphics[width=\textwidth]{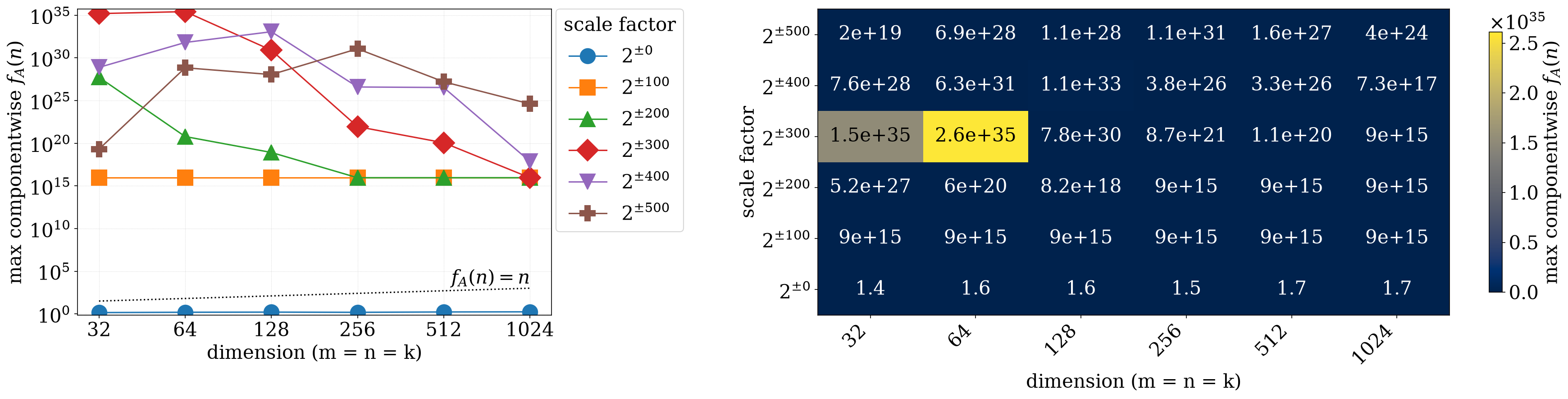}}}
        \vspace{-10mm}\caption{Test\_2c Sweep}\vspace{6mm}
    \end{subfigure}
    \begin{subfigure}{\textwidth}
        \centering
        \textcolor{lightgray}{\dashbox{\includegraphics[width=\textwidth]{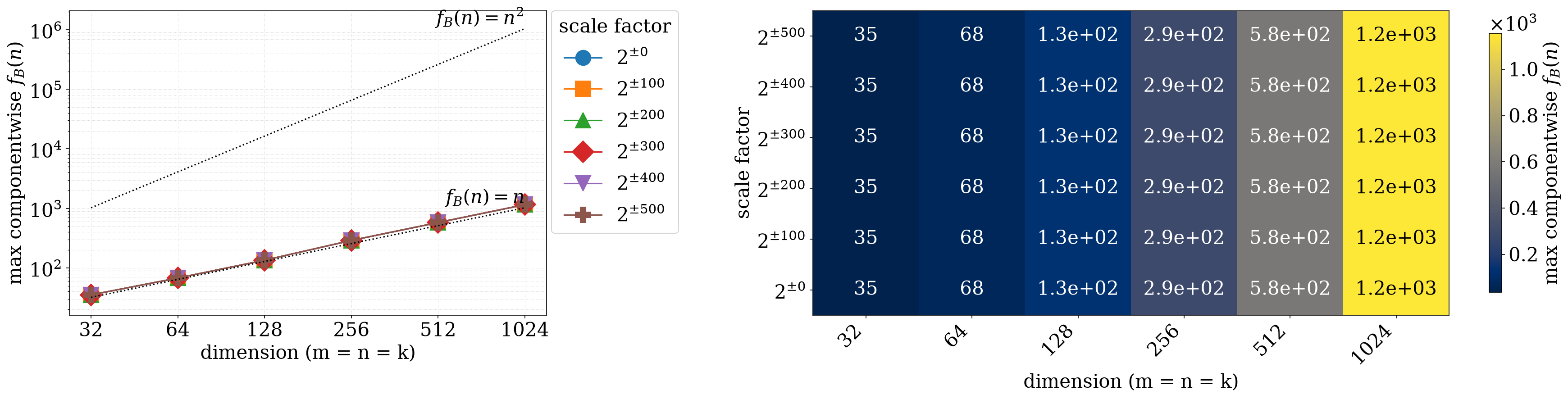}}}
        \vspace{-10mm}\caption{Test\_3 Sweep}\vspace{6mm}
    \end{subfigure}
    \begin{subfigure}{\textwidth}
        \centering
        \textcolor{lightgray}{\dashbox{\includegraphics[width=\textwidth]{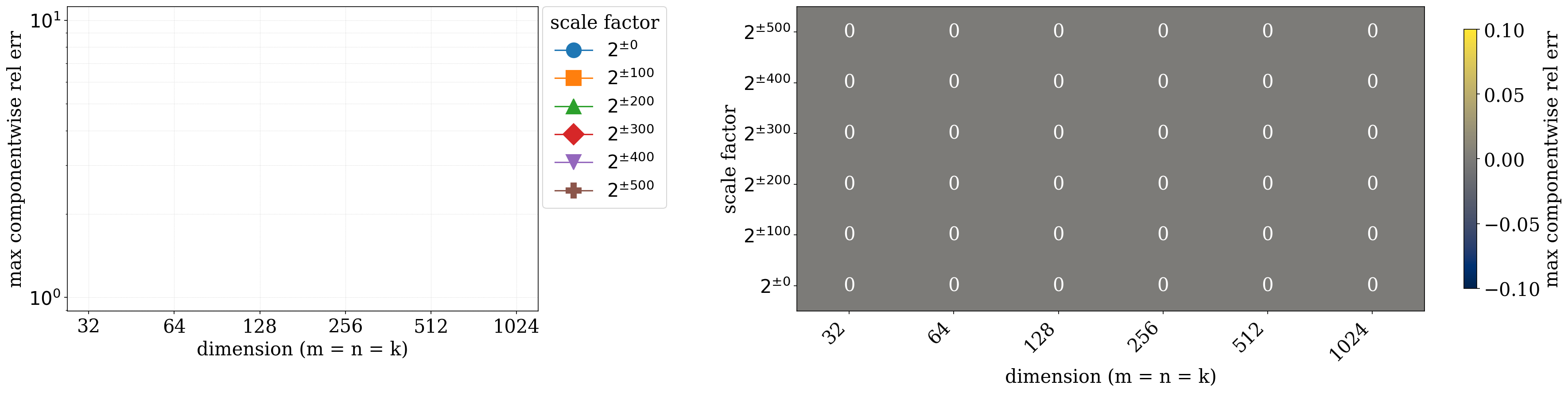}}}
        \vspace{-10mm}\caption{Test\_4b Sweep}\vspace{6mm}
    \end{subfigure}
    \caption{Ozaki II}
    \label{fig:sweeps_ozaki2}
\end{figure}

\begin{figure}[p]
    \centering
    \begin{subfigure}{\textwidth}
        \centering
        \includegraphics[width=0.5\textwidth]{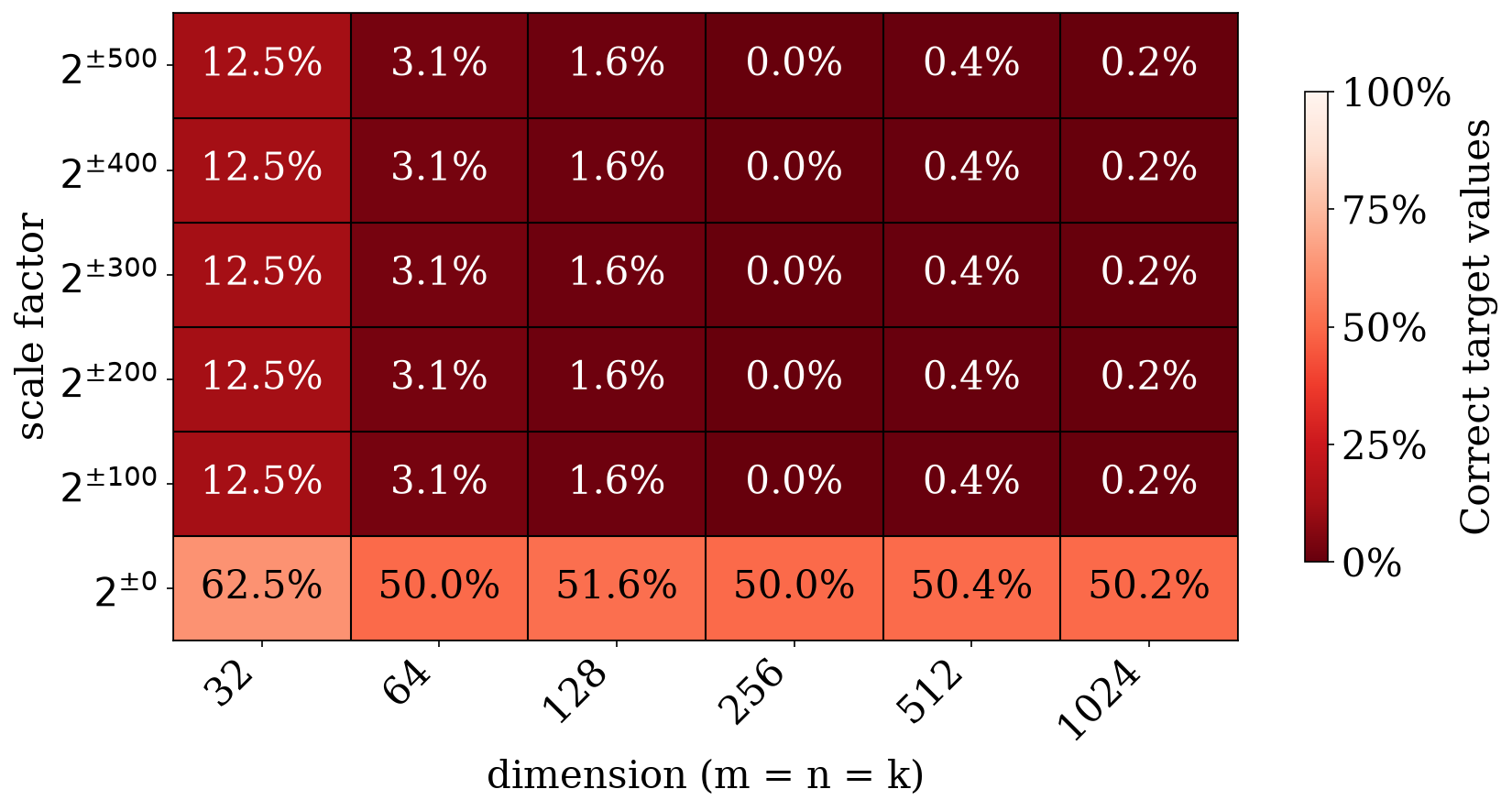}
        \caption{Test\_1c Sweep}\vspace{6mm}
    \end{subfigure}
    \begin{subfigure}{\textwidth}
        \centering
        \textcolor{lightgray}{\dashbox{\includegraphics[width=\textwidth]{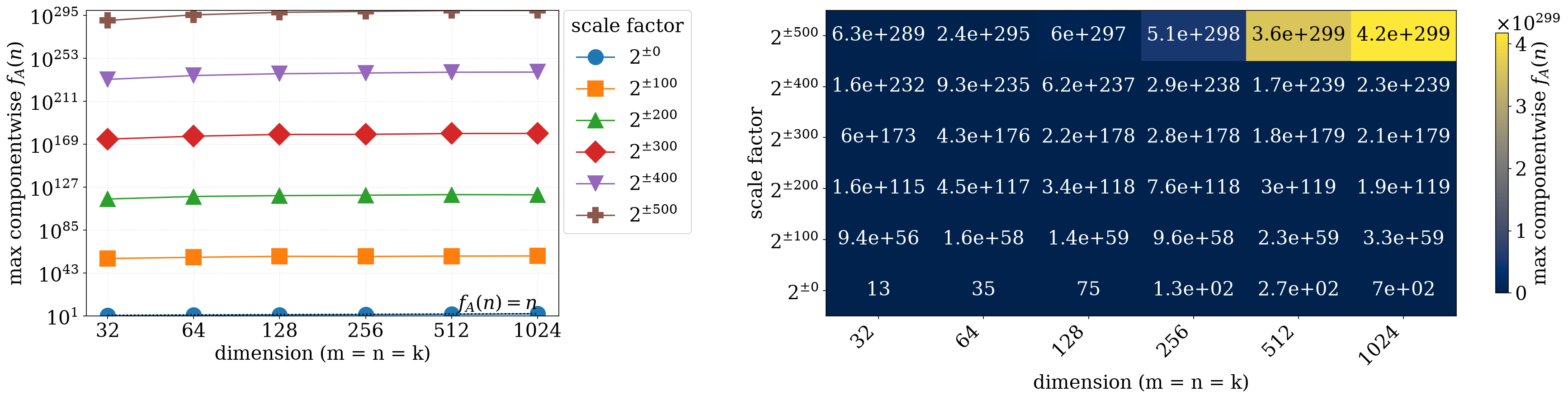}}}  
        \vspace{-10mm}\caption{Test\_2c Sweep}\vspace{6mm}
    \end{subfigure}
    \begin{subfigure}{\textwidth}
        \centering
        \textcolor{lightgray}{\dashbox{\includegraphics[width=\textwidth]{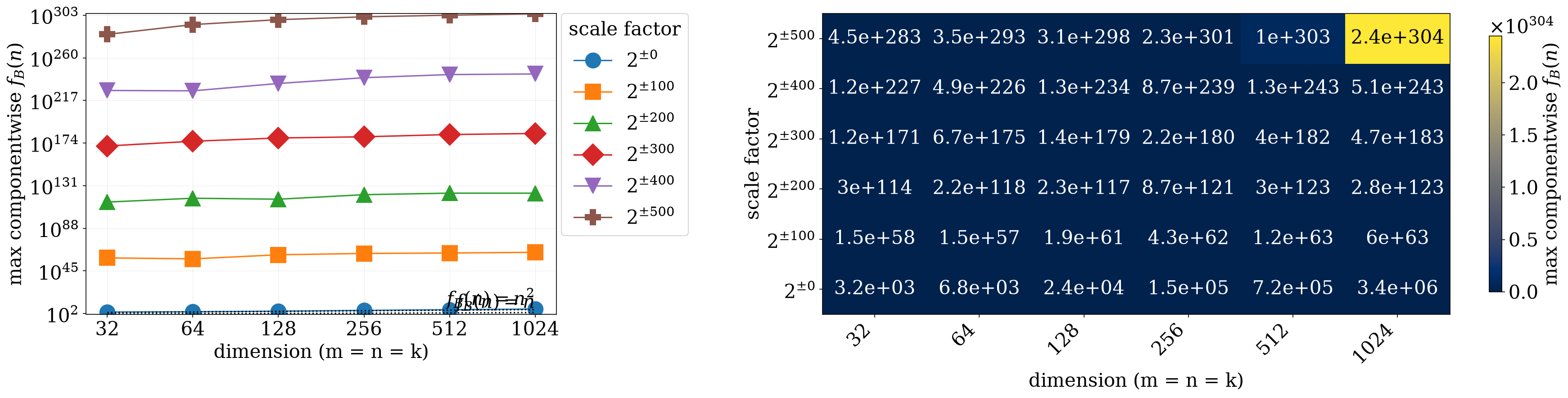}}}
        \vspace{-10mm}\caption{Test\_3 Sweep}\vspace{6mm}
    \end{subfigure}
    \begin{subfigure}{\textwidth}
        \centering
        \textcolor{lightgray}{\dashbox{\includegraphics[width=\textwidth]{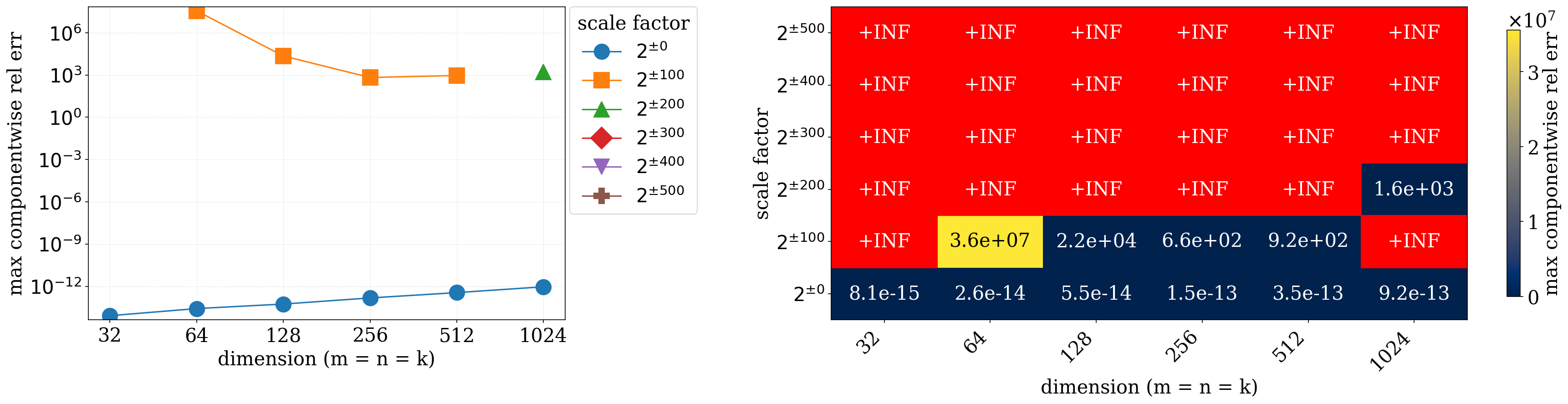}}}
        \vspace{-10mm}\caption{Test\_4b Sweep}\vspace{6mm}
    \end{subfigure}
    \caption{Strassen}
    \label{fig:sweeps_strassen}
\end{figure}

For this demonstration, we target Netlib BLAS (3.12.1), BLIS (2.0), OpenBLAS (0.3.30), oneMKL (2026.0, both sequential and multi-threaded, and run on a Dual Socket Intel Xeon Gold 6444Y (16 cores/32 threads per socket) server with a single oneMKL thread per core) and cuBLAS (from CUDA Toolkit 13.3 and run on a single NVIDIA Blackwell B200 GPU).
For reference, we also target the Naive Triple Loop implementation described in Section \ref{sec:testresultspart1}. As we will see below, all of these earn an ``A'' grade  across all the tested dimensions and scaling factors. In order to demonstrate what test results for non-``A'' DGEMM implementations might look like, we also target the less-robust research implementations of Ozaki I, Ozaki II, and Strassen described in Section \ref{sec:testresultspart1}. Results are presented in Figures \ref{fig:sweeps_naive_triple_loop}-\ref{fig:sweeps_strassen}.

\paragraph{First, we consider the production DGEMM implementations.} Our grading code reports ``A'' for all of Netlib BLAS, BLIS, OpenBLAS, oneMKL and cuBLAS. 
This indicates that ``A''-grade accuracy was observed across all the dimensions and scaling factors tested. This is depicted in the results of Test\_1c and Test\_2c. In the case of the former, look to the uniform 100\% values in their (a) subfigures. In the case of the latter, look to the subunit slope of all lines depicted in the left plots of the (b) subfigures. For a view of all production DGEMM implementations on the same axis (evaluated at their largest scaling factor of $2^{\pm 500}$ and across problem sizes), see Figure \ref{fig:aggregate_f_of_n_test2c} which plainly shows all implementations safely below the threshold of $f_A(n) < n$. Furthermore, this figure highlights a pattern in BLIS/OpenBLAS/oneMKL which have a distinct elbow at dimension 256 which the others do not. This is due to K-blocking, for which the sweep of the dimension parameters covers the threshold at which blocking occurs in these implementations. On the other hand, the naive triple loop and Netlib BLAS have no K-blocking and the threshold for cuBLAS is likely higher than the dimensions tested here. 

While no further testing beyond this is necessary to derive a DGEMM grade, sweeps of Test\_3 and Test\_4b are included for completeness and discussed below.

One of the key signals from Test\_3 is the invariance with respect to scaling factor which fits expectation for dot-product-based DGEMM implementations. Note that the line $f_B(n) = n^2$ in the (c) subfigures is the failure threshold and $f_B(n) = n$ is included for reference only. For a view of all production DGEMM implementations on the same axis (evaluated at their largest scaling factor of $2^{\pm 500}$ and across problem sizes), see Figure \ref{fig:aggregate_f_of_n_test3} which plainly shows all implementations safely below the threshold of $f_B(n) < n^2$ as well as the effects of K-blocking in BLIS/OpenBLAS/oneMKL.

The sweeps of Test\_4b depicted in the (d) subfigures show the K-blocking effects more starkly while also offering an additional insight: Recall that Test\_4b involves calculating relative errors induced by permutation of dot product summands; here, we observe relative errors that grow as $\sqrt{N}$ rather than all zero or catastrophic as we might expect with emulation and Strassen respectively, thus corroborating the ``A'' grade determination from Test\_1c and Test\_2c \cite{HighamMary2019}.

\paragraph{Next we consider the less robust research implementations of the emulation and Strassen approaches described in Section \ref{sec:testresultspart1} for comparison.}
Our grading code reports ``B'' for Ozaki I and Ozaki II and ``C'' for Strassen. We begin by looking at the results of Test\_1c and Test\_2c which determine an ``A'' grade. For all tested cases with non-trivial scaling factors, these implementations do not meet the bar: the percentage of correct target-values in Test\_1c dips well below 100\% and $f_A(n)$ soars far above $n$ in Test\_2c. For the trivial unit scaling factor, the picture is slightly different. In this regime, both Ozaki schemes earn an ``A''; this fits expectations as these emulation methods struggle with representing inputs with a wide dynamic range. Subsequent testing using Test\_1c or Test\_2c that zooms in between the scaling factors of $2^0$ (pass) and $2^{\pm 100}$ (fail) would shed more light on the boundary condition on the dynamic range to maintain ``A''-grade accuracy. On the other hand, results for Strassen with trivial unit scaling factor are mixed: the percentage of correct target values in Test\_1c drops below 100\% indicating a non-``A'' grade, but the values of $f_A(n)$ measured in Test\_2c come in under the $f_A(n) < n$ threshold for an ``A''. Still, it should be noted that, while under the threshold, the magnitude of these values is still much greater than those observed in the aforementioned production DGEMM implementations. In the face of such mixed results, our grading code will not assign an ``A''.

Subsequently, Test\_3 marks the decision point between ``B'' and ``C'' grades. As was seen in Section \ref{sec:testresultspart1}, there is a clear discrepancy between the emulation-based DGEMM implementations and the Strassen-based DGEMM implementation with respect to the magnitudes of $f_B(n)$. Note too that Ozaki I and Ozaki II both maintain invariance with respect to scaling factors expected from a dot-product-based DGEMM while Strassen's $f_B(n)$ varies directly. Accordingly, Ozaki I and Ozaki II earn ``B'' grades and Strassen gets a ``C''.

As before with the ``A'' DGEMM's, the results of Test\_4b corroborate the assigned grades. The ``B'' grade emulation methods incur zero error under permutation of dot product summands because the dot product is exact whereas the ``C'' grade Strassen method incurs catastrophic errors.

\ignore{
\paragraph{K-blocking is the likely explanation for the decreasing error in BLIS and OpenBLAS for DGEMMs of size greater than 256.} \jackson{todo: can likely verify this by inspecting the source}

\paragraph{For the trivial unit scaling factor, Strassen earns an ``A'' on Test\_2c but does not on Test\_1c}. \jackson{I could use some help interpreting this}

\jwd{In the case of a unit scaling factor,
Test\_2a and Test\_2b set all entries of
$A$ and $B$ to random numbers in $[1,2]$,
for which all algorithms compute the entries
of $C$ to high relative accuracy. For Test\_2c,
the small normwise error bound for Strassen is
also enough to compute the entries of $C$
which are sums of positive numbers to high
relative accuracy. In Test\_1c, Strassen 
(without Zeroing) fails to compute all the
exact zero entries of $C$ because of inexact  
cancellation.}
\jwd{Test\_1c refers to $\sqrt{OV}$ and $\sqrt{UN}$. Could  
you say how these correspond to the ``scaling factors''
in your tests? Ditto for Test\_2c.}
}

\clearpage
\newgeometry{margin=0.5in}
\clearpage
\restoregeometry

\clearpage
\newgeometry{margin=0.5in}
\clearpage
\restoregeometry

\clearpage
\newgeometry{margin=0.5in}
\clearpage
\restoregeometry

\clearpage
\newgeometry{margin=0.5in}
\clearpage
\restoregeometry

\clearpage
\newgeometry{margin=0.5in}
\clearpage
\restoregeometry

\clearpage
\newgeometry{margin=0.5in}
\clearpage
\restoregeometry

\clearpage
\newgeometry{margin=0.5in}
\clearpage
\restoregeometry

\clearpage
\newgeometry{margin=0.5in}
\clearpage
\restoregeometry

\clearpage
\newgeometry{margin=0.5in}
\clearpage
\restoregeometry

\clearpage
\newgeometry{margin=0.5in}
\clearpage
\restoregeometry

\section{Grading triangular solve: \trsm}
\label{sec:trsm-residual-gemm}

\trsm differs slightly from \gemm in what error metric is the most useful.
In the scalar case, $b = q x + r$ is an exact relationship in faithfully rounded arithmetic for real valued quotients $q$ and remainders $r$.
This, too, is a convenient form for the matrix version,
\begin{equation}
  \label{eq:trsm-res}
    B = T \Xhat + R \text{ with } \Xhat = \dotrsm{T}{B},
\end{equation}
where $T$ is a triangular matrix. Here we call $R$ the residual.

We propose \trsm grades analogous to the proposed \gemm grades but phrased in terms of the residual.
The ``A,'' ``B,'' and ``C'' grading scales are listed in Table~\ref{tab:notation-grades}.
The loosest bound, a ``C'' grade, is the normwise residual bound derived from Equation~\eqref{eq:trsm-res}.
Given an $n \times n$ triangular matrix $T$,  $g_C(n) \leq n$ for standard \trsm. The form relying on Strassen's or Winograd's \gemm{}s has $g_C(n) \leq O(n^{\log_2 12})$ where $\log_2 12 > 3.6$ \cite{10.1145/98267.98290}.
The stricter ``A'' grade follows the traditional componentwise bounds, and analysis generally shows that $g_A(n) \leq n$
\cite{Higham2002}.

We extend these traditional measures with an intermediate grade ``B'':
\begin{equation}
  \label{eq:trsm-mixed}
    | \,R \, | \leq g_B(n) \, \epsilon \rowmax(|T|) \cdot \colmax(|\Xhat|).
\end{equation}
This bound treats each column of $\Xhat$ separately, hence $T X = B$ behaves as a separate solution per column.
Future work will examine $g_B$ as well as the convergence of iterative refinement using a ``B'' \trsm and \gemm combination.


\subsection{Tests for \trsm implementations}

We assume a non-adversarial role and do not consider much active ``gaming.''
We start by re-using the numerical tests from Section~\ref{sec:testing}.
Solving the system $[I, -A; 0, I] X = [0; B]$ needs to compute $\Xhat = [\dogemm{A}{B}; B]$, so then \gemm tests have natural \trsm counterparts.
All the ``gaming'' tricks apply if implementations scan for this structure.

We also generalize Test\_1 for \trsm by identifying
structurally exact zero entries.
Repeat the following steps a few times to introduce structural zeros into the computed $\Xhat$:
\begin{enumerate}
\item Pick a column $j$ in $B$.
\item Select a random set of rows $\iI_j$ and set $\midx{B}{\iI_j}{j} = 0$.
\item For all $k \in \operatorname{complement}(\iI_j)$, set $\midx{L}{\iI_j}{k} = 0$. 
\end{enumerate}
This effectively breaks the \emph{column elimination tree}~\cite{10.1145/356004.356006} into a column elimination forest.
The tree starts as a linear graph (a trivial tree).
The procedure above separates subtrees from the root vertex (the first column).
%
Strassen-like methods in general will introduce nonzeros where we
placed (structural) zeros so long as some of those introduced zeros straddle
the base-case blocks.

Constructing more complicated tests along the lines of Section~\ref{sec:testing} is feasible
but unlikely to show interesting behavior.
Reductions from \trsm to \gemm and \gemm to \trsm establish forms of equivalence groupings that imply that the grading tests also are equivalent, although practical implementations would not rely on the reductions.




\subsection{\gemm and \trsm grades are equivalent}
\label{sec:gemm-trsm-grades}

Any \gemm that achieves a certain grade can produce a \trsm of the same grade,
and the same for \trsm producing \gemm. The gradations are separate categories in a sense.
A \trsm that does not achieve a grade of ``B'' cannot
produce a \gemm of ``B'' or better, and similarly for \gemm to \trsm.
Given the derivation of BLAS from \gemm and \trsm\cite{K_gstr_m_1998,10.1145/98267.98290},
we expect the three grades to be separate across all routines.

\subsubsection{Using \gemm to generate a same-grade \trsm}
\label{sec:gemm-to-trsm}

Without loss of generality, we will focus on \emph{upper-triangular}, $n \times n$ systems,
\begin{equation}
  \begin{bmatrix}
        \matblock{U}{1}{1} & \matblock{U}{1}{2} \\  & \matblock{U}{2}{2}
  \end{bmatrix}
  \begin{bmatrix}
        \vecaccess{X}{1} \\  \vecaccess{X}{2}
  \end{bmatrix} =
  \begin{bmatrix}
        \vecaccess{B}{1} \\ \vecaccess{B}{2}
  \end{bmatrix}.
\end{equation}
Let $\matblock{U}{1}{1}$ be an upper-triangular, $n_1 \times n_1$ matrix, and let $\matblock{U}{2}{2}$ be an upper-triangular, $n_2 \times n_2$ with $n_1 + n_2 = n$.
The base case considers only a scalar system $u x = b$ solved by scalar division.
This incurs a single rounding error outside of division-by-zero and other exceptional cases.

Following implementations from \cite{10.1145/98267.98290,K_gstr_m_1998}, we compute $\Xhat$ in a partitioned manner,
\begin{align}
    \vecaccess{\Xhat}{2} &= \ldotrsm{\matblock{U}{2}{2}}{\vecaccess{B}{2}}, \text{and} \\
    \vecaccess{\Xhat}{1} &= \ldotrsm{\matblock{U}{1}{1}}{\left(\vecaccess{B}{1} \ominus \ldogemm{\matblock{U}{1}{2}}{\vecaccess{\Xhat}{2}}\right)} .
\end{align}

If $\matblock{U}{2}{2}$ is scalar and we compute all of $\vecaccess{\Xhat}{1}$ as a vector operation, this corresponds to the \axpy version of forward substitution\cite{golub-van-loan}. If instead we expand $\vecaccess{\Xhat}{1}$ row by row then this corresponds to the \blasdot version.
These differ only in the order of summation.
We will lump these differences into the polynomial factors $g(n)$ in our evaluation, but note that operation order can affect cancellation dramatically.
If we extend in either fashion and then switch to block operations, then this can be an iterative routine as opposed to a recursive definition\footnote{See \url{https://github.com/Reference-LAPACK/lapack/blob/master/SRC/VARIANTS/lu/REC/dgetrf.f} for a similar iterative expansion of Toledo's recursive LU\cite{doi:10.1137/S0895479896297744}}.

The grade of this kind of \trsm implementation then depends on the grade of the used \gemm. We sketch enough of the analysis from the appendix of~\cite{10.1145/98267.98290} to demonstrate \emph{how} the \gemm and \trsm grades are linked. The computed solution $\Xhat$ satisfies
\begin{equation}
  \begin{bmatrix}
        \matblock{U}{1}{1} & \matblock{U}{1}{2} \\ & \matblock{U}{2}{2}
  \end{bmatrix}
  \begin{bmatrix}
        \vecaccess{\Xhat}{1} \\  \vecaccess{\Xhat}{2}
  \end{bmatrix}
  +
  \begin{bmatrix}
        \vecaccess{R}{1} \\ \vecaccess{R}{2}
  \end{bmatrix}
  =
  \begin{bmatrix}
        \vecaccess{B}{1} \\ \vecaccess{B}{2}
  \end{bmatrix}.
\end{equation}
We assume that the subproblem $\vecaccess{R}{2} = \vecaccess{B}{2} - \matblock{U}{2}{2} \vecaccess{\Xhat}{2}$ is of the form of Equation~\eqref{eq:trsm-res} computed either recursively or iteratively as explained above.

\begin{equation}\label{eq:trsm-errs}
  \begin{split}
    \matblock{U}{1}{1} \vecaccess{\Xhat}{1} &= \left(\vecaccess{B}{1} \ominus \left(\matblock{U}{1}{2} \vecaccess{\Xhat}{2} + \dbox{M}\right) + S\right) + \vecaccess{R}{1}' \\
    &= (\vecaccess{B}{1} - \matblock{U}{1}{2} \vecaccess{\Xhat}{2}) - \dbox{M} + S + \vecaccess{R}{1}' = Z + \vecaccess{R}{1} .
  \end{split}
\end{equation}
Here $\vecaccess{R}{1}'$ is the residual from the subproblem of solving
$\matblock{U}{1}{1} \vecaccess{\Xhat}{1} = Z$, $S$ is the error in the subtraction $\ominus$, and $M$ is
the error in $\dogemm{\matblock{U}{1}{2}}{\vecaccess{\Xhat}{2}}$. $M$ is highlighted because this is the
role of \gemm in the routine. The residual $R'_1$ is bound by the (recursive)
assumption on the final \trsm routine.

The grades follow immediately from the bounds on each error contribution.
See~\cite{10.1145/98267.98290} for the details. The key pieces are to bound the
$\vecaccess{B}{1}$ contribution in the subtraction using Equation~\eqref{eq:trsm-errs}
itself and extend piece-wise bounds on submatrices to the entire matrix, \textit{e.g.,} $\colmax(|\matblock{U}{\dimall}{2}|) \geq \colmax(|\matblock{U}{2}{2}|)$.

\subsubsection{Using \trsm to generate a same-grade \gemm}
\label{sec:trsm-grade-to-gemm}

A \trsm with a given grade also can generate a \gemm with the same grade. This
works through a construction rather than an implicit induction. We consider a
$3n \times 3n$ upper triangular matrix
\begin{equation}
  T =
  \begin{bmatrix}
    I & -D_A^{-1} A & \\
      & I & -B D_B^{-1} \\
      &   & I
  \end{bmatrix}
\end{equation}
and solve $T X = [0; 0; I]$.
We further assume that the entries of $D_A$ and $D_B$ are powers of the radix such that the scaling operation is exact.
In exact arithmetic,
\begin{equation}
  \label{eq:trsm-exact-soln}
  X =
  \begin{bmatrix}
    (D_A^{-1} A) (B D_B^{-1}) \\ B D_B^{-1}\\ I
  \end{bmatrix}
  .
\end{equation}
The computed solution,
\begin{equation}
  \Xhat =
  \begin{bmatrix}
        \vecaccess{\Xhat}{1}\\ \vecaccess{\Xhat}{2}\\ \vecaccess{\Xhat}{3}
  \end{bmatrix} =
  \begin{bmatrix}
    \opeval{(D_A^{-1} A)(B D_B^{-1})} \\ B D_B^{-1} \\ I
  \end{bmatrix},
\end{equation}
has exact entries $\vecaccess{\Xhat}{2}$ and $\vecaccess{\Xhat}{3}$.
We define
$\dogemm{(D_A^{-1} A)}{(B D_B^{-1})} = \vecaccess{\Xhat}{1}$ as our resulting \gemm
routine.

Partition the residual $R = [0; 0; I] - T \Xhat$ row-wise in the same manner. Then $\vecaccess{R}{2} = \vecaccess{R}{3} = 0$ because those computations are exact.
However, $\vecaccess{R}{1} = (D_A^{-1} A) (B D_B^{-1}) - \dogemm{(D_A^{-1} A)}{(B D_B^{-1})} = E_S$, the error in the \gemm. We denote it $E_S$ here to emphasize that it is scaled by $D_A$ and $D_B$.
We will use different scaling matrices along with the \trsm bounds in Table~\ref{tab:notation-grades} to produce equivalently-graded \gemm computations.

\paragraph*{Grade ``C''}

The most direct case is when \trsm is graded ``C.'' Let $D_A = \diag(\|A\|),$ and $D_B = \diag(\|B\|).$
Scaling by $D_A^{-1}$ and $D_B^{-1}$ is the same as multiplying the matrix by the scalars $1/\|A\|$ and $1/\|B\|$ respectively.
Then $\|D_A^{-1} A\| = 1,$ and $\|B D_B^{-1}\| = 1.$
So
$\|E_S\| = 1 / (\|A\| \, \|B\|) \|E\| \leq g_C(n)$ and
$\|E\| \leq f_C(n) \|A\|\,\|B\|$ as desired.

\paragraph*{Grade ``B''}
Let $D_A = \diag(\rowmax(|A|)),$ and $D_B = \diag(\colmax(|B|))$.
Then $|D_A^{-1} A| = \mathbf{1}_n,$ and $|B D_B^{-1}| = \mathbf{1}_n^T,$ where $\mathbf{1}_n$ is an $n$-vector with unit entries.
So $|E_S| = |D_A^{-1} E D_B^{-1}| = D_A^{-1} |E| D_B^{-1} \leq g_B(n) \mathbf{1}_n \mathbf{1}_n^T,$
and scaling from the left with $D_A$ and the right with $D_B$ produces
$|E| \leq g_B(n) \allowbreak \rowmax(|A|) \, \colmax(|B|)$.

\paragraph*{Grade ``A''}
This grade does not rely on scaling., $D_A = I$, and $D_B = I$. Then $\vecaccess{\Xhat}{1} = \dogemm{A}{B}$, and $\vecaccess{R}{1} = \dogemm{A}{B} - AB = E$.
In the bound for the ``A'' \trsm grade,
\begin{equation}
  |T| \cdot |\Xhat| =
  \begin{bmatrix}
    |\vecaccess{\Xhat}{1}| + |A| \cdot |B| \\ 2 |B| \\ I
  \end{bmatrix} .
\end{equation}
Using a running error bound\cite{Higham2002}, $|\vecaccess{\Xhat}{1}| \leq (2n-1) |A| \cdot |B|$,
\begin{align}
    |E| &\leq g_A(3n) \epsilon \left(|\vecaccess{\Xhat}{1}| + |A| \cdot |B|\right) \\
  &\leq (2n-1) g_A(3n) \epsilon |A| \cdot |B| = f_A(n) \epsilon |A| \cdot |B| .
\end{align}

\section{Other BLAS 3-like operations}\label{sec:otherBLAS}

The other BLAS 3 operations (\trmm, \gemmtr, \syrk, \syrtwok) and
structured forms (symmetric, banded, \textit{etc.}) reduce to \gemm and
\trsm\cite{K_gstr_m_1998,10.1145/98267.98290}.
Each of the symmetric or triangular matrices are \emph{structurally} symmetric or triangular.
No errors place non-zeros in the wrong triangle or make $S$ non-symmetric.
All but \trsm reduce to the same grades and tests as \gemm in Section~\ref{sec:grading}.
Later work will expand the test code for \uplo (upper/lower-triangular), \trans
(transpose / Hermitian), and other BLAS options. This should leverage the same
work as needed to test consistent exception handling
(Section~\ref{sec:Exceptions}).

Other BLAS 3-like operations, \textit{e.g.,} LAPACK helper routines, also need
future work. The $WY$ form for aggregated Householder reflections (\larfb,
\larzb) apply products of Householder orthogonal transformations in \gemm
operations. The backward error bounds differ between routines based on grade
``A'' and ``C'' \gemm implementations differ only in constant
factors\cite{Higham2002}. Other blocked or aggregated orthogonal transformations
may behave similarly on their own, but the usage in higher-level routines like
eigensolvers deserve future attention. We are not aware of existing work on the
impact of emulated precisions with limited range upon eigensolvers and other
orthogonal decompositions.


\section{Impact on Higher Level Algorithms}\label{sec:LUQRChol}


The importance of the grade that a BLAS implementation receives ultimately depends on its
impact on higher level algorithms using the BLAS.
In this section we discuss the impact on three
important linear algebra algorithms, Cholesky
(for solving a symmetric positive definite (spd)
linear system $Ax=b$), LU (for solving a general
linear system) and QR (for solving a least squares
problem ${\rm argmin}_x \|Ax-b\|_2$). The classical
(backward) error analyses of these problems
is done normwise, \textit{i.e.,} it shows that the computed
result is the exact solution for a problem
with $A + E$ in place of $A$, where $\|E\| = O(\epsilon) \|A\|$, with error bounds proportional
to $\epsilon$ times a condition number, \textit{e.g.,}
$\kappa(A) = \|A\| \cdot \|A^{-1}\|$ for solving
a linear system.

To see why we should expect to do better,
consider solving $Ax=b$ using LU, implemented
straightforwardly using conventional floating point without any BLAS calls.
It is easy to see that the algorithm is
``column scaling invariant,'' \textit{i.e.,} if $D$ is
diagonal and nonsingular, changing $A$ to $AD$
changes $A=PLU$ to $AD = PL(UD)$, \textit{i.e.,} the
column scale factors are absorbed into the $U$ factor.
And solving $(AD)\xhat = b$ instead of $Ax=b$
yields $\xhat = D^{-1}x$. In fact, if the
diagonal entries of $D$ are powers of 2, then
(ignoring over/underflow, and assuming a
deterministic order of execution), exactly the same
rounding errors are made throughout the algorithm,
\textit{i.e.,} all intermediate results have the same mantissas. In other words, the error (in $\xhat = D^{-1}x$) should
depend on the condition number of $AD$, where $D$ is
(implicitly) chosen to minimize the condition
number $\kappa(AD)$ of $AD$. It is possible
that $\kappa(AD) \ll \kappa(A)$ if $A$ has
columns of very different norms, and $D$ is
chosen to make all the column norms (nearly) equal.

As an example, suppose each entry of
$A$ is an independent uniform random
number in the range $[-1,1]$ (so all
following statements about accuracy are true with high probability, w.h.p.). Let $D$ be diagonal
with diagonal entries that are powers
of 2, some very large and some very small.
Let $b$ also be random with independent uniform
entries from $[-1,1]$. Then (w.h.p.)
$A$ is well-conditioned, the entries
of the true solution of $Ax=b$ are all
$\Theta(1)$, and the computed solution 
$\tilde{x}$ of $Ax=b$ has a small relative
error in each component.
The true solution $\xhat  = D^{-1}x$ of $(AD)\xhat=b$ will have some very large
and very small entries (corresponding to
very small and very large diagonal entries
of $D$). Finally, the computed solution
$D^{-1}\tilde{x}$ of $(AD)\xhat=b$
will have
some very large and very small entries,
{\em all of which are computed to high
relative accuracy}. 
This is a much stronger
error bound than the usual one depending
just on the very large condition number 
of $AD$, which can only guarantee a very
large error bound on all entries of the
computed solution, large or small.
We consider the first, tighter
error bound above as deserving an ``A'' grade, and the
second, looser error bound a ``C'' grade. 
Cholesky and QR are analogous.

Our question is what happens when using
blocked/tiled versions of these algorithms,
using BLAS to accelerate them, as in LAPACK
\cite{LAPACK:1999}
and related libraries. 
A natural conjecture is 
that if the underlying BLAS get an
``A'' or ``B'' grade, then LU, QR and
Cholesky get an ``A'' grade, but if the
BLAS get a ``C'', so do the higher level
routines; we prove this in detail for 
one version of blocked Cholesky below.

Given the large design space of ways to 
tile LU, QR, and Cholesky, we will restrict
ourselves to Cholesky, tiled using only
the most straightforward way, defined below
in Algorithm~\ref{alg:Cholesky}.
But first, we state known error bounds for
Cholesky, corresponding to ``A'' and ``C'' grades,
and applicable to the unblocked algorithm
using conventional floating point.
We first define our notation: $A$ is a real $n \times n$ positive definite matrix, and we want to solve
$Ax=b$. We can write $A = DHD$ where $D$ is
diagonal with $\matentry{D}{i}{i} = \matentry{A}{i}{i}^{1/2}$,
$\matentry{H}{i}{i} = 1$, and $H$ is also positive definite.
When $D$ is ill-conditioned, it is possible
that the condition number $\kappa(H) \ll \kappa(A)$, and in fact
$\kappa(H) \leq n \cdot \min_D \kappa (DAD)$
where is the minimum is over all nonsingular
diagonal $D$ \cite{vandersluis}.

\noindent
{\bf Theorem~1:} \cite[Sec 10.1]{Higham2002}\cite{Demmel_Cholesky_89,Wilkinson_Cholesky_68}. Suppose the unblocked Cholesky factorization
of $A$ with machine precision $\epsilon$
runs to completion, yielding a triangular
factor $\hat{L}$, and a computed solution $\xhat$
of $Ax=b$.
Let $E = \hat{L}\hat{L}^T - A$, and let the
diagonal matrix $D$ satisfy $\matentry{D}{i}{i} = \matentry{A}{i}{i}^{1/2}$.
Then
the results satisfy the ``C'' grade error bounds:
\begin{equation}\label{eqn:CholBndC}
\|E\|_2 = O(n \epsilon) \|A\|_2
\hspace{.3in} {\rm and} \hspace{.3in}
\frac{\|x - \xhat\|_2}{\|x\|_2} =O(n^2 \epsilon)  \kappa(A)
\end{equation}
as well as the stronger ``A'' grade error bounds:
\begin{equation}\label{eqn:CholBndA}
|\matentry{E}{i}{j}| = O(n \epsilon) (\matentry{A}{i}{i} \matentry{A}{j}{j})^{1/2}
\hspace{.3in} {\rm and} \hspace{.3in}
\frac{\|D(x - \xhat)\|_2}{\|Dx\|_2} =O(n^2 \epsilon)  \kappa(H)
\end{equation}

We now give
a ``random example'' in the same style as for
LU decomposition above, to show when 
bounds (\ref{eqn:CholBndA}) provide much
tighter error bounds than (\ref{eqn:CholBndC}). 
Let each entry of an $n \times n$
matrix $B$ be an independent uniform random number
in the range $[-1,1]$, so w.h.p. each diagonal
entry of $B^T B$ will be $\Theta(1)$, and
$B^TB$ will be well-conditioned. As above,
let $\hat{D}$ be diagonal with some very large
and very small diagonal entries, and let 
$A=\hat{D}B^TB\hat{D}$. 
$A$ will be very ill-conditioned w.h.p., with 
$\kappa(A) = \Theta( ( \max_i \matentry{D}[\hat]{i}{i} / \min_i \matentry{D}[\hat]{i}{i} )^2)$. Let $b$ also be random with
independent uniform entries in $[-1,1]$.
Then the true solution of $Ax=b$ will have some
very large and very small entries (corresponding to
very small and very large diagonal entries of $D$).
Based on (\ref{eqn:CholBndA}) all the computed entries of $x$ will be accurate
to high relative accuracy, independent of their
magnitudes. This is in contrast to the bound (\ref{eqn:CholBndC}), which is proportional to 
$\kappa(A)$,
which is much larger for all entries.

Our goal is to show that tiled Cholesky attains
the same error bounds as in (\ref{eqn:CholBndA}), as long
as the implementations of \gemm (including \syrk) and \trsm used
get either an ``A'' or a ``B'', and otherwise
the bounds in (\ref{eqn:CholBndC}).

\begin{algorithm}[H]
\caption{Tiled Cholesky to solve $Ax=b$}
\label{alg:Cholesky}
\begin{algorithmic}[1]
\Require Symmetric Positive Definite $n \times n$ 
matrix $A$ (lower triangle only), block size $b$ where $b | n$
\Ensure $x$ satisfies $Ax=b$, and $A = L \cdot L^T$,
where $A$ is on input, and $L$ is stored in the lower
triangle of $A$ on output (both modulo roundoff)

\algrule
\Statex Partition $A$ into $N \times N$ blocks of size $b \times b$
    with $\matblock{A}{i}{j}$ the $(i,j)$ submatrix.
        \For{$i = 1, \dots, N$}
        \State Factor $\matblock{A}{i}{i} = L \cdot L^T$ using conventional floating point without BLAS; overwrite $\matblock{A}{i}{i}$ with $L$
        \If {$i<N$}
        \State $\matblock{A}{\dimrange{i+1}{N}}{i} =\dotrsmT{\matblock{A}{i}{i}}{\matblock{A}{\dimrange{i+1}{N}}{i}}$
        \State $\matblock{A}{\dimrange{i+1}{N}}{\dimrange{i+1}{N}} = \dosyrk{\matblock{A}{\dimrange{i+1}{N}}{\dimrange{i+1}{N}}}{\matblock{A}{\dimrange{i+1}{N}}{i}}$ \Comment{Update lower triangle only}
        \EndIf
        \EndFor
        \State Solve $L \cdot y = b$ using conventional floating point, overwrite $b$ with $y$
        \State Solve $L^T \cdot x = y$ using conventional floating point, overwrite $b$ with $x$
\end{algorithmic}
\end{algorithm}

\noindent
{\bf Theorem~2:}
If both \trsm and \syrk get either ``A'' or ``B'' grades,
then the output of Algorithm~1 satisfies bound~(\ref{eqn:CholBndA}).
If \trsm and \syrk only get (at least) a ``C'' grade,
then the output of Algorithm~1 only satisfies
bound~(\ref{eqn:CholBndC}).

\begin{proof}
To simplify the error bounds,
we use $f(n)$ to denote any positive function bounded by a polynomial in $n$, since this is all we need
in an error bound to determine a grade.
We use induction based on
the following block factorization:
\[
\aligned
A & = \begin{bmatrix} \matblock{A}{1}{1} & \matblock{A}{2}{1}^T \\ \matblock{A}{2}{1} & \matblock{A}{2}{2}
\end{bmatrix} = 
\begin{bmatrix}
    \matblock{L}{1}{1}\matblock{L}{1}{1}^T & \matblock{A}{2}{1}^T \\ \matblock{A}{2}{1} & \matblock{A}{2}{2}
\end{bmatrix} =
\begin{bmatrix}
    \matblock{L}{1}{1} & 0 \\ \matblock{L}{2}{1} & I
\end{bmatrix} \cdot
\begin{bmatrix}
    \matblock{L}{1}{1}^T & \matblock{L}{2}{1}^T \\ 0 & \matblock{A}{2}{2} - \matblock{L}{2}{1} \cdot \matblock{L}{2}{1}^T
\end{bmatrix} \\
& =
\begin{bmatrix}
    \matblock{L}{1}{1} & 0 \\ \matblock{L}{2}{1} & I
\end{bmatrix} \cdot
\begin{bmatrix}
    \matblock{L}{1}{1}^T & \matblock{L}{2}{1}^T \\ 0 &  \matblock{L}{2}{2} \cdot \matblock{L}{2}{2}^T
\end{bmatrix} =
\begin{bmatrix}
    \matblock{L}{1}{1} & 0 \\ \matblock{L}{2}{1} & \matblock{L}{2}{2}
\end{bmatrix} \cdot
\begin{bmatrix}
    \matblock{L}{1}{1}^T & \matblock{L}{2}{1}^T \\ 0 &  \matblock{L}{2}{2}^T
\end{bmatrix}
\endaligned
\]
where 
$\matblock{A}{1}{1}$ is $b \times b$ and
$\matblock{L}{2}{1} = \matblock{A}{2}{1} \cdot \matblock{L}{1}{1}^{-T}$.
$\matblock{A}{1}{1} = \matblock{L}{1}{1} \cdot \matblock{L}{1}{1}^T$ is
factored in line 2 of Algorithm~\ref{alg:Cholesky},
$\matblock{L}{2}{1}$ is computed in line~4 and
$\matblock{A}{2}{2} - \matblock{L}{2}{1} \cdot \matblock{L}{2}{1}^T$ in line~5.

Denote computed quantities by $\matblock{L}[\hat]{i}{j}$
in lines 2 and 4, and
$\matblock{A}[\hat]{2}{2}$ in line 5.
Denote quantities that we need to bound by
\[
\aligned
\matblock{E}{1}{1} & = \matblock{A}{1}{1} - \matblock{L}[\hat]{1}{1} \cdot \matblock{L}[\hat]{1}{1}^T \\
\matblock{E}{2}{1} & = \matblock{A}{2}{1} - \matblock{L}[\hat]{2}{1} \cdot \matblock{L}[\hat]{1}{1}^T \\
\matblock{E}{2}{2} & = \matblock{A}[\hat]{2}{2} - (\matblock{A}{2}{2} - \matblock{L}[\hat]{2}{1} \cdot \matblock{L}[\hat]{2}{1}^T) \\
\matblock{E}{2}{2}' & = \matblock{A}[\hat]{2}{2} - \matblock{L}[\hat]{2}{2} \cdot \matblock{L}[\hat]{2}{2}^T
\endaligned
\]
Since $\matblock{L}[\hat]{1}{1}$ is computed using
conventional floating point without BLAS,
$\matblock{E}{1}{1}$ is bounded by
\begin{equation}\label{eqn_E11bnd}
|\matentry{\matblock{E}{1}{1}}{i}{j}| \leq f(n) \epsilon (\matentry{\matblock{A}{1}{1}}{i}{i}\matentry{\matblock{A}{1}{1}}{j}{j})^{1/2}
\end{equation}
based on (\ref{eqn:CholBndA}).
This is also the base case of the induction.

First, we assume that $\syrk$ and $\trsm$ both
get an ``A'' or ``B'' grade. We note that for
Cholesky to complete, 
$\|\matentry{L}[\hat]{i}{1:j}\|_2^2 \leq \matentry{A}{i}{i}$ up to roundoff.
Based on (\ref{eq:trsm-mixed})
\[
|\matentry{\matblock{E}{2}{1}}{i}{j}| \leq f(n) \epsilon \| \matentry{\matblock{L}[\hat]{2}{1}}{i}{\dimall} \| \cdot \| \matentry{\matblock{L}[\hat]{1}{1}^T}{\dimall}{j} \|
= f(n) \epsilon \| \matentry{\matblock{L}[\hat]{2}{1}}{i}{\dimall} \|_2 \cdot \| \matentry{\matblock{L}[\hat]{1}{1}}{j}{\dimall} \|_2 \leq f(n) \epsilon (\matentry{A}{i}{i} \matentry{A}{j}{j})^{1/2}
\]
Next,
\[
\matblock{E}{2}{2}'' \equiv \matblock{A}{2}{2} - \matblock{L}[\hat]{2}{1}\cdot \matblock{L}[\hat]{2}{1}^T -
\matblock{L}[\hat]{2}{2} \cdot \matblock{L}[\hat]{2}{2}^T = \matblock{E}{2}{2}' - \matblock{E}{2}{2}
\]
so we can bound
\begin{equation*}
\aligned
|\matentry{\matblock{E}{2}{2}}{i}{j}| & \leq f(n) \epsilon (|\matentry{\matblock{A}{2}{2}}{i}{j}| +
\|\matentry{\matblock{L}[\hat]{2}{1}}{i}{\dimall}\|_2 \cdot \|\matentry{\matblock{L}[\hat]{2}{1}^T}{\dimall}{j}\|_2)
\\
& = f(n) \epsilon (|\matentry{\matblock{A}{2}{2}}{i}{j}| +
\|\matentry{\matblock{L}[\hat]{2}{1}}{i}{\dimall}\|_2 \cdot \|\matentry{\matblock{L}[\hat]{2}{1}}{j}{\dimall}\|_2) \\
& \leq f(n) \epsilon (\matentry{\matblock{A}{2}{2}}{i}{i} \matentry{\matblock{A}{2}{2}}{j}{j})^{1/2} \\
& \leq f(n) \epsilon (\matentry{A}{i+b}{i+b} \matentry{A}{j+b}{j+b})^{1/2}
\endaligned
\end{equation*}
and, by induction,
\[
|\matentry{\matblock{E}{2}{2}'}{i}{j}| \leq f(n) \epsilon (\matentry{\matblock{A}[\hat]{2}{2}}{i}{i}
\matentry{\matblock{A}[\hat]{2}{2}}{j}{j})^{1/2})
\leq f(n) \epsilon (\matentry{A}{i+b}{i+b} \matentry{A}{j+b}{j+b})^{1/2})
\]
so adding the bounds on $|\matblock{E}{2}{2}|$ and $|\matblock{E}{2}{2}'|$ yields
\[
|\matentry{\matblock{E}{2}{2}''}{i}{j}| \leq f(n) \epsilon (\matentry{A}{i+b}{i+b} \matentry{A}{j+b}{j+b})^{1/2}
\]
Altogether, we have the desired error bound in (\ref{eqn:CholBndA}) 
for $E = A - \hat{L}\cdot \hat{L}^T$:
\[
|\matentry{E}{i}{j}| \leq
f(n) \epsilon (\matentry{A}{i}{i}\matentry{A}{j}{j})^{1/2} .
\]
The bound for $x - \xhat$ in (\ref{eqn:CholBndA})
follows from this 
\cite[Sec 10.1]{Higham2002}\cite{Demmel_Cholesky_89,Wilkinson_Cholesky_68}.

Second, we assume that $\syrk$ and $\trsm$ only get a
``C'' grade. First, we note that
\[
\| \matblock{E}{1}{1} \|_2 \leq f(n) \epsilon \|\matblock{A}{1}{1}\|_2 \leq f(n) \epsilon \|A\|_2
\]
follows directly from (\ref{eqn_E11bnd}).
Similarly,
\[
\| \matblock{E}{2}{1} \|_2 \leq
f(n) \epsilon \| \matblock{L}[\hat]{2}{1} \|_2 \cdot
\| \matblock{L}[\hat]{1}{1} \|_2
\leq f(n) \epsilon 
\| A \|_2^{1/2} \cdot 
\| A\|_2^{1/2} =
f(n) \epsilon \|A\|_2
\]
follows from the last row of Table~1
and $\|\matblock{L}{i}{j}\|_2 \leq \|L\|_2 = \|A\|_2^{1/2}$
(the difference between $\|L\|_2$ and $\|\hat{L}\|_2$
is small and can be absorbed in the $f(n)$ factor).
The bound $\|\matblock{E}{2}{2}\|_2 \leq f(n) \epsilon \|L\|_2^2
= f(n) \epsilon \|A\|_2$ similarly follows from the
definition of a ``C'' grade for $\gemm$ and
analogously $\syrk$. The bound
$\|\matblock{E}{2}{2}'\|_2 \leq f(n) \epsilon \| \matblock{A}[\hat]{2}{2} \|_2
\leq f(n) \epsilon \|A\|_2$ holds by induction.
Adding the bounds on $\|\matblock{E}{1}{1}\|_2$, $\|\matblock{E}{2}{1}\|_2$,
$\|\matblock{E}{2}{2}\|_2$ and $\|\matblock{E}{2}{2}'\|_2$ yields the
desired  bound 
\[
\|A - \hat{L} \hat{L}^T \|_2 \leq f(n) \epsilon \|A\|_2
\]
needed for a ``C'' grade.
\end{proof}

Finally, we note that just as it is possible to
convert a grade ``C'' \gemm into a grade ``B''
by (1) scaling $A' = D_1 \cdot A$ and $B' = B \cdot D_2$ to have rows (resp. columns) of unit norm,
(2) multiplying $C' = A' \cdot B'$ and 
(3) unscaling to get $C = D_1^{-1} \cdot C' \cdot D_2^{-1}$, it is possible to do an analogous scaling
trick to improve the grade of Cholesky from ``C'' to
``A'': (1) scale $A' = D \cdot A \cdot D$ to have
unit diagonal entries $A'(i,i)=1$,
(2) solve $A' \cdot x' = D \cdot b$ using ``C'' 
grade Cholesky, and
(3) unscale $x = D \cdot x'$. 
This adds an extra $O(n^2)$ cost to the overall solve
cost.

\section{Consistent Exception Handling}
\label{sec:Exceptions}

We briefly review the proposed requirements for how the
BLAS should handle floating point exceptions
\cite{ExceptionHandling-1}
and how to assign grades by testing. At a high level, the
proposed rule is that exceptional values (both \Inf and
\NaN) that are either input to a routine, or created
internally, should propagate to the output unless
there is a well-defined mathematical reason they should not.

We do not prescribe the conditions under which an \Inf versus
a \NaN should propagate, because there are situations in
which either is a plausible outcome. For example, computing
the sum of 5 finite numbers $[x, x, 3, -x, -x]$, where
$x$ is greater than half the overflow threshold, could yield
$+\Inf$, $-\Inf$, \NaN, 0 or 3 (the right answer), depending
on the order of summation. Also, in the case
of complex arithmetic, programming languages
and compilers do not all agree on the definition
of complex multiplication, and how it should
propagate exceptions
\cite{ExceptionHandling-1}.

We can also evaluate to which locations an exceptional
value propagates. For example, an \Inf or \NaN in row $i$
of $A$ will only cause the entire row $i$ of $C = A \cdot B$ to 
contain \Infs or \NaNs when a conventional $O(n^3)$ algorithm
is used, but propagate to more entries if a Strassen-like
algorithm is used.

\gemm is an example where \Infs and \NaNs do not need to
propagate in some cases. Recall that \gemm computes
$C = \alpha \cdot A \cdot B + \beta \cdot C$ where
$\alpha$ and $\beta$ are scalars. If $\alpha$
(resp. $\beta$) is zero, then the input values of
$A$ and $B$ (resp. $C$) are not even accessed, and
it is not expected that any input \Infs or \NaNs propagate
to the output. There are analogous rules for other BLAS
\cite{ExceptionHandling-1}.

Our grading will be Pass/Fail (P/F), meaning
that an implementation either adheres to the
above rules, or there is at least one example
where it does not.

Regarding testing, it is straightforward to
test with input \Infs and \NaNs and check whether
one (or more) appears in the output when it
should, or does not appear when it shouldn't. 
And BLAS routines (certainly \gemm) are simple
enough to devise finite inputs that should create
overflows, and so \Infs, and possibly \NaNs, as in
the example above. Alternatively, one could
consider ``spoofing'' \cite{EXCVATE2025}, \textit{i.e.,}
artificially inserting \Infs or \NaNs into
intermediate results of the computation.
Implementing tests for consistent exception
handling is future work.

\section{Grading Fairly}
\label{sec:grading_fairly}
One contribution of this work is to establish whether a BLAS package is as accurate as needed/expected. The tests in this research help determine the accuracy and the previous section addresses functionality in terms of exceptions.

Because of the following points, methods that achieve the grade of A are acceptable
for nearly all purposes:
\begin{itemize}
\item IEEE-754 (2019) specifies the accuracy for basic computational units (FMA/Add/etc) on scalars, but not dot products or BLAS. One merely has to worry about the final accuracy of a BLAS library, even one using emulation, not whether it's IEEE-754 compliant or not
(modulo exception handling, discussed next).
\item Exception Handling can be inconsistent even in the NETLIB BLAS
\cite{ExceptionHandling-1}.
As pointed out by examples in the previous section, 
outputs can vary among 0, $\pm\Inf$, \NaN, or the correct answer (independent of how
underflow is handled).
\cite{ExceptionHandling-1} points out that
there is inconsistency across 
IEEE 754 implementations, compilers and
programming languages about issues like
how underflow is handled (gradual vs abrupt underflow), and how overflow is handled in
complex arithmetic. These are longstanding
inconsistencies that users (and testers) of
the BLAS, and higher level libraries,
must acknowledge and live with.
\item Even if the correctly rounded
answer could in principle be computed
(at very high cost),
there is no one right numerical answer for a BLAS like DDOT or DGEMM. Because dot products can be computed in any order, 
with round-off, floating-point answers can differ. This is precisely why we need the contributions of this paper: we need to set a standard for the BLAS based on accuracy and theory, since there's no reason to expect that a DGEMM from one library should get the same answer as a DGEMM from another library, even when both computations were done correctly. 
The subject of exact bit-wise reproducible results (which may be reproducible and inexact) is another topic and has its own possibly high costs \cite{ahrens2020reproblas}.
\item Emulation may do more work in exact arithmetic because of fixed-point calculations and because many of the algorithms must scale the inputs to fit in
a fixed point format.
This means that emulation schemes can
sometimes be more accurate than using
floating point, for example when the
exponent range of the input is narrow 
enough, even though in the worst case
it only gets a ``B'' instead of an
``A'' grade.
\end{itemize}

We conclude by saying that fair testing means that alternate approaches to the BLAS (such as using emulation or Strassen) should be considered reliable as long as they achieve an accuracy standard that we describe in this paper. 
Concerns about how exceptions are handled
(including underflow, overflow, and invalid)
can be important for some applications,
and testing this is part of our future work.

\section{Related Work}
\label{sec:related-work}

Teasing out intricacies of actual floating-point implementations has been a concern for quite some time~\cite{Paranoia,paranoia-code}.


NVIDIA uses similar hybrid techniques with emulation that Strassen users do: for Strassen, there is a size $n$ where anything smaller is computed with a regular DGEMM. Similarly, if one gets an estimate for the floating point exponent range required in the DGEMM computation, this can be used to provide enough splits for Ozaki-I, or moduli for Ozaki-II to cover that range. This can lead to an effective hybrid technique that's not dependent on the problem size so much as it's dependent on the data input range \cite{10.1145/3773656.3773670}. For performance on small problems, one might choose FP64 DGEMM over an emulation DGEMM. This combination of problem size and data input range can be used to provide a hybrid algorithm that chooses the fastest algorithm (in many cases emulation) if it's safe, and calls FP64 DGEMM otherwise.

Starting with NVIDIA cuBLAS 13.0u2, FP64 DGEMM emulation techniques like Ozaki-1 have been incorporated. Emulation with FP32 via bfloat16 was done even earlier in cuBLAS 12.9. But the hybrid FP64 technique we’ve been discussing has been employed since 13.0u2. If it doesn’t seem safe to run the faster FP64 emulation kernels, it can default to the standard FP64 floating-point kernels. There’s no guarantee that the final FP64 accuracy will be better or worse, but one can still rely on the strictest error bounds. This way, a user doesn’t have to worry about accuracy of a performant emulation technique.



The accuracy of the cuBLAS
has already been shown in Figure~\ref{fig:sweeps_cublas}.
This was run on a single 
$\text{NVIDIA}^{\textregistered}$ 
Blackwell B200 GPU and Ozaki-1 emulation was enabled 
(the environment variables 
$\text{CUBLAS\_EMULATE\_DOUBLE\_PRECISION}\allowbreak =\allowbreak 1$ and
$\text{CUBLAS\_EMULATION\_STRATEGY}\allowbreak =\allowbreak \text{EAGER}$
were set). Even when the strategy was set to EAGER, it used the FP64 DGEMM path for the Test 2 excessive range. But when incorrectly told a bogus range with the Test 2 input, errors match what was shown for Ozaki-1 (Figure~\ref{fig:sweeps_ozaki1}.) However, when run according to the documentation, instead of getting the huge errors expected in Ozaki-1, the cuBLAS automatically detects the wide range of Test 2 and defaults to the FP64 code path and still gives ``A'' accuracy.

One might ask what exponent range criterion needs to be satisfied to make the ``A'' grade error bounds apply. How does one correctly compute or estimate the required data input range? Tests on the reliability of these techniques are beyond the scope of this paper.

\section{Conclusions and Open Problems}\label{sec:Conclusions}

Motivated by the growing use of low-precision accelerators in applications that historically have depended on conventional floating point arithmetic, with a large and growing design space of different algorithms for matrix multiplication and related linear algebra problems, our goal is to simplify the error analysis of these algorithms and applications that use them. To this end, we have identified a simple way to ``grade'' their accuracy into a few categories (``A'', ``B'', and ``C''), provided publicly available tests to assign these grades, that cannot be ``gamed'' by closed-source implementations, and verified their behavior on a variety of matrix multiply implementations. We have also shown how the same grading scheme can be extended to BLAS routines beyond matrix multiplication, and some higher level algorithms like Cholesky.

A number of open problems remain, which we describe below.

\begin{enumerate}
    \item We should extend our test code to non-square matrices, since these may result in different algorithms being used, and so different grades. For example, there are Strassen-like divide-and-conquer algorithms that apply to rectangular matrices, instead of square matrices whose dimensions are powers of 2. 
    \item We already mentioned the possibility of an ``A+'' grade for accumulating dot products in higher precision. One can imagine changing our grading scheme to be
    parameterized by the internal precision used before
    rounding the result to the lower output precision.
    \item We should extend our grades and tests to complex matrices.
    Complex matrix multiplication has a larger design
    space than real matrix multiplication, both because of
    more algorithms, like ``Gauss's trick'' for reducing
    the cost from 4 to 3 real matrix multiplications, and more recently to
    2 real matrix multiplications \cite{Caday2026}, and
    different programming languages having different 
    requirements for multiplying complex numbers (see Section 2.2, ``Programming languages and compilers,'' of the extended report~\cite{ExceptionHandling-1}).
    \item We should extend our analysis of Cholesky to 
      include more higher-level algorithms that use the BLAS.
    \item We should devise test matrices that can provide
    better estimates of the worst-case $f(n)$ factors in the error bounds. 
    \item We should provide test matrices for \trsm, other BLAS, Cholesky, and other routines, and continue to update our publicly available test suite to make these and other improvements available.
    \item We should design and implement tests for consistent exception handling.
    \item We should design tests that determine how underflow is handled.
    This includes the error bounds
    satisfied when subnormal numbers are
    provided as inputs, or underflows
    (would) occur during execution. 
    This can also have an impact on
    higher level algorithms 
    \cite{Demmel_Underflow_1984}.
    \item Making Ozaki I and Ozaki II schemes hit the accuracy requirements for grade A involves the implementation determining the exponent range on the input and output. For instance, if Ozaki I has enough slices, or Ozaki II has the right number of moduli, there's no reason it wouldn't get an A grade. Consider the Test 1 sweep in Figure~\ref{fig:sweeps_ozaki1} and Figure~\ref{fig:sweeps_ozaki2}. For the $2^{\pm 0}$ unscaled case, the exponent range was small enough and the tests passed. Future work might be creating tests to see if these exponent ranges are computed correctly.
    \item It would be interesting to extend this
    work to BLAS and other algorithms using interval arithmetic, 
    but we leave this work to others.
\end{enumerate}

\newpage
\printbibliography

@misc{perlmutter,
  author       = {{National Energy Research Scientific Computing Center (NERSC)}},
  title        = {Perlmutter Supercomputer},
  howpublished = {HPE Cray EX system at Lawrence Berkeley National Laboratory},
  year         = {2026},
  url          = {https://www.nersc.gov/systems/perlmutter/},
  note         = {Retrieved 15 May 2026},
}

@inproceedings{NvidiaPASC25,
    author = {Bayraktar, H.},
    title = {Precision redefined: unlocking and delivering the full power of modern {GPU}s for scientific computing},
    booktitle = {{PASC25} Minisymposium on Matrix Multiplication}, 
    url = {https://scicomp.leeds.ac.uk/pasc25-minisymposium-on-matrix-multiplication/},
    year = {2025}}

@inproceedings{10.1145/3773656.3773670,
author = {Schwarz, Angelika and Anders, Anton and Brower, Cole and Bayraktar, Harun and Gunnels, John and Clark, Kate and Xu, RuQing G. and Rodriguez, Samuel and Cayrols, Sebastien and Tabaszewski, Pawel and Podlozhnyuk, Victor},
title = {Guaranteed {DGEMM} accuracy while using reduced precision tensor cores through extensions of the {Ozaki} scheme},
year = {2026},
isbn = {9798400720673},
publisher = {Association for Computing Machinery},
address = {New York, NY, USA},
url = {https://doi.org/10.1145/3773656.3773670},
doi = {10.1145/3773656.3773670},
booktitle = {Proceedings of the Supercomputing Asia and International Conference on High Performance Computing in Asia Pacific Region},
pages = {91–101},
numpages = {11},
location = {
},
series = {SCA/HPCAsia '26}
}

@inproceedings{BensonBallardSPPPP15,
author = {A. Benson and G. Ballard},
title = {A framework for practical  
         parallel fast matrix multiplication},
booktitle = {Proc. 20th ACM SIGPLAN
            Symp. on Principles and Practice of Parallel Programming},
year = {2015},
doi = {10.1145/2688500.2688513}}

@unpublished{PASC25,
  author = {Mary, T. and Mikaitis, M.},
  title = {{PASC25} minisymposium on fast and
           accurate numerical linear algebra on
           low-precision hardware: algorithms and
           error analysis},
  url = {https://scicomp.leeds.ac.uk/pasc25-minisymposium-on-matrix-multiplication/}}

@misc{blasnetlib,
  author =       {Netlib},
  title =        {{BLAS} (Basic Linear Algebra Subprograms)},
  url =          {http://www.netlib.org/blas/},
  note =         {Retrieved 31 May 2026},
  year =         2026,
}

@book{LAPACK:1999,
  doi = {10.1137/1.9780898719604},
  url = {https://doi.org/10.1137/1.9780898719604},
  year = {1999},
  month = jan,
  publisher = {Society for Industrial and Applied Mathematics},
  author = {E. Anderson and Z. Bai and C. Bischof and L. S. Blackford and J. Demmel and J. Dongarra and J. Du Croz and A. Greenbaum and S. Hammarling and A. McKenney and D. Sorensen},
  title = {{LAPACK} Users{\textquotesingle} Guide}
}

@article{ahrens2020reproblas,
author={P. Ahrens and J. Demmel and H. D. Nguyen},
title={Efficient reproducible floating point summation and {BLAS}},
journal={ACM Trans. Math. Software},
volume=46,
number=3,
monty=July,
year=2020,
doi={10.1145/3389360},
}

@techreport{Caday2026,
    author = {Caday, P.},
    title = {The {2M} {M}uliplication {A}lgorithm for {C}omplex {M}atrices},
    institution = {NVIDIA},
    year = {2026},
    url = {https://research.nvidia.com/publication/2026-06\_2m-multiplication-algorithm-complex-matrices}
}

@book{Higham2002,
   author = {Higham, N.},
   title = {{Accuracy and Stability of Numerical Algorithms}},
   publisher = {{SIAM}},
   edition = {2nd},
   year = {2002},
   doi = {10.1137/1.9780898718027}}

@techreport{Demmel_Cholesky_89,
  author = {Demmel, J.},
  title = {On floating point errors in {C}holesky},
  institution = {U. of Tennessee, Knoxville},
  type = {{LAPACK Working Note}},
  number = {14},
  year = {1989},
  url = {http://www.netlib.org/lapack/lawnspdf/lawn14.pdf}}

@incollection{Wilkinson_Cholesky_68,
  author = {Wilkinson, J.},
  title = {A priori error analysis of algebraic processes},
  booktitle = {{Proceedings of the International
     Congress of Mathematicians, Moscow 1966}},
  editor = {I. G. Petrovsky},
  pages = {629-640},
  publisher = {Mir Publishers},
  year = {1968}}

@ARTICLE{vandersluis,
      AUTHOR = {Van Der Sluis, A.}, 
      TITLE = {Condition numbers and equilibration of matrices},
      JOURNAL = {Num. Math.},
      YEAR = {1969},
      VOLUME = {14},
      PAGES = {14--23},
      DOI = {10.1007/BF02165096} }

@article{DemmelDumitriuHoltzKleinberg2007,
author = {J. Demmel and I. Dumitriu and O. Holtz
          and R. Kleinberg},
title = {Fast Matrix Multiplication is Stable},
journal = {Numerische Mathematik},
volume = {106},
number = {2},
month = {April},
year = {2007},
doi = {10.1007/s00211-007-0061-6}}

@article{Bini_Lotti1980,
author = {D. Bini and G. Lotti},
title = {Stability of fast algorithms for matrix multiplication},
journal = {Numerische Mathematik},
volume = {36},
pages = {63--72},
year = {1980},
doi = {10.1007/BF01395989}}

@article{Demmel_Underflow_1984,
    author = {J. Demmel},
    title = {Underflow and the reliability of numerical software},
    journal = {{SIAM J. Sci. Stat. Comput.}},
    volume = {5},
    number = {4},
    month = {December},
    year = {1984},
    doi = {10.1137/0905062}
}

@article{BallardBensonDruinskyLipshitzSchwartz2016,
author = {G. Ballard and A. Benson and A. Druinsky and B. Lipshitz and O. Schwartz},
title = {Improving the numerical stability of fast
        matrix multiplication},
journal = {{SIAM J. Mat. Anal. Appl.}},
volume = {37},
number = {4},
pages = {1382--1418},
year = {2016},
doi = {10.1137/15M1032168}}

@misc{Henry_Tang_Heinecke_2019,
   author = {G. Henry and P. T. P. Tang and A. Heinecke},
   title = {Leveraging the bfloat16 artificial intelligence datatype
            for higher-precision computations},
   month = {Apr},
   year = {2019},
   url = {https://arxiv.org/abs/1904.06376}}

@misc{Ootomo_Yokota_2023,
   author = {H. Ootomo and R. Yokota},
   title = {Recovering single precision accuracy from {T}ensor {C}ores
            while surpassing the {FP32} theoretical peak performance},
   month = {Oct},
   year = {2023},
   url = {https://arxiv.org/abs/2203.03341}}

@article{Uchino_Ozaki_Imamura_2024,
  title =        {Performance enhancement of the {Ozaki} Scheme on integer matrix
                  multiplication unit},
  volume =       {39},
  ISSN =         {1741-2846},
  url =          {http://dx.doi.org/10.1177/10943420241313064},
  DOI =          {10.1177/10943420241313064},
  number =       {3},
  journal =      {The International Journal of High Performance Computing
                  Applications},
  publisher =    {SAGE Publications},
  author =       {Uchino, Yuki and Ozaki, Katsuhisa and Imamura, Toshiyuki},
  year =         {2025},
  month =        Jan,
  pages =        {462–476}
}

@article{Ogita_Rump_Oishi_2005,
   author = {Ogita, T. and Rump, S. and Oishi, S.},
   title = {Accurate sum and dot product},
   journal = {SIAM J. Sci. Comput.},
   volume = {26},
   year = {2005},
   doi = {10.1137/030601818}}

@article{Rump_Ogita_Oishi_1_2008,
   author = {Rump, S. and Ogita, T. and Oishi, S.},
   title = {Accurate floating point summation part {I}: faithful rounding},
   journal = {SIAM J. Sci. Comput.},
   volume = {31},
   number = {1},
   pages = {189--224},
   year = {2008},
   doi = {10.1137/050645671}}

@article{Rump_Ogita_Oishi_2_2008,
   author = {Rump, S. and Ogita, T. and Oishi, S.},
   title = {Accurate floating point summation part {II}: sign, k-fold faithful and rounding to nearest},
   journal = {SIAM J. Sci. Comput.},
   volume = {31},
   number = {2},
   pages = {1269--1302},
   year = {2008},
   doi = {10.1137/07068816X}}

@article{Ootomo_Ozaki_Yokota_2024,
author = {Hiroyuki Ootomo and Katsuhisa Ozaki and Rio Yokota},
title ={{DGEMM} on integer matrix multiplication unit},
journal = {The International Journal of High Performance Computing Applications},
volume = {38},
number = {4},
pages = {297-313},
year = {2024},
doi = {10.1177/10943420241239588},
}

@article{Ozaki_Ogita_Oishi_Rump_2012,
  author = {K. Ozaki and T. Ogita and S. Oishi and S. Rump},
  title = {Error-free transformations of matrix multiplication by using
          fast routines of matrix multiplication and its applications},
  journal = {{N}umerical {A}lgorithms},
  volume = {59},
  issue = {1},
  pages = {95--118},
  month = {Jan},
  year = {2012},
  doi = {10.1007/s11075-011-9478-1}}

@article{TPU17,
  author =       {Jouppi, Norman P. and Young, Cliff and Patil, Nishant \textit{et al.}},
  ALTauthor =    {Jouppi, Norman P. and Young, Cliff and Patil, Nishant and
                  Patterson, David and Agrawal, Gaurav and Bajwa, Raminder and
                  Bates, Sarah and Bhatia, Suresh and Boden, Nan and Borchers,
                  Al and Boyle, Rick and Cantin, Pierre-luc and Chao, Clifford
                  and Clark, Chris and Coriell, Jeremy and Daley, Mike and Dau,
                  Matt and Dean, Jeffrey and Gelb, Ben and Ghaemmaghami, Tara
                  Vazir and Gottipati, Rajendra and Gulland, William and
                  Hagmann, Robert and Ho, C. Richard and Hogberg, Doug and Hu,
                  John and Hundt, Robert and Hurt, Dan and Ibarz, Julian and
                  Jaffey, Aaron and Jaworski, Alek and Kaplan, Alexander and
                  Khaitan, Harshit and Killebrew, Daniel and Koch, Andy and
                  Kumar, Naveen and Lacy, Steve and Laudon, James and Law, James
                  and Le, Diemthu and Leary, Chris and Liu, Zhuyuan and Lucke,
                  Kyle and Lundin, Alan and MacKean, Gordon and Maggiore,
                  Adriana and Mahony, Maire and Miller, Kieran and Nagarajan,
                  Rahul and Narayanaswami, Ravi and Ni, Ray and Nix, Kathy and
                  Norrie, Thomas and Omernick, Mark and Penukonda, Narayana and
                  Phelps, Andy and Ross, Jonathan and Ross, Matt and Salek, Amir
                  and Samadiani, Emad and Severn, Chris and Sizikov, Gregory and
                  Snelham, Matthew and Souter, Jed and Steinberg, Dan and Swing,
                  Andy and Tan, Mercedes and Thorson, Gregory and Tian, Bo and
                  Toma, Horia and Tuttle, Erick and Vasudevan, Vijay and Walter,
                  Richard and Wang, Walter and Wilcox, Eric and Yoon, Doe Hyun},
  title =        {In-datacenter performance analysis of a tensor processing
                  unit},
  year =         {2017},
  issue_date =   {May 2017},
  publisher =    {Association for Computing Machinery},
  address =      {New York, NY, USA},
  volume =       {45},
  number =       {2},
  issn =         {0163-5964},
  url =          {https://doi.org/10.1145/3140659.3080246},
  doi =          {10.1145/3140659.3080246},
  journal =      {SIGARCH Comput. Archit. News},
  month =        jun,
  pages =        {1–12},
  numpages =     {12}
}

@article{Rump2009,
  author = {S. Rump},
  title = {Ultimately Fast Accurate Summation},
  journal = {{SIAM. J. Sci. Comp.}},
  volume = {21},
  number = {5},
  year = {2009},
  doi = {10.1137/080738490}}

@article{DemmelHigham1992,
author = {J. Demmel and N. Higham},
title = {Stability of block algorithms with fast level 3 {BLAS}},
journal = {ACM Trans. Math. Software.},
volume = {18},
number = {3},
year = {1992},
doi = {10.1145/131766.131769}
}

@article{10.1145/98267.98290,
author = {Higham, N.},
title = {Exploiting fast matrix multiplication within the level 3 {BLAS}},
year = {1990},
issue_date = {Dec. 1990},
publisher = {Association for Computing Machinery},
address = {New York, NY, USA},
volume = {16},
number = {4},
issn = {0098-3500},
url = {https://doi.org/10.1145/98267.98290},
doi = {10.1145/98267.98290},
journal = {ACM Trans. Math. Softw.},
month = dec,
pages = {352–368},
numpages = {17}
}

@article{Miller1972,
  author =       {W. Miller},
  title =        {Computational complexity and numerical stability},
  journal =      {{SIAM} J. Computing},
  volume =       {4},
  number =       {2},
  month =        {June},
  year =         {1975},
  doi =          {10.1137/0204009}
}

@article{HighamMary2019,
    author = {N. Higham and T. Mary},
    title = {A new approach to probabilistic rounding error analysis},
    journal = {{SIAM J. Sci. Comput.}},
    volume = {41},
    number = {5},
    year = {2019},
    doi = {10.1137/18M1226312}}

@ARTICLE{Strassen69,
      AUTHOR = {Strassen, V.},
      TITLE = {{Gaussian} elimination is not optimal},
      JOURNAL = {Numerische Mathematik},
      VOLUME = {13},
      NUMBER = {4},
      PAGES = {354--356},
      YEAR = {1969},
      DOI = {10.1007/BF02165411}
}

@inproceedings{ExceptionHandling-1,
  author =       {James Demmel and Jack J. Dongarra and Mark Gates and Greg
                  Henry and Julien Langou and Xiaoye S. Li and Piotr Luszczek
                  and Weslley Pereira and E. Jason Riedy and Cindy
                  Rubio{-}Gonz{\'{a}}lez},
  title =        {Proposed Consistent Exception Handling for the {BLAS} and
                  {LAPACK}},
  booktitle =    {Proc. Intern. Workshop on Software Correctness for HPC
                  Applications (CORRECTNESS'22)},
  month =        {Nov},
  year =         {2022},
  note =         {Extended version},
  url =          {https://arxiv.org/abs/2207.09281},
}

@article{xblas:toms,
     AUTHOR = {X. S. Li and J. W. Demmel and D. H. Bailey and G. Henry
              and Y. Hida and J. Iskandar and W. Kahan and S. Y. Kang and
              A. Kapur and M. C. Martin and B. J. Thompson and T. Tung
              and D. J. Yoo},
     title = {Design, implementation and testing of
              extended and mixed precision {BLAS}},
     journal = {{ACM} Trans. Math. Softw.},
     volume = 28,
     number = 2,
     pages = {152-205},
     YEAR = {2002},
     DOI = {10.1145/567806.567808}
}

@article{DemmelHida2003,
  author = {J. Demmel and Y. Hida},
  title = {Accurate and efficient floating point summation},
  journal = {{SIAM J. Sci. Comp.}},
  volume = {25},
  number = {4},
  year = {2003},
  doi = {10.1137/S1064827502407627}}

@INPROCEEDINGS{EXCVATE2025,
  author =       {Vanover, J. and Li, X. S. and Rubio-Gonzalez, C. and Demmel,
                  J.},
  booktitle =    {Proceedings 32nd {IEEE} Symposium on Computer Arithmetic.
                  {ARITH}'25, 2025},
  title =        {{EXCVATE}: Spoofing exceptions and solving constraints to test
                  exception handling in numerical libraries},
  year =         {2025},
  doi =          {10.1109/ARITH64983.2025.00026},
}

@article{Paranoia,
author={Richard Karpinski},
title={Paranoia: A floating-point benchmark},
journal={Byte Magazine},
volume=10,
number=2,
month={February},
year=1985,
pages={223-235},
url={https://people.math.sc.edu/Burkardt/c_src/paranoia/paranoia.html}
}

@Misc{paranoia-code,
  author =       {W. M. Kahan},
  title =        {A  {PARANOID}  {PROGRAM}  to  {DIAGNOSE}  {FLOATING-POINT}  {ARITHMETIC} in  {BASIC}  v. {A1.10}},
  howpublished = {Netlib},
  month =        {May},
  year =         {1982},
  url = {https://www.netlib.org/paranoia/},
  note =         {Pascal port by B. A. Wichmann. C port by David M. Gay and Thos Sumner. Bug fixes from David Hough.},
  OPTannote =    {}
}

@Article{K_gstr_m_1998,
  author       = {Kågström, Bo and Ling, Per and van Loan, Charles},
  title        = {{GEMM}-based level 3 {BLAS}: {H}igh-performance model
                  implementations and performance evaluation benchmark},
  journal      = {ACM Transactions on Mathematical Software},
  year         = 1998,
  volume       = 24,
  number       = 3,
  month        = sep,
  pages        = {268–302},
  issn         = {1557-7295},
  doi          = {10.1145/292395.292412},
  url          = {http://dx.doi.org/10.1145/292395.292412},
  publisher    = {Association for Computing Machinery (ACM)}
}

@Article{Zhang_2025,
  author       = {Zhang, Yu and Lu, Lu and Yang, Zhanyu and Liang, Zhihong and
                  Suo, Siliang},
  title        = {{LE-GEMM}: A lightweight emulation-based {GEMM} with precision
                  refinement on {GPU}},
  journal      = {Journal of Systems Architecture},
  year         = 2025,
  volume       = 160,
  month        = mar,
  pages        = 103336,
  issn         = {1383-7621},
  doi          = {10.1016/j.sysarc.2025.103336},
  url          = {http://dx.doi.org/10.1016/j.sysarc.2025.103336},
  publisher    = {Elsevier BV}
}

@inproceedings{uchino2025high,
  title={High-performance and power-efficient emulation of matrix multiplication using {INT8} matrix engines},
  author={Uchino, Yuki and Ozaki, Katsuhisa and Imamura, Toshiyuki},
  booktitle={Proceedings of the SC'25 Workshops of the International Conference for High Performance Computing, Networking, Storage and Analysis},
  pages={1824--1831},
  year={2025},
  doi={10.1145/3731599.3767539}
}

@book{golub-van-loan,
  author =       {Golub, Gene H. and Van Loan, Charles F.},
  title =        {Matrix Computations - 4th Edition},
  publisher =    {Johns Hopkins University Press},
  year =         {2013},
  doi =          {10.1137/1.9781421407944},
  address =      {Philadelphia, PA},
  edition =      {},
  URL =          {https://epubs.siam.org/doi/abs/10.1137/1.9781421407944},
  eprint =       {https://epubs.siam.org/doi/pdf/10.1137/1.9781421407944}
}

@misc{Hopper,
  author =       {NVIDIA},
  title =        {{NVIDIA H100 Tensor Core GPU Architecture}},
  url =          {https://resources.nvidia.com/en-us-hopper-architecture/nvidia-h100-tensor-c},
  note =         {Retrieved 31 May 2026},
}

@misc{Rubin,
  author =       {Kyle Aubrey},
  title =        {Inside the {NVIDIA} {Rubin} Platform: Six New Chips, One {AI}
                  Supercomputer},
  year =         2026,
  url =          {https://developer.nvidia.com/blog/inside-the-nvidia-rubin-platform-six-new-chips-one-ai-supercomputer/},
  note =         {Retrieved 31 May 2026},
}

@inproceedings{Gall2014,
author = {Le Gall, Fran\c{c}ois},
title = {{Powers of tensors and fast matrix multiplication}},
year = {2014},
isbn = {9781450325011},
publisher = {Association for Computing Machinery},
address = {New York, NY, USA},
doi = {10.1145/2608628.2608664},
booktitle = {Proceedings of the 39th International Symposium on Symbolic and Algebraic Computation},
pages = {296–303},
numpages = {8},
location = {Kobe, Japan},
series = {ISSAC '14}
}

@article{Coppersmith-Winograd1990,
author = {Coppersmith, Don and Winograd, Shmuel},
title = {{Matrix multiplication via arithmetic progressions}},
year = {1990},
issue_date = {March 1990},
publisher = {Academic Press, Inc.},
address = {USA},
volume = {9},
number = {3},
issn = {0747-7171},
doi = {10.1016/S0747-7171(08)80013-2},
journal = {J. Symb. Comput.},
month = mar,
pages = {251–280},
numpages = {30}
}

@INPROCEEDINGS{chatBLAS2024,
  author={Valero-Lara, Pedro and Godoy, William F. and Teranishi, Keita and Balaprakash, Prasanna and Vetter, Jeffrey S.},
  booktitle={SC24-W: Workshops of the International Conference for High Performance Computing, Networking, Storage and Analysis}, 
  title={{ChatBLAS}: The first {AI}-generated and portable {BLAS} library}, 
  year={2024},
  volume={},
  number={},
  pages={19-24},
  doi={10.1109/SCW63240.2024.00010}}

@article{dumitrescu1998accuracy,
  title =        {Improving and estimating the accuracy of {Strassen's}
                  algorithm},
  author =       {Dumitrescu, Bogdan},
  journal =      {Numerische Mathematik},
  volume =       {79},
  number =       {4},
  pages =        {485--499},
  year =         {1998},
  publisher =    {Springer},
  doi =          {10.1007/s002110050348}
}

@inproceedings{almanduanwilliamsxuxuzhou25,
   author = {Alman, J. and Duan, R. and
       Vassilevska Williams, V. and Xu, Y. and Xu, Z. and Zhou, R.},
   title = {More asymmetry yields faster matrix multiplication},
   booktitle = {Proc. 2025 ACM-SIAM Symp. on
      Discrete Algorithms (SODA)},
   year = {2025},
   doi = {10.1137/1.9781611978322.63}}

@misc{FalconGEMM,
author={Zhu, H. and Cao, J. and Shao, J. and Feng, S. and Qiu, Q. and Chen, P. and Zhang, X. and Zhou, Y. and Yiu, M. L. and Ji, G. and Deng, M. and Zhu, W. and Meng, J.},
title={{FalconGEMM}: Surpassing hardware peaks with lower-complexity matrix multiplication},
month={May},
year={2026},
url={https://arxiv.org/pdf/2605.06057},
}

@article{AlphaTensor,
   author = {Fawzi, A. and Balog, M. and Huang, A. and Hubert, T. and Romera-Paredes, B. and Barekatain, M. and Novikov, A. and Ruiz, F. J. R. and Schrittwieser, J. and Swirszcz et. al.},
   title = {Discovering faster matrix multiplication algoroithms with reinforcement learning},
   journal = {Nature},
   volume = {610},
   number = {7930},
   year = {2022},
   doi = {10.1038/s41586-022-05172-4}}

@article{10.1145/356004.356006,
  author =       {Schreiber, Robert},
  title =        {A New Implementation of Sparse {Gaussian} Elimination},
  year =         {1982},
  issue_date =   {Sept. 1982},
  publisher =    {Association for Computing Machinery},
  address =      {New York, NY, USA},
  volume =       {8},
  number =       {3},
  issn =         {0098-3500},
  url =          {https://doi.org/10.1145/356004.356006},
  doi =          {10.1145/356004.356006},
  journal =      {ACM Trans. Math. Softw.},
  month =        sep,
  pages =        {256–276},
  numpages =     {21}
}

@article{doi:10.1137/S0895479896297744,
  author =       {Toledo, Sivan},
  title =        {Locality of Reference in {LU} Decomposition with Partial
                  Pivoting},
  journal =      {{SIAM} Journal on Matrix Analysis and Applications},
  volume =       {18},
  number =       {4},
  pages =        {1065-1081},
  year =         {1997},
  doi =          {10.1137/S0895479896297744},
  URL =          { https://doi.org/10.1137/S0895479896297744 },
}

\end{document}